\documentclass[11pt]{amsart}         
\usepackage{amscd}                   
\usepackage{amssymb}

\newtheorem{thm}{Theorem}[section]
\newtheorem{lem}[thm]{Lemma}

\newtheorem{cor}[thm]{Corollary}
\newtheorem{claim}[thm]{Claim}
\renewcommand{\thestep}{}
\newtheorem{prop}[thm]{Proposition}

\theoremstyle{definition}

\renewcommand{\thecase}{}

\theoremstyle{remark}

\makeatletter
\def\alphenumi{
  \def\theenumi{\alph{enumi}}
  \def\p@enumi{\theenumi}
  \def\labelenumi{(\@alph\c@enumi)}}
\makeatother

\makeatletter
\def\thecase{\@arabic\c@case}
\makeatother

\numberwithin{equation}{section}

\makeatletter
\def\thestep{\@arabic\c@step}
\makeatother

\newenvironment{pf}{\begin{proof}[\proofname]}{\end{proof}}
\newenvironment{pf*}[1]{\begin{proof}[#1]}{\end{proof}}

\newcommand\barB{{\bar{B}}}

\newcommand\barM{{\bar{M}}}

\newcommand\EE{\mathbb{E}}

\newcommand\RR{\mathbb{R}}

\newcommand\bPhi{{\boldsymbol{\Phi}}}

\newcommand\bPsi{{\boldsymbol{\Psi}}}

\newcommand\bB{{\mathbf{B}}}

\newcommand\bE{{\mathbf{E}}}

\newcommand\bF{{\mathbf{F}}}

\newcommand\bK{{\mathbf{K}}}

\newcommand\bS{{\mathbf{S}}}

\newcommand\bx{{\mathbf{x}}}

\newcommand{\cov}{\nabla}

\newcommand{\rd}{\partial}

\newcommand\half{{\textstyle{\frac{1}{2}}}}

\newcommand\quarter{{\textstyle{\frac{1}{4}}}}

\newcommand\fB{{\mathfrak{B}}}

\newcommand\fg{{\mathfrak{g}}}

\newcommand\al{\alpha}
\newcommand\be{\beta}
\newcommand\De{\Delta}
\newcommand\de{\delta}
\newcommand\eps{\varepsilon}
\newcommand\ga{\gamma}
\newcommand\Ga{\Gamma}
\newcommand\la{\lambda}
\newcommand\La{\Lambda}

\newcommand\om{\omega}
\newcommand\Om{\Omega}

\newcommand\gl{{\mathfrak{g}\mathfrak{l}}}

\newcommand\so{{\mathfrak{s}\mathfrak{o}}}
\newcommand\stab{{\mathfrak{s}\mathfrak{t}\mathfrak{a}\mathfrak{b}}}
\newcommand\su{{\mathfrak{s}\mathfrak{u}}}

\newcommand\GL{\operatorname{GL}}

\newcommand\PU{\operatorname{PU}}

\newcommand\SO{\operatorname{SO}}

\newcommand\U{\operatorname{U}}

\newcommand\less{\setminus}

\newcommand{\8}{\infty}

\newcommand\Ad{{\operatorname{Ad}}}

\newcommand\Aut{\operatorname{Aut}}

\newcommand\Center{\operatorname{Center}}

\newcommand\dist{\operatorname{dist}}
\newcommand\divg{\operatorname{div}}

\newcommand\End{\operatorname{End}}
\newcommand\Exp{\operatorname{Exp}}

\newcommand\Hom{\operatorname{Hom}}

\newcommand\Imag{\operatorname{Im}}

\newcommand\Ker{\operatorname{Ker}}

\newcommand\Ric{\operatorname{Ric}}

\newcommand\Rm{\operatorname{Rm}}

\newcommand\Stab{\operatorname{Stab}}

\newcommand\supp{\operatorname{supp}}

\newcommand\Sym{\operatorname{Sym}}

\newcommand\tr{\operatorname{tr}}

\newcommand\id{{\mathrm{id}}}

\newcommand\sA{{\mathcal{A}}}
\newcommand\sB{{\mathcal{B}}}

\newcommand\sG{{\mathcal{G}}}

\newcommand\sL{{\mathcal{L}}}

\newcommand\tg{{\tilde g}}

\newcommand\tS{{\widetilde S}}
\newcommand\tT{{\widetilde T}}

\begin{document}
\title[Critical-exponent Sobolev Norms and the Slice Theorem]
{Critical-exponent Sobolev Norms and
the Slice Theorem for the Quotient Space of Connections}
\author[Paul M. N. Feehan]{Paul M. N. Feehan}
\address{Department of Mathematics\\
Harvard University\\
One Oxford Street\\
Cambridge, MA 02138}
\email{feehan@math.harvard.edu}
\date{Completed October 2, 1996. dg-ga/9711004. 
This version: February 1, 1999.}
\thanks{The author was supported in part by an NSF Mathematical 
Sciences Postdoctoral Fellowship under grant DMS 9306061}
\maketitle


\section{Introduction}
The use of certain `critical-exponent' Sobolev norms is an important
feature of methods employed by Taubes to solve the anti-self-dual and related
non-linear elliptic partial differential equations
\cite{TauFrame,TauStable,TauConf}. Indeed, the estimates one can obtain
using these critical-exponent norms appear to be the best possible when one
needs to bound the norm of a Green's operator for a Laplacian, depending on
a connection varying in a non-compact family, in terms of minimal data such
as the first positive eigenvalue of the Laplacian or the $L^2$ norm of the
curvature of the connection. Despite their utility, particularly in
applications where an optimal analysis is required for gluing or
degeneration problems (for example, when considering degenerating families
of anti-self-dual connections or stable, holomorphic vector bundles --- see
Section \ref{subsec:Applic} below), these methods are not widely known.
Following Taubes \cite{TauPath,TauFrame,TauStable,TauConf} we describe a
collection of critical-exponent Sobolev norms and general Green's operator
estimates depending only on first positive eigenvalues or the $L^2$ norm of
the connection's curvature. These estimates are especially useful both for
the construction of gluing maps, in the case of either anti-self-dual
connections \cite{TauStable} or, more recently, in the case of $\PU(2)$
monopoles \cite{FL3,FL4}, and for analyzing their asymptotic behavior with
respect to Uhlenbeck limits of the underlying gluing data.  We apply them
here to prove an optimal slice theorem for the quotient space of
connections. The result is `optimal' in the sense that if a point $[A]$ in
the quotient space is known to be just $L^2_1$-close enough to a reference
point $[A_0]$ (see below for the precise statement), then $A$ can be placed
in Coulomb gauge relative to $A_0$, with all constants depending at most on
the first positive eigenvalue of the covariant Laplacian defined by $A_0$
and the $L^2$ norm of the curvature of $A_0$. Such slice theorems are
particularly advantageous when analyzing gluing maps and their
differentials in situations (such as those of
\cite{FeehanLeness,MorganOzsvath} and \cite{FL3}) where the underlying
gluing data is allowed to `bubble'. In this paper we shall for simplicity
only consider connections over four-dimensional manifolds, but the methods
and results can adapted to the case of manifolds of arbitrary dimension, as
in \cite{UhlLp}, to prove slice theorems applicable to cases where the
reference connection is allowed to degenerate.

\subsection{Critical-exponent Sobolev norms and the slice theorem}
Suppose that $X$ is a closed, Riemannian 
four-manifold, that $G$ is a compact Lie
group, and that $\sB_E^{k,p}=\sA_E^{k,p}/\sG_E^{k+1,p}$ is the quotient
space of $L^p_k$ connections on a $G$ bundle $E$ modulo the Banach Lie
group of $L^p_{k+1}$ gauge transformations. Here, the integer $k\ge 1$ and
the Sobolev exponent $1<p<\8$ obey the constraint $(k+1)p>4$, so
$L^p_{k+1}(X)\subset C^0(X)$ and gauge transformations in $\sG_E^{k+1,p}$ are
continuous. When $(k+1)p=4$ we have the `borderline', `critical', or
`limiting case' of the Sobolev embedding theorem: $L^p_{k+1}(X)\subset L^q(X)$
for all $q<\8$ but not $q=\8$.

A connection $A\in \sA_E^{k,p}$ is in {\em Coulomb gauge}
relative to a reference connection $A_0$ if $d_{A_0}^*(A-A_0)=0$ and it is
a standard result that $\bS_{A_0} = A_0+\Ker
d_{A_0}^*\subset
\sA_E^{k,p}$ provides a {\em slice\/} for the action of the gauge group
$\sG_E^{k+1,p}$  \cite{AHS,DK,FU,Lawson,MitterViallet,
MorganGTNotes,Parker,Singer}. 
More exactly, if $\bB_{A_0}^{k,p}(\eps)$ is the $L^p_k$ ball in
$\bS_{A_0}$ with center $A_0$ and $L^p_{k,A_0}$-radius $\eps$ and
$\Stab_{A_0}\subset \sG_E^{k+1,p}$ is the stabilizer of $A_0$, then
the projection $\pi:\bB_{A_0}^{k,p}(\eps)/\Stab_{A_0}\to \sB_E^{k,p}$ is a
homeomorphism onto its image and thus contains a small enough $L^p_k$ ball
$$
B_{[A_0]}^{k,p}(\eta) = \{[A]\in\sB_E^k:\dist_{L^p_{k,A_0}}([A],[A_0]) <
\eta\}, 
$$
where the gauge-invariant distance function on the quotient is defined by
$$
\dist_{L^p_{k,A_0}}([A],[A_0])
= \inf_{u\in\sG_E^{k+1,p}}\|u(A)-A_0\|_{L^p_{k+1,A_0}}.    
$$
One unsatisfactory aspect of the standard slice theorem concerns the
dependence of the constants $\eps([A_0],k,p)$ and $\eta([A_0],k,p)$ above
on the orbit $[A_0]$ --- in particular on the curvature $F_{A_0}$ --- when
$k$ and $p$ are large enough that gauge transformations in $\sG_E^{k+1,p}$
are continuous. Even in the minimal cases, $k=1$ and $p>2$ or $k=2$ and
$p=2$, the constants $\eps,\eta$ depend unfavorably on $[A_0]$ when the
curvature $F_{A_0}$ `bubbles' and $[A_0]$ approaches an ideal point in the
Uhlenbeck compactification $\barM_E$ of the moduli space $M_E$ of
anti-self-dual connections (that is, a point in $\barM_E\less M_E$). This makes
it difficult to analyze the asymptotic behavior of Taubes' gluing maps
\cite{TauSelfDual,TauIndef,TauFrame,TauStable} and their differentials on
neighborhoods of ideal points in $\barM_E$, since the balls
$\bB_{A_0}^{k,p}(\eps)$ and $B_{[A_0]}^{k,p}(\eta)$ tend to shrink as
$[A_0]$ approaches a point in $\barM_E\less M_E$.  For example, if the
connection $A_0$ is anti-self-dual, then its {\em energy\/} is bounded by a
constant depending only on the topology of $E$ via the Chern-Weil identity
$$
-{\frac{1}{4\pi^2}}\int_X\tr(F_{A_0}\wedge F_{A_0}) = p_1(\fg_E), 
$$
whereas $\|F_{A_0}\|_{L^p}$ (with $p>2$) or $\|F_{A_0}\|_{L^2_{1,A_0}}$
tends to infinity as the curvature of $A_0$ becomes concentrated and
$[A_0]$ approaches the Uhlenbeck boundary. 

Our main purpose in this article is to prove a global analogue, 
Theorem \ref{thm:GaugeFixing}, of Uhlenbeck's
local Coulomb gauge-fixing theorem \cite[Theorems 1.3 \& 2.1]{UhlLp} and a
corresponding slice theorem, Theorem \ref{thm:Slice}, where the radii of
the coordinate balls on the quotient $\sB_E^{k,p}$ depend only on
$\|F_{A_0}\|_{L^2}$ and the {\em least positive eigenvalue\/} $\nu_0[A_0]$
of the Laplacian $d_{A_0}^*d_{A_0}$ on $\Om^0(\fg_E)$. 
The key difficulty in establishing Theorem
\ref{thm:GaugeFixing} is to ensure that the constants depend at most on
$\|F_{A_0}\|_{L^2}$ and $\nu_0[A_0]$: To guarantee this minimal dependence,
we employ critical-exponent Sobolev norms (defined below) to circumvent the
fact that when $(k+1)p=4$ the standard Sobolev embedding and multiplication
theorems fall just short of what one needs to give the quotient
$\sB_E^{k,p}=\sA_E^{k,p}/\sG_E^{k+1,p}$ a manifold structure (see Section
\ref{sec:SharpSobolev}). Such norms were introduced by Taubes for
related purposes in \cite{TauFrame}.

\subsection{Statement of results}
\label{subsec:State}
For clarity, we now fix $p=2$ and $k\ge 2$ and 
define the following distance functions on 
the quotient space $\sB_E^k=\sA_E^k/\sG_E^{k+1}$ of $L^2_k$
connections modulo $L^2_{k+1}$ gauge transformations,
\begin{align*}
\dist_{\sL^{\sharp,2}_{1,A_0}}([A],[A_0])
&= \inf_{u\in\sG_E^{k+1}}\left(\|u(A)-A_0\|_{L^{2\sharp,4}}
+ \|d_{A_0}^*(u(A)-A_0)\|_{L^{\sharp,2}}\right), \\
\dist_{L^{\sharp,2}_{1,A_0}}([A],[A_0])
&= \inf_{u\in\sG_E^{k+1}}\left(\|u(A)-A_0\|_{L^2_{1,A_0}}
+ \|d_{A_0}^*(u(A)-A_0)\|_{L^{\sharp,2}}\right),
\end{align*}
where the Sobolev norms are defined by:
\begin{align*}
\|a\|_{L^{\sharp}(X)} &= \sup_{x\in X}\|\dist^{-2}(x,\cdot)|a|\|_{L^1(X)},\\
\|a\|_{L^{\sharp,2}(X)} &= \|a\|_{L^2(X)} + \|a\|_{L^{\sharp}(X)}, \\
\|a\|_{L^{2\sharp}(X)} &= \sup_{x\in X}\|\dist^{-1}(x,\cdot)|a|\|_{L^2(X)},\\
\|a\|_{L^{2\sharp,4}(X)} &= \|a\|_{L^4(X)} + \|a\|_{L^{2\sharp}(X)}, \\
\|a\|_{L^2_{1,A_0}(X)} &= 
\left(\|a\|_{L^2(X)}^2 + \|\cov_{A_0}a\|_{L^2(X)}^2\right)^{1/2},
\end{align*}
for any $a\in\Om^1(\fg_E)$; here, $\dist(x,y)$ denotes the geodesic
distance between points $x,y\in X$.
Like the $L^4$ norm, the $L^{2\sharp}$ norm on
one-forms is scale-invariant. Our first result is the following global
analogue of Uhlenbeck's theorem and complements results of Taubes in
\cite[\S 6]{TauFrame}:

\begin{thm}\label{thm:GaugeFixing}
Let $X$ be a closed, smooth four-manifold with metric $g$
and let $G$ be a compact Lie
group. Then there are positive constants $c,z$ with the following
significance.  Let $E$ be a $G$ bundle over $X$ and suppose that $k\ge 2$
is an integer.  Given a point $[A_0]$ in $\sB_E^k$, let $\nu_0[A_0]$ be the
least positive eigenvalue of the Laplacian $d_{A_0}^*d_{A_0}$
on $\Om^0(\fg_E)$ and set $K_0 =
(1+\nu_0[A_0]^{-1})(1+\|F_{A_0}\|_{L^2})$. Let $\eps_1$ be a 
constant satisfying $0<\eps_1\le zK_0^{-2}(1+\nu_0[A_0]^{-1/2})^{-1}$. 
Then the following hold:
\begin{enumerate}
\item For any $[A]\in\sB^k_E$ with
$\dist_{\sL^{\sharp,2}_{1,A_0}}([A],[A_0]) < \eps_1$,
there is a gauge transformation $u\in \sG_E^{k+1}$, unique up to an element
of the stabilizer $\Stab_{A_0}\subset\sG_E^{k+1}$, such that
\begin{enumerate}
\item $d_{A_0}^*(u(A)-A_0) = 0$, 
\item $\|u(A)-A_0\|_{L^{2\sharp,4}} \le
cK_0\dist_{\sL^{\sharp,2}_{1,A_0}}([A],[A_0])$. 
\end{enumerate}
\item For any $[A]\in\sB^k_E$ with
$\dist_{L^{\sharp,2}_{1,A_0}}([A],[A_0]) < \eps_1$,
there is a gauge transformation $u\in \sG_E^{k+1}$, unique up to an element
of the stabilizer $\Stab_{A_0}\subset\sG_E^{k+1}$, such that
\begin{enumerate}
\item $d_{A_0}^*(u(A)-A_0) = 0$, 
\item $\|u(A)-A_0\|_{L^{2\sharp,4}} \le
cK_0\dist_{\sL^{\sharp,2}_{1,A_0}}([A],[A_0])$, 
\item $\|u(A)-A_0\|_{L^2_{1,A_0}} \le
cK_0\dist_{L^{\sharp,2}_{1,A_0}}([A],[A_0])$. 
\end{enumerate}
\end{enumerate}
\end{thm}

Theorem \ref{thm:GaugeFixing} is especially well-adapted to establishing
the existence of transformations to Coulomb gauge when the point
$[A]\in\sB_E^k$ has the form $A=A_0+d_{A_0}^{+,*}v$, with $A_0$
approximately anti-self-dual (so $F_{A_0}^+$ is small)
and $v\in L^2_{k+1}(\La^+\otimes\fg_E)$, since
$d_{A_0}^*d_{A_0}^{+,*}v = (F_{A_0}^+)^*v$. Points of this form in the
moduli space of anti-self-dual connections, $M_E$, are constructed by
Taubes' gluing maps \cite{TauSelfDual,TauIndef,TauStable}. 

The distance function $\dist_{L^{\sharp,2}_{1,A_0}}([A],[A_0])$ is bounded
by scale invariant norms,
$$
\|a\|_{L^4(X,g)} + \|\cov_A^ga\|_{L^2(X,g)} 
+ \sup_{x\in X}\|\dist_g^{-2}(x,\cdot)
|d_{A_0}^{*_g}a|\|_{L^1(X,g)},\qquad a\in \Om^1(\fg_E),
$$
since the $L^{4/\ell}$ norm on $\otimes^\ell(T^*X)$ is
conformally invariant, while the third term is invariant under constant
rescalings $g\mapsto\tg=\la^{-2}g$ of the metric, as $d_{A_0}^{*_\tg}a =
\la^2d_{A_0}^{*_g}a$, $\dist_\tg^{-2}(x,y) = \la^2\dist_g^{-2}(x,y)$ and
$dV_\tg = \la^{-4}dV_g$. Similarly for
$\dist_{\sL^{\sharp,2}_{1,A_0}}([A],[A_0])$. 

In Theorem 2.1 of \cite{UhlLp} the $L^2$ norm of the curvature $F_A$ of a
local connection matrix $A$ over the unit ball in $\RR^4$ provides a natural
(gauge-invariant) measure of the distance from $[A]$ to $[\Ga]$, where
$\Ga$ is the product connection. Uhlenbeck's theorem guarantees the existence 
of an $L^p_{k+1}$ gauge transformation $u$ taking
an $L^p_k$ connection $A$ on the product bundle over the unit four-ball,
with product connection $\Ga$, to a connection
$u(A)$ satisfying
$d_\Ga^*(u(A)-\Ga)=0$ and $\|u(A)-\Ga\|_{L^2_1}\le
c\|F_A\|_{L^2}$; one only requires that $\|F_A\|_{L^2}$ be smaller than
a universal constant. 

We next have the following refinement of the standard slice theorem for the
quotient space $\sB_E^k$. The observation that an $L^4$-ball in $\Ker
d_{A_0}^*$ provides a slice for $\sG_E^{k+1}$ was pointed out to us Mrowka;
that slightly smaller $L^{2\sharp,4}$ and $L^2_{1,A_0}$ balls provide
slices follows from the second of our two proofs of Theorem
\ref{thm:GaugeFixing} in Section \ref{sec:SliceII}.
For any $\eps>0$, define
\begin{align*}
B_{[A_0]}^{1,\sharp,2}(\eps) &=
\{[A]\in\sB_E^k:\dist_{L^{\sharp,2}_{1,A_0}}([A],[A_0]) < \eps\}
\subset \sB_E^k, \\
B_{[A_0]}^{1,*,2}(\eps) &=
\{[A]\in\sB_E^k:\dist_{\sL^{\sharp,2}_{1,A_0}}([A],[A_0]) < \eps\} 
\subset \sB_E^k,\\ 
\bB^4_{A_0}(\eps) &= \{A\in \sA_E^k: d_{A_0}^*(A-A_0)=0\text{ and
}\|A-A_0\|_{L^4(X)}<\eps\}
\subset \bS_{A_0},
\end{align*}
where $\bS_{A_0}=\{A_0\}+\Ker(d_{A_0}^*|_{L^2_k})\subset\sA_E^k$ is the slice
through $A_0$. 

\begin{thm}\label{thm:Slice}
Let $X$ be a closed, smooth four-manifold with metric $g$ and let $G$ be a
compact Lie group. Then there are positive constants $c_1,c_2,z$ with the
following significance.  Let $E$ be a $G$ bundle over $X$, let $k\ge 2$ be
an integer, and suppose that $[A_0]\in\sB_E^k$.  Then the following hold:
\begin{enumerate}
\item For any constant $\eps_0$ satisfying 
$0<\eps_0<z(1+\nu_0[A_0]^{-1/2})^{-1}$,
the projection $\pi:\bB^4_{A_0}(\eps_0)/\Stab_{A_0}\to \sB_E^k$ given by
$A\mapsto [A]$ is a homeomorphism onto an open neighborhood of $[A_0]\in
\sB_E^k$ and a diffeomorphism on the open subset where 
$\Stab_{A_0}/\Center(G)$ acts freely;
\item For any constant $\eps_1$ satisfying
$0<\eps_1\le zK_0^{-2}(1+\nu_0[A_0]^{-1/2})^{-1}$ we have the following
inclusions of open neighborhoods in $\sB_E^k$:
$$
B_{[A_0]}^{1,\sharp,2}(\eps_1)
\subset B_{[A_0]}^{1,*,2}(c_1\eps_1)
\subset \pi(\bB^4_{A_0}(c_2K_0\eps_1)).
$$
\end{enumerate}
\end{thm}

Let $\sA_E^{*,k}\subset\sA_E^k$ be the subspace of connections $A$ with
minimal stabilizer $\Stab_A=\Center(G)$ and let $\sB_E^{*,k} =
\sA_E^{*,k}/\sG_E^{k+1}$. It is well-known that the quotient space
$\sB_E^k$ is Hausdorff, that
the subspace $\sB_E^{*,k}\subset \sB_E^k$ is an
open, $C^\8$ Banach manifold, and that the projection $\pi:\sA_E^{*,k} \to
\sB_E^{*,k}$ is a $C^\8$ principal $\sG_E^{k+1}/\Center(G)$ bundle. See
Proposition \ref{prop:Slice} for detailed statements.

That sharper versions of the standard slice theorem (as in
\cite{DK,FU,Lawson}, for example) would hold is suggested by related
results of Taubes, namely \cite[Lemma A.1]{TauPath} and \cite[Lemma
6.5]{TauFrame}: for example, they show that if $u$ is an $L^2_2$ gauge
transformation intertwining $L^2_1$ connections $A_i$, $i=1,2$, obeying a
slice condition $d_{A_0}^*(A_i-A_0)=0$ defined by an $L^2_1$ connection
$A_0$, then $u$ is necessarily in $C^0$. Moreover, transition functions
relating neighborhoods of the origin in
$\Ker (d_{A_0}^*|_{L^2_1})$ and $\Ker (d_{A_0+a}^*|_{L^2_1})$,
where $a$ is $L^2_{1,A_0}$-small, are constructed in \cite[Lemma
6.5]{TauFrame}; the 
constants depend only on $\|F_{A_0}\|_{L^2}$ and $\nu_0[A_0]$. (See
\cite[\S 6]{TauFrame} for detailed statements and related results.)  The
proof of Theorem \ref{thm:GaugeFixing} makes use of methods developed in
\cite{TauPath,TauFrame,TauStable}.

\subsection{Outline of the proofs}
\label{subsec:Outline}
Assertion (1) of Theorem \ref{thm:Slice} is proved in Section
\ref{sec:SliceI}. The proof that the projection map
$\pi:\bB_{A_0}^4(\eps_0)\to\sB_E^k$ is a local diffeomorphism away from
connections with non-minimal stabilizer essentially follows Uhlenbeck's
verification of `openness' in her proof of Theorem 2.1 in \cite{UhlLp} 
via the method of continuity (see Lemma
\ref{lem:LocalDiffeo}).  The proof that the $L^4$ ball 
$\bB_{A_0}^4(\eps_0)$ injects into the quotient (see Lemma
\ref{lem:injective}) was suggested to us by Mrowka.  
The remainder of our article is 
taken up with the 
proof of Theorem \ref{thm:GaugeFixing} and hence Assertion (2) of Theorem
\ref{thm:Slice}.

In Section \ref{sec:SharpSobolev} we introduce the family of 
critical-exponent Sobolev norms, $L^{\sharp,2}_{k,A_0}$, $k=0,1,2$, 
used to complete the proof of Theorem
\ref{thm:GaugeFixing} and in Section \ref{sec:Green} we describe the crucial
embedding theorems enjoyed by those Sobolev spaces, as well as estimates
for the Green's operator of the Laplacian $d_{A_0}^*d_{A_0}$.  In particular,
$L^p_{k,A_0}\subset L^{\sharp,2}_{k,A_0}$, for every $p>2$ while, in the
other direction, $L^{\sharp,2}_{2,A_0}\subset C^0$. The latter embedding is
the key motivation for the definition of these norms and it greatly
facilitates the derivation of Green's operator estimates, in a wide number
of applications in gauge theory \cite{TauFrame,TauStable}, with minimal
dependence on the curvature of the connection $A_0$.  The main
ideas and embedding results in Sections \ref{sec:SharpSobolev} and
\ref{sec:Green} are due to Taubes
\cite{TauPath,TauFrame,TauStable,TauConf}, so these sections are essentially
expository. An earlier exposition from a somewhat different perspective,
due to Donaldson, of Taubes' methods and some applications appears in
\cite{DonApprox}.  The estimates of Section \ref{sec:Green} are stated only
in the four-dimensional case. While we might expect all of them to hold, in
some form, for higher dimensions we confine our attention to dimension four
as our intended applications are primarily concerned with smooth
four-manifold topology. In essence, the critical-exponent norms make a
virtue out of necessity of the familiar fact that while the Green's
operator of the Laplacian $d^*d$ on $C^\8(X)$ maps $L^p(X)$ into
$L^{2p/(2-p)}(X)$ for $1<p<2$, it does not map $L^2(X)$ into $L^\8(X)$
\cite[Chapter V]{Stein}. We recall that an
Orlicz space $L_\varphi$ can be used to provide the `best target space' for
an embedding of $L^2_2(X)$ \cite[Chapter 8]{Adams}. Here, we may instead
view $L^{\sharp,2}_2(X)$ as providing the `best domain space' for an
embedding into $L^\8(X)$, since $L^p_2(X)\subset L^{\sharp,2}_2(X)\subset
L^\8(X)$ for all $p>2$.

We give two proofs of Theorem \ref{thm:GaugeFixing}. For our first proof,
in Section \ref{sec:Continuity}, we essentially follow the strategy of
Uhlenbeck \cite{UhlLp} and apply the method of continuity. The difficult
step here (in establishing `openness' --- see Section
\ref{subsec:Openness}) is to prove that the intrinsic, gauge-invariant
$\sL^{\sharp,2}_{1,A_0}$ and $L^{\sharp,2}_{1,A_0}$ distances in the
quotient $\sB_E^k$ bound the $L^{2\sharp,4}$ and $L^2_{1,A_0}$ norms in the
slice $\bS_{A_0}\subset\sA_E^k$: this is the point in our first proof where
we use the critical-exponent estimates derived in Section \ref{sec:Green}
to control gauge transformations. The proof of `closedness' uses a
compactness argument and is given in Section \ref{subsec:Closedness}.  
 
Our second proof of Theorem \ref{thm:GaugeFixing} occupies 
Sections \ref{sec:SharpGaugeGroup} and \ref{sec:SliceII}. 
In Section \ref{sec:SharpGaugeGroup} we show that the exponential map
$\Exp:\Om^0(\fg_E)\to \sG_E$ extends to a continuous map
$\Exp:L^{\sharp,2}_2(\fg_E)\to L^{\sharp,2}_2(\fg_E)$ and that the
resulting space of $L^{\sharp,2}_2$-gauge transformations
$\sG_E^{2,\sharp,2}$ is a Banach Lie group. In particular,
$L^{\sharp,2}_2$-gauge transformations are {\em continuous\/} and are
contained in $\sG_E^{2,p}$ for every $p>2$. The Sobolev multiplication and
composition results for the critical-exponent norms then allow us to apply
the inverse function theorem directly in Section \ref{sec:SliceII}, while
still ensuring that all constants depend at most on $\nu_0[A_0]$ and
$\|F_{A_0}\|_{L^2}$. We
first use the compactness result of Section \ref{subsec:DistFunction} to
establish the existence of gauge transformations $w$ in $\sG_E^3$ which
minimize the $\sL^{\sharp,2}_{1,A_0}$ and $L^{\sharp,2}_{1,A_0}$ distances
in the quotient $\sB_E^k$. Then, assuming the norm
$\|w(A)-A_0\|_{L^{2\sharp,4}}$ or
$\|w(A)-A_0\|_{L^2_{1,A_0}}$ is sufficiently small, we use the Sobolev
embedding and multiplication theorems of Sections
\ref{sec:SharpSobolev}, \ref{sec:Green}, and \ref{sec:SharpGaugeGroup}
and a quantitative version of the inverse function theorem to prove the
existence of a gauge transformation $v\in\sG_E^3$ such that
$d_{A_0}^*(u(A)-A_0)=0$, $u=vw\in\sG_E^{k+1}$, and
$\|u(A)-A_0\|_{L^{2\sharp,4}}$ and $\|u(A)-A_0\|_{L^2_{1,A_0}}$ are
controlled by $\dist_{\sL^{\sharp,2}_{1,A_0}}([A],[A_0])$ and
$\dist_{L^{\sharp,2}_{1,A_0}}([A],[A_0])$, respectively.

\subsection{Applications and extensions}
\label{subsec:Applic}
Applications of the methods and results of the present article occur in
situations where connections, metrics, or holomorphic structures are
allowed to degenerate and uniform estimates are required for elliptic
operators whose coefficients depend on these degenerating geometric structures.

\subsubsection{Degeneration of anti-self-dual connections} 
While the standard slice theorem is adequate for many applications in
smooth four-manifold topology \cite{DK,FU,FrM}, one finds that it is rather
less adequate for constructing gluing maps and analyzing their asymptotic
behavior in sufficient generality to be useful in approaches to the
Kotschick-Morgan conjecture \cite{FeehanLeness,MorganOzsvath}. That
conjecture asserts that the Donaldson invariants of a four-manifold with
$b^+(X)=1$, computed using metrics lying in different chambers of the
positive cone of $H^2(X;\RR)/\RR^*$, differ by terms depending only the
homotopy type of $X$ \cite{Goettsche,KoM}. 

Taubes' gluing maps can be used \cite{MorganOzsvath} to
construct links of arbitrary lower-level reducibles in
$\barM_\kappa^w\less M_\kappa^w$, where $E$ is now a $\U(2)$ bundle and
$M_\kappa^w$ is the moduli space of anti-self-dual connections on $\su(E)$,
with $w_2(\su(E)) = c_1(E)\pmod{2}$ and $\kappa=-\quarter p_1(\su(E))$.
Any obvious approach to the conjecture places considerable demands
on gluing theory since one must, at least in principle, be able to describe
links of all ideal reducibles and not just isolated special cases.

For a lower-level reducible $[A_1,\bx]\in
M_{\kappa-\ell}^w\times\Sym^\ell(X)$, the Laplacian $d_{A_1}^*d_{A_1}$ has
a kernel of dimension equal to that of $\Stab_{A_1}$, so the Laplacian
$d_A^*d_A$ of connection $A$ necessarily has $\dim \Stab_{A_1}$ `small
eigenvalues' tending to zero as $[A]$ approaches $[A_1,\bx]$ and the
remaining eigenvalues of $d_A^*d_A$ are bounded below by one-half the least
positive eigenvalue of $d_{A_1}^*d_{A_1}$. A similar phenomenon arises in
\cite{TauIndef,TauStable} with the Laplacian
$d_A^+d_A^{+,*}$ on $\Om^+(\fg_E)$ when the `background' connection
$A_1$ has $\Ker d_{A_1}^+d_{A_1}^{+,*}\ne 0$: the small eigenvalues
of $d_A^+d_A^{+,*}$
represent an obstruction to perturbing approximately anti-self-dual
to anti-self-dual connections near $[A_1,\bx]$. 
For such neighborhoods in $\barM_\kappa^w$, Theorems 
\ref{thm:GaugeFixing} and \ref{thm:Slice} are easily modified by replacing
$\nu_0[A_0]$ with a suitable small eigenvalue cutoff, by analogy with
\cite{TauIndef,TauStable}, and proceeding along the lines of \cite[\S 6,
Part 3]{TauFrame}. (Without such modifications, the coordinate balls of
Theorem \ref{thm:Slice} shrink if $[A_0]\in M_\kappa^w$ approaches a
reducible point in $\barM_\kappa^w\less M_\kappa^w$.)  To illustrate
applications of the methods of Sections \ref{sec:SharpSobolev} and
\ref{sec:Green} and to point to possible generalizations of the estimates
in present article, we derive some elliptic estimates for $d_A^++d_A^*$ in
Section \ref{subsec:EllipticEstimatedA+}.

\subsubsection{Degeneration of $\PU(2)$ monopoles} We recall that Viktor
Pidstrigach and Andrei Tyurin proposed a method \cite{PTCambridge,
PTLocal} to prove Witten's conjecture concerning the relation between the
Donaldson and Seiberg-Witten invariants of smooth four-manifolds
\cite{DonSW,Witten}. Their proposal uses a moduli space of 
solutions to the $\PU(2)$ monopole equations, which are a
natural generalization of the $\U(1)$ monopole equations of Seiberg and
Witten and the anti-self-dual equation for $\SO(3)$ connections, to
construct a cobordism between links of compact moduli spaces of $\U(1)$
monopoles of Seiberg-Witten type and the moduli space of anti-self-dual
connections, which appear as singularities in this larger moduli space.

The problems one encounters in attempting to implement this program have
much in common with those encountered in previous attempts to prove the
Kotschick-Morgan conjecture \cite{FL1,FLGeorgia,FL2}.
For the purposes of proving the general gluing
theorem for $\PU(2)$ monopoles (which one needs to describe the
above-mentioned links) we require, among other things, a useful criterion
to detect when a monopole near an ideal point in the Uhlenbeck
compactification of the moduli space of $\PU(2)$ monopoles lies in the
image of a gluing map. An $L^2_{k,A_0}$ (with $k\ge 2$) or even an
$L^p_{1,A_0}$ (with $p>2$) measure of closeness to the image of the
approximate gluing map is not useful for this purpose as the radius of such
balls tend to shrink as the curvature of the connection bubbles,
as in the case of anti-self-connections, and simpler implicit-function
theorem arguments used to show that a point lies in the image of a gluing
map would fail in general. In \cite{FL3} we use 
analogues of Theorems \ref{thm:GaugeFixing} and \ref{thm:Slice} and 
their methods of proof to facilitate arguments that the $\PU(2)$ monopole
gluing maps are `surjective' (in the above sense) and are diffeomorphisms onto
their images. The surjectivity property of Taubes' gluing maps for
anti-self-dual connections is a special case of a more general gluing
result for critical points of the Yang-Mills functional
\cite[Proposition 8.2]{TauFrame}. The gluing maps for anti-self-dual
connections constructed by Donaldson and employed in
\cite{DonConn,DonHCobord,DonPoly,DK,DS} are shown to be surjective in
\cite{DonConn,DK,DS}.

\subsubsection{Degeneration of stable, holomorphic vector bundles}
Via work of Donaldson \cite{DonNS,DonASD} one may view the degeneration
of stable holomorphic vector bundles as a special case of the degeneration
of anti-self-dual connections.  Comparisons of these degenerations have
been given by Li \cite{Li} and Morgan
\cite{MorganComparison} in their work relating the the Gieseker and Uhlenbeck
compactifications of the moduli space of anti-self-dual connections
over a complex, K\"ahler surface. The notion of stability is broadened in
\cite{Bradlow,OTQuaternion,TelemanNonabelian} by considering the vortex and
nonabelian monopole equations over complex, K\"ahler manifolds. In
particular, as explained to the author by Daskalopoulos
\cite{Daskal}, there are potential applications for the types of estimates
considered here when applying methods from gauge theory to study the
degeneration of holomorphic vector bundles near the boundary of
Teichm\"uller space; see
\cite{BradlowDask} for a description of the moduli space of stable,
holomorphic vector bundles over Riemann surfaces.  Such gauge-theoretic
methods are in turn expected to have applications to three-manifold
topology \cite{Daskal}.

\subsubsection*{Acknowledgements}
We warmly thank Tom Mrowka for helpful and stimulating
discussions on gauge theory and his suggestion that $L^4$-balls inject into
the quotient, and warmly thank Cliff Taubes for very helpful discussions 
and especially for explaining his ideas and methods in
\cite{TauPath,TauFrame}. We 
would also like to thank the Harvard Mathematics Department and the
National Science Foundation for their generous hospitality and support
during the preparation of this article.
 
\section{Preliminaries}
\label{sec:Prelim}
We assume throughout this article that $X$ is a closed, connected, smooth,
four-manifold with Riemannian metric $g$. Let $G$ be a compact Lie group
with matrix representation $\rho:G\subset \SO(\EE) = \SO(r)$ where
$\EE\simeq\RR^r$ as a real inner product space, let $P$ be a principal $G$
bundle, and let $E=P\times_\rho\EE$ be the corresponding Riemannian vector
bundle associated to $P$ by the representation $\rho$. Let $\fg_E\subset
\gl(E)$ be the bundle of Lie algebras associated to $P$ via the adjoint
representation $\Ad:G\to \Aut(\fg)$ of $G$ on its Lie algebra $\fg$ and
viewed as a subbundle of $\gl(E)$ via the induced representation
$\rho_*:\fg\subset \so(\EE)$.

Given the covariant derivative $\cov_A:C^\8(E)\to C^\8(T^*X\otimes E)$, we
define the exterior covariant derivative $d_A:\Om^i(E)\to\Om^{i+1}(E)$ in the
usual way by setting $d_A=\cov_A$ on $\Om^0(E)=C^\8(E)$ and extending $d_A$
to $\Om^i(E)=C^\8(\La^i\otimes E)$, where $\La^i := \La^i(T^*X)$, 
according to the rule 
$d_A(\om\wedge v) = d\om\wedge v + (-1)^i\om\wedge d_Av$ for
$\om\in\Om^i(X)$ and $v\in \Om^j(E)$.  

For any integer $k\ge 0$, exponent $1\le p\le\8$, and $L^p_k$ connection
$A_0$ on $E$ we define the $L^p_k$ Sobolev completion,
$L^p_k(\La^\ell\otimes E)$, of $\Om^\ell(E)$ with respect to the norm
$$
\|s\|_{L^p_{k,A_0}(X)} 
:= \left(\sum_{j=0}^k\|\cov_{A_0}^j s\|_{L^p(X)}^p\right)^{1/p}.
$$
We define the action of a $C^\8$ gauge transformation $u\in\sG_E$ on a
$C^\8$ connection $A$ on the bundle $E$ by pushforward, so
$u(A) := A - (d_Au)u^{-1}$. Fix a connection $A_0\in \sA_E$, let
$\sA_E^k = A_0 + L^2_k(\La^1\otimes\fg_E)$, and define
$$
\sG_E^{k+1} := \{u\in L^2_{k+1}(\gl(E)): u\in G\text{ a.e.}\} 
\subset L^2_{k+1}(\gl(E)).
$$
The space $\sG_E^{k+1}$ is a Banach Lie group, with Lie algebra
$T_{\id_E}\sG_E^{k+1}=L^2_{k+1}(\fg_E)$, and acts smoothly on $\sA_E^k$ with
quotient $\sB_E^k := \sA_E^k/\sG_E^{k+1}$ endowed with the quotient $L^2_k$
topology. 

The stabilizer subgroup $\Stab_A\subset\sG_E^{k+1}$ for a connection $A$ on
$E$ always contains the center $Z(G)\subset G$. We let
$\sA_E^{*,k}\subset\sA_E^k$ denote the space of connections $A\in \sA_E^k$
with minimal stabilizer $\Stab_A = Z(G)$ and let $\sB_E^{*,k} =
\sA_E^{*,k}/\sG_E^{k+1}$.  As usual, the stabilizer subgroup
$\Stab_A\subset\sG_E$ can be identified with a closed subgroup of $G\subset
\GL(E|_{x_0})$ for any point $x_0\in X$ by parallel translation with
respect to the connection $A$. Let $\stab_A$ denote the Lie algebra of
$\Stab_A$, so $\stab_A = \Ker\{d_A:L^2_{k+1}(\fg_E)\to
L^2_k(\La^1\otimes\fg_E)\}$.

Throughout the article, we use $c$ or $z$ to denote positive constants which
depend at most on the Riemannian manifold $(X,g)$ and the group $G$;
constants may increase from one line to the next and are not renamed
unless clarity demands otherwise.

\section{The slice theorem} 
\label{sec:SliceI}
In this section we prove the first assertion of Theorem \ref{thm:Slice} ---
see Proposition \ref{prop:Slice} below --- namely, that a small enough
$L^4$-ball $\bB^4_{A_0}(\eps_0)/\Stab_{A_0}$ provides a slice for the
action of $\sG_E^{k+1}$. The proof that the projection
$\pi:\bB^4_{A_0}(\eps_0)/\Stab_{A_0} \to \sB_E^k$ is injective 
(Lemma \ref{lem:injective}) was suggested to us by Mrowka.

Let $k\ge 2$ be an integer. The Banach Lie group $\sG_E^{k+1}$ has Lie
algebra $T_{\id_E}\sG_E^{k+1}=L^2_{k+1}(\fg_E)$ and exponential map
$\Exp:L^2_{k+1}(\fg_E)\to
\sG_E^{k+1}$ given by $\zeta\mapsto u=\Exp\zeta$.
Recall that $\Stab_A=\{\ga\in
\sG_E^{k+1}:\ga(A) = A\}$ may be identified with a Lie subgroup of $G$
and has Lie algebra $\stab_A = \Ker(d_A|_{L^2_{k+1}})$.
The operator $d_A^*:L^2_{k+2}(\La^1\otimes\fg_E)\to L^2_{k+1}(\fg_E)$
has closed range and we have an $L^2$-orthogonal decomposition
\begin{align}
T_{\id_E}\sG_E^{k+1} &= L^2_{k+1}(\fg_E) 
\label{eq:GaugeLieAlgOrthogDecomp}\\
&= (\Ker(d_A|_{L^2_{k+1}}))^\perp \oplus \Ker(d_A|_{L^2_{k+1}}) \notag\\
&= \Imag(d^*_A|_{L^2_{k+2}}) \oplus \Ker(d_A|_{L^2_{k+1}}) \notag\\
&= (\Ker(d_A|_{L^2_{k+1}}))^\perp \oplus \stab_A. \notag
\end{align}
Let $\Stab^\perp_A = L^2_{k+1}\cap \Stab^{\perp}_A
=\Exp((\Ker d_A|_{L^2_{k+1}})^\perp)$,
the second equality following from the Sobolev composition lemma.
The subspace $\Stab^\perp_A\subset \sG_E^{k+1}$ is closed and is a
Banach submanifold of $\sG_E^{k+1}$ with codimension $\dim\stab_A$. {}From
Claim \ref{claim:GaugeGroupSlice} below we see that $\Stab^\perp_A$ is
a slice near $\id_E\in \sG_E^{k+1}$ for the right action of $\Stab_A$ on
$\sG_E^{k+1}$.  

The map $d_A:L^2_{k+1}(\fg_E)\to L^2_k(\La^1\otimes\fg_E)$
has closed range and so
we have an $L^2$-orthogonal decomposition 
\begin{align}
T_A\sA_E^k &= L^2_k(\La^1\otimes\fg_E) \label{eq:TangentPairOrthogDecomp}\\
&= \Imag(d_A|_{L^2_{k+1}}) \oplus \Ker(d_A^*|_{L^2_k}) \notag\\
&= \Imag(d_A|_{L^2_{k+1}}) \oplus \bK_A,\notag
\end{align}
of the tangent space to the space of $L^2_k$ connections at the point $A$,
where we set where $\bK_A  = \Ker (d^*_A|_{L^2_k})$.

The {\em slice\/} $\bS_A\subset
\sA_E^k$ through a connection $A$ is given by $\bS_A =
A + \bK_A$. If
$\pi$ is the projection from $\sA^k_E$ onto $\sB_E^k =
\sA^k_E/\sG_E^{k+1}$, denoted by $A\mapsto [A]$, we let 
\begin{align*}
\bB_A(\eps) &= \{A_1\in \bS_A:\|A_1-A\|_{L^2_{k;A}}< \eps\} \\
&= A + \{a\in \bK_A:\|a\|_{L^2_{k;A}}< \eps\}
\end{align*}
be the open $L^2_k$-ball in $\bS_A$ with center $A$ and $L^2_{k,A}$-radius
$\eps$. Similarly, we let 
\begin{align*}
\bB^4_A(\eps) &= \{A_1\in \bS_A:\|A_1-A\|_{L^4}< \eps\} \\
&= A + \{a\in \bK_A:\|a\|_{L^4}< \eps\}
\end{align*}
be the open ball in $\bS_A$ with center $A$ and $L^4$-radius $\eps$.

The proof that the quotient space $\sA_E^k$ is Hausdorff makes
use of the following well-known technical result 
\cite[Proposition A.5]{FU}. Note that the space $\sG_E^2$ is
neither a Banach Lie group nor does it act smoothly on $\sA_E^k$ for $k\ge 1$.

\begin{lem}\label{lem:GaugeSequence}
Let $E$ be a Hermitian bundle over a Riemannian manifold $X$ and let $k\ge
2$ be an integer.  Suppose $\{A_\al\}$ and $\{B_\al\}$ are sequences of
$L^2_k$ unitary connections on $E$ and that $\{u_\al\}$ is a sequence in
$\sG_E^2$ such that $u_\al(A_\al) = B_\al$. Then the following hold.
\begin{enumerate}
\item The sequence $\{u_\al\}$ is in $\sG_E^{k+1}$;
\item If $\{A_\al\}$ and $\{B_\al\}$ converge in 
$\sA_E^k$ to limits $A_\8$, $B_\8$, then  
there is a subsequence $\{\al'\}\subset\{\al\}$ such that
$\{u_{\al'}\}$ converges in $\sG_E^{k+1}$ to $u_\8$ and $B_\8 = u_\8(A_\8)$.
\end{enumerate}
\end{lem}

We shall need the following quantitative version of the inverse function
theorem here and especially in Section \ref{sec:SliceII}:

\begin{thm}\label{thm:InverseFT}
Let $\bPhi:\bE\to \bF$ be a $C^\ell$ map of Banach spaces, for some
$\ell\ge 1$,  
such that the differential
$(D\bPhi)_{x_0}:\bE\to \bF$ has a continuous inverse
$(D\bPhi)_{x_0}^{-1}:\bF\to \bE$ 
satisfying
$$
\|(D\bPhi)_{x_0}^{-1}\| \le K \quad\text{and}\quad 
\|(D\bPhi)_x-(D\bPhi)_{x_0}\| \le \half K^{-1},
\quad\text{if}\quad \|x-x_0\|\le\de,
$$
for some positive constants $K$ and $\de$. Then the following hold:
\begin{enumerate}
\item The restriction of $\bPhi$ to the ball $U=B^\bE(x_0,\de)$ is
injective and 
$\bPhi(U)=V$ is an open set in $\bF$ containing the ball
$B^\bF(\bPhi(x_0),\de/(2K))$; 
\item The inverse map $\bPhi^{-1}:V\to U$ is $C^\ell$;
\item If $x_1,x_2\in B^\bE(x_0,\de)$, then $\|x_1-x_2\|\le
2K\|\bPhi(x_1)-\bPhi(x_2)\|$.  
\end{enumerate}
\end{thm}

For quantitative comparisons in this section, the following elementary fact
will suffice:

\begin{lem}\label{lem:OneSidedInverse}
Let $\bE$, $\bF$ be Banach spaces and let $T\in\Hom(\bE,\bF)$ have
a right (left) inverse $S$. If $\tT\in\Hom(\bE,\bF)$ satisfies
$\|\tT-T\|<\|S\|^{-1}$, then $\tT$ also has a right (left) inverse.
\end{lem}

\begin{pf}
If $S\in\Hom(\bF,\bE)$ is a right inverse for $T$, so $TS=\id_\bF$, then
$\|(\tT-T)S\|\le \|\tT-T\|\|S\| < 1$ and $\id_\bE+(\tT-T)S$ is an invertible
element of the Banach algebra $\End(\bE)$. Define
$\tS=S(1+(\tT-T)S)^{-1}$, so $\tS T=\id_\bE$ and $\tS$ is a right inverse for
$\tT$. Similarly for left inverses.
\end{pf}

This consequence of the usual characterization of
invertible elements of a Banach algebra will be invoked in the proof of
Lemma \ref{lem:LocalDiffeo}. 

\begin{prop}\label{prop:Slice}
Let $X$ be a closed, Riemannian four-manifold.  Then there is a
positive constant $z$ with the following significance.  Let $E$ be a
$G$ bundle over $X$. Suppose that $k\ge 2$ is an integer.
Given $[A_0]$ in $\sB_E^k$, let
$\nu_0[A_0]$ be the least positive eigenvalue of the Laplacian
$\De_{A_0}^0$ and let $\eps_0$ be a constant satisfying 
$0<\eps_0<z(1+\nu_0[A_0]^{-1/2})^{-1}$. 
Then the following hold:
\begin{enumerate}
\item The space $\sB_E^k$ is Hausdorff;
\item The subspace $\sB_E^{*,k}\subset \sB_E^k$ is
open and is a $C^\8$ Banach manifold with local parametrizations given by
$\pi:\bB^4_{A_0}(\eps_0)\to \sB_E^{*,k}$;
\item The projection $\pi:\sA_E^{*,k} \to
\sB_E^{*,k}$ is a $C^\8$ principal $\sG_E^{k+1}/\Center(G)$ bundle; 
\item The projection  
$\pi:\bB^4_{A_0}(\eps_0)/\Stab_{A_0}\to \sB_E^k$ is a
homeomorphism onto an open neighborhood of $[A_0]\in
\sB_E^k$ and a diffeomorphism on the and a diffeomorphism on the open
subset where $\Stab_{A_0}/\Center(G)$ acts freely.
\end{enumerate}
\end{prop}

\begin{pf}
The stabilizer $\Stab_{A_0}$  acts freely on $\sG_E^{k+1}$ and thus
on the Banach manifold $\sG_E^{k+1}\times 
\bS_{A_0}^k$ by $(u,A)\mapsto \ga\cdot(u,A) = 
(u\ga^{-1},\ga(A))$ and so the
quotient $\sG_E^{k+1}\times_{\Stab_{A_0}}\bS_{A_0}$ is
again a Banach manifold. We define a smooth map   
\begin{equation}
\bPsi: \sG_E^{k+1}\times_{\Stab_{A_0}} \bS_{A_0}
\rightarrow \sA_E^k, \qquad
[u,A]  \mapsto u(A). \label{eq:PsiDefn}
\end{equation}
Our main task is to show that the map $\bPsi$ is (i) a local
diffeomorphism onto its image and (ii) injective upon restriction to a
sufficiently small neighborhood $\sG_E^{k+1}\times_{\Stab_{A_0}}
\bB^4_{A_0}(\eps_0)$. 
Given $\de_0>0$, let $B_{\id_E}(\de_0)$ be the ball
$\{u\in \sG_E^{k+1}: \|u-\id_E\|_{L^2_{k+1;A_0}}<\de_0\}$ 
and let $B^\perp_{\id_E}(\de_0)=
B_{\id_E}(\de_0)\cap\Stab^\perp_{A_0}$. 

\begin{claim}\label{claim:GaugeGroupSlice}
For small enough $\de=\de(A_0,k)$, the ball $B_{\id_E}(\de)$ is
diffeomorphic to an open neighborhood in $B^\perp_{\id_E}(\de)\times
\Stab_{A_0}$, with inverse map given by $(u_0,\ga)\mapsto u =
u_0\gamma$.
\end{claim}  

\begin{pf}
The differential of the multiplication map
$$
\Stab^\perp_{A_0}\times \Stab_{A_0}\rightarrow \sG_E^{k+1}, 
\qquad (u_0,\ga)\mapsto u_0\ga,
$$ 
at $(\id_E,\id_E)$ is given by 
$$
\Ker(d_{A_0}|_{L^2_{k+1}})^\perp
\oplus \stab_{A_0}\rightarrow L^2_{k+1}(\fg_E)
\qquad (\zeta,\chi)\mapsto u_0\zeta\ga + u_0\ga\chi,
$$
and so is just the identity map with respect to the $L^2$-orthogonal
decomposition \eqref{eq:GaugeLieAlgOrthogDecomp} of the range.  Hence, the
Banach space implicit function theorem implies that there is a
diffeomorphism from an open neighborhood of $(\id_E,\id_E)$ onto an open
neighborhood of $\id_E\in
\sG_E^{k+1}$. For small enough $\de$, we may suppose that if
$u\in B_{\id_E}(\de)$, then $u$ can
be written uniquely as $u=u_0\gamma$ with 
$u_0\in  B^\perp_{\id_E}(\de)$ and
$\gamma\in\Stab_{A_0}$.  
\end{pf}

\begin{lem}\label{lem:LocalDiffeo}
For any $0<\eps_0<\half(1+\nu_0[A_0]^{-1/2})^{-1}$, 
the map $\bPsi$ is a local
diffeomorphism from $\sG_E^{k+1}\times_{\Stab_{A_0}}
\bB^4_{A_0}(\eps_0)$ onto its image in $\sA_E^k$.
\end{lem}

\begin{pf}
We first restrict the map $\bPsi$ to a neighborhood 
$B_{\id_E}(\de_0)\times_{\Stab_{A_0}}\bS_{A_0}$, which is
diffeomorphic to the neighborhood
$B^\perp_{\id_E}(\de)\times\bS_{A_0}$ in  
$\Stab^\perp_{A_0}\times \bS_{A_0}$ by Claim
\ref{claim:GaugeGroupSlice}. The differential of the induced map
\begin{equation}
\bPsi:\Stab^\perp_{A_0} \times \bS_{A_0}
\to \sA_E^k, \qquad
(u,A) \mapsto u(A), \label{eq:PsiDefnPerp}
\end{equation}
at $(\id_E,A) := (\id_E,A_0+a_0)$ is given by
\begin{align*}
(D\bPsi)_{(\id_E,A)}: &T_{\id_E}\Stab^\perp_{A_0}
\oplus T_{A}\bS_{A_0}
\to T_{A}\sA_E^k, \\
&(\zeta,a) \mapsto -d_A\zeta + a 
= -d_{A_0}\zeta - [a_0,\zeta] + a,
\end{align*}
where we recall that $T_{A}\bS_{A_0} = \bK_{A_0} 
= \Ker (d_{A_0}^*|_{L^2_k})$ and
$$
T_{\id_E}\Stab^\perp_{A_0} = (\Ker(d_{A_0}|_{L^2_{k+1}}))^\perp
= \Imag(d^*_{A_0}|_{L^2_{k+2}}).
$$
Using the $L^2$-orthogonal
decomposition \eqref{eq:TangentPairOrthogDecomp} of the range we see that
the map 
$$
-d_{A_0}\oplus \id_E: (\Ker(d_{A_0}|_{L^2_1}))^\perp
\oplus \Ker (d_{A_0}^*|_{L^2_1}) \to \Imag(d_{A_0}|_{L^2_1})
\oplus \Ker (d_{A_0}^*|_{L^2_1})
$$
given by $(\zeta,b) \mapsto -d_{A_0}\zeta + b$ is a Hilbert space
isomorphism. More explicitly, the operator 
$$
d_{A_0}:(\Ker(d_{A_0}|_{L^2_1}))^\perp \to \Imag(d_{A_0}|_{L^2_1}) =
(\Ker(d_{A_0}^*|_{L^2}))^\perp 
$$
has a two-sided inverse 
$$
G^0_{A_0}d^*_{A_0}: \Imag(d_{A_0}|_{L^2_1})\to (\Ker(d_{A_0}|_{L^2_1}))^\perp,
$$
where $G^0_{A_0}$ is the Green's operator for the Laplacian $\De^0_{A_0} =
d^*_{A_0}d_{A_0}$: indeed, $G^0_{A_0}d^*_{A_0}d_{A_0} = G^0_{A_0}\De^0_{A_0}$
is the $L^2$-orthogonal projection $\Pi^0_{A_0}$
from $L^2_1(\La^1\otimes\fg_E)$ onto $(\Ker(d_{A_0}|_{L^2_1}))^\perp$
and $d_{A_0}G^0_{A_0}d^*_{A_0}$ is the $L^2$-orthogonal projection 
$\Pi^{1,\perp}_{A_0}=\id-\Pi^1_{A_0}$ from $L^2(\La^1\otimes\fg_E)$
onto $(\Ker(d_{A_0}^*|_{L^2}))^\perp$, as
$$
d_{A_0}^*(\id - d_{A_0}G^0_{A_0}d^*_{A_0}) = 0.
$$ 
For $\zeta\in (\Ker(d_{A_0}|_{L^2_1}))^\perp$ and 
$b=d_{A_0}\zeta\in \Imag(d_{A_0}|_{L^2_1})$, we have
\begin{align*}
\|G^0_{A_0}d^*_{A_0}b\|_{L^2_{1,A_0}} 
&= \|G^0_{A_0}\De^0_{A_0}\zeta\|_{L^2_{1,A_0}} 
= \|\Pi^0_{A_0}\zeta\|_{L^2_{1,A_0}}
= \|\zeta\|_{L^2_{1,A_0}} \\
&\le \|d_{A_0}\zeta\|_{L^2} + \|\zeta\|_{L^2} 
\le (1+\nu_0^{-1/2})\|d_{A_0}\zeta\|_{L^2} = (1+\nu_0^{-1/2})\|b\|_{L^2}
\end{align*}
and so $G^0_{A_0}d^*_{A_0}$ has $\Hom(L^2,L^2_{1,A_0})$ operator norm bound
$$
\|G^0_{A_0}d^*_{A_0}\| \le 1+\nu_0^{-1/2}.
$$
The Sobolev embedding $L^2_1\subset L^4$ and Kato's inequality
imply that
$$
\|d_A\zeta - d_{A_0}\zeta\|_{L^2}
\le \|[a_0,\zeta]\|_{L^2}
\le 2\|a_0\|_{L^4}\|\zeta\|_{L^4}
\le 2\|a_0\|_{L^4}\|\zeta\|_{L^2_{1,A_0}},
$$
and so $d_A-d_{A_0}$ has $\Hom(L^2_{1,A_0},L^2)$ operator norm bound
$$
\|d_A - d_{A_0}\| \le 2\|a_0\|_{L^4}.
$$
In particular, we see that 
$(D\bPsi)_{(\id_E,A_0)}^{-1} = G^0_{A_0}d^*_{A_0}\oplus \id
= G^0_{A_0}d^*_{A_0}\oplus G^0_{A_0}\De^0_{A_0}$ satisfies
$$
\|(D\bPsi)_{(\id_E,A_0)}^{-1}\| \le 1+\nu_0^{-1/2} \quad\text{and}\quad
\|(D\bPsi)_{(\id_E,A)}-(D\bPsi)_{(\id_E,A_0)}\| \le 2\|a_0\|_{L^4}.
$$
Hence, Lemma \ref{lem:OneSidedInverse} implies that if $\|a_0\|_{L^4} <
\half(1+\nu_0^{-1/2})^{-1}$, then the operator
$$
(D\bPsi)_{(\id_E,A)}:
(\Ker(d_{A_0}|_{L^2_1}))^\perp
\times \Ker (d_{A_0}^*|_{L^2_1}) \to L^2(\La^1\otimes\fg_E)
$$
is an isomorphism from $L^2_1$ to $L^2$ and
restricts to a bounded linear map from $L^2_{k+1}$ to $L^2_k$.
Provided $(D\bPsi)_{(\id_E,A)}:L^2_{k+1}\to L^2_k$ is bijective, the open
mapping theorem guarantees the existence of a bounded inverse
$(D\bPsi)_{(\id_E,A)}^{-1}:L^2_k\to L^2_{k+1}$.
If $(D\bPsi)_{(\id_E,A)}(\zeta,a)=0$ for $(\zeta,a)\in L^2_{k+1}$, then 
$(\zeta,a)$ is zero in $L^2_1$ and thus zero in $L^2_{k+1}$, so
$(D\bPsi)_{(\id_E,A)}$ is injective. If $b\in
L^2_k(\La^1\otimes\fg_E)$, then $b = (D\bPsi)_{(\id_E,A)}(\zeta,a) = -d_A\zeta +
a$ for some $(\zeta,a)\in (\Ker(d_{A_0}|_{L^2_1}))^\perp
\times \Ker (d_{A_0}^*|_{L^2_1})$. As $d_{A_0}^*a=0$, we have
$$
d_{A_0}^*d_A\zeta = -d_{A_0}^*b \in L^2_{k-1}
$$
and $d_{A_0}^*d_A:L^2_{k+1}\to L^2_{k-1}$ is an elliptic operator with
$L^2_{k-1}$ coefficients. Thus, $\zeta\in L^2_{k+1}$,  so $a=b+d_A\zeta
\in L^2_k$, and $(D\bPsi)_{(\id_E,A)}$ is surjective.

Combining the above observations, we see that the operator
$$
(D\bPsi)_{(\id_E,A)}:(\Ker(d_{A_0}|_{L^2_{k+1}}))^\perp
\oplus \Ker (d_{A_0}^*|_{L^2_k}) \to L^2_k(\La^1\otimes\fg_E),
$$
is an isomorphism for all $A=A_0+a_0$ with
$\|a_0\|_{L^4} < \eps_0 = \half(1+\nu_0^{-1/2})^{-1}$. 
So, by the Banach space implicit function theorem,
there are positive constants $\eps=\eps(A,k)$ and $\de=\de(A,k)$ and an open
neighborhood $U_A\subset
\sA_E^k$ such that the map
$$
\bPsi: B^\perp_{\id_E}(\de) 
\times \bB_A(\eps)\to U_A, \qquad
(u,A_1)\mapsto u(A_1),
$$ 
with $\bB_A(\eps)\subset \bB^4_{A_0}(\eps_0)$,
gives a diffeomorphism from an open neighborhood of
$(\id_E,A)$ onto an open neighborhood of $A$. 
In particular,
we obtain a map $U_A\to \Stab^\perp_{A_0}$, given by
$A_1 \mapsto u = u_{A_1}$, such that
$$
\bPsi^{-1}(A_1) = (u,u^{-1}(A_1))
\in B^\perp_{\id_E}(\de) \times \bB_{A}(\eps)
\subset \Stab^\perp_{A_0}\times \bB^4_{A_0}(\eps_0).
$$
Hence, for any $A_1\in U_A$ there is a unique
$u\in  B^\perp_{\id_E}(\de)$ such that
$u^{-1}(A_1)-A_0\in\bK_{A_0}$: 
\begin{equation}
d^*_{A_0}(u^{-1}(A_1)-A_0) = 0.
\label{eq:Coulomb}
\end{equation}
The neighborhood $\bB^4_{A_0}(\eps_0)$ is $\Stab_{A_0}$-invariant: if $A\in
\bB_{A_0}(\eps)$ and 
$\ga\in\Stab_{A_0}$, then 
$$
\|\ga(A)-A_0\|_{L^4} = \|A-\ga^{-1}(A_0)\|_{L^4} = \|A-A_0\|_{L^4} <\eps,
$$
and
\begin{align*}
d^*_{A_0}\left(\ga(A)-A_0\right) 
&= \ga\left(d^*_{\ga^{-1}(A_0)}(A-\ga^{-1}(A_0))\right) \\
&= \ga\left(d^*_{A_0}(A-A_0)\right) = 0,
\end{align*} 
so $\ga(A)\in \bB_{A_0}(\eps)$.

The group $\sG_E^{k+1}$ acts on $\sG_E^{k+1}\times
\bS_{A_0}$ by $(u,A)\mapsto (vu,A)$, and so gives a
diffeomorphism  
$$
B_{\id_E}(\de) \times \bB^4_{A_0}(\eps_0)
\to B_v(\de) \times \bB^4_{A_0}(\eps_0), \quad
(u,A) \to (vu,A),
$$
and as this action commutes with the given action of
$\Stab_{A_0}$, it descends to a diffeomorphism 
$$
B_{\id_E}(\de) \times_{\Stab_{A_0}} \bB^4_{A_0}(\eps_0)
\to B_v(\de) \times_{\Stab_{A_0}} \bB^4_{A_0}(\eps_0), \quad
[u,A] \to [vu,A],
$$
for each $v\in \sG_E^{k+1}$. Consequently, the 
$\sG_E^{k+1}$-equivariant map
$$
\sG_E^{k+1}\times_{\Stab_{A_0}}
\bB^4_{A_0}(\eps_0)\to \sA_E^k
$$
is a local diffeomorphism onto its image, as desired.
\end{pf}

Plainly, $[\ga(A)] =
[A]$ for each $\ga\in\Stab_{A_0}$ and $A \in
\bB^4_{A_0}(\eps_0)$  
and hence, the projection $\pi:\bB^4_{A_0}(\eps_0)\to
\sA_E^k$ factors through
$\bB^4_{A_0}(\eps)/\Stab_{A_0}$. 

\begin{lem}\label{lem:injective}
There is a positive constant $z$ with the following significance.
Let $\nu_0[A_0]$ be the least positive eigenvalue of the Laplacian
$\De_{A_0}^0$. Then for any constant $\eps_0$ satisfying
$0<\eps_0<z(1+\nu_0[A_0]^{-1/2})^{-1}$, the projection map
$\pi:\bB^4_{A_0}(\eps_0)/\Stab_{A_0}\to\sB_E^k$ is injective.
\end{lem}

\begin{pf}
Suppose $A_i\in\bB^4_{A_0}(\eps_0)$ for $i=1,2$ and that
$[A_1]=[A_2]\in \sB_E^k$, so 
$u(A_1)=A_2$ for some $u\in \sG_E^{k+1}$.
Since $u(A_0) = A_0 - (d_{A_0}u)u^{-1}$, we see that $u\in\Stab_{A_0}$ if
and only $d_{A_0}u = 0$. Here, we view $u\in L^2_{k+1}(\gl(E))$
via the isometric embedding $\sG_E^{k+1}\subset L^2_{k+1}(\gl(E))$
and write
$$
u = u_0 - \ga,
$$
where $u_0 \in (\Ker d_{A_0})^\perp$ and $\ga \in \Ker
d_{A_0}$. We claim that $u_0=0$, so $u=\ga\in\Stab_{A_0}$.  

Since $u(A_1):= A_1 - (d_{A_1}u)u^{-1}=A_2$, we have
$A_2u = A_1u - d_{A_1}u = A_1u - d_{A_0}u - [A_1-A_0,u]$, and therefore
$$
d_{A_0}u_0 = d_{A_0}u = u(A_1-A_0) - (A_2-A_0)u.
$$
Since $d_{A_0}^*(A_i-A_0)=0$ for $i=1,2$, we obtain
\begin{align*}
d_{A_0}^*d_{A_0}u_0 &= -*(d_{A_0}u\wedge *(A_1-A_0)) + ud_{A_0}^*(A_1-A_0) \\
&\qquad - (d_{A_0}^*(A_2-A_0))u + *(*(A_2-A_0)\wedge d_{A_0}u) \\
&= -*(d_{A_0}u_0\wedge *(A_1-A_0)) + *(*(A_2-A_0)\wedge d_{A_0}u_0).
\end{align*}
Integrating by parts gives 
$$
\|d_{A_0}u_0\|_{L^2}^2 = (d_{A_0}^*d_{A_0}u_0,u_0)_2 \le
\|d_{A_0}^*d_{A_0}u_0\|_{L^{4/3}}\|u_0\|_{L^4}.
$$ 
Kato's inequality and the embedding $L^2_1\subset L^4$ gives
$\|u_0\|_{L^4}\le c(\|d_{A_0}u_0\|_{L^2} + \|u_0\|_{L^2})$, so the
eigenvalue estimate 
$\|u_0\|_{L^2}\le \nu_0^{-1/2}\|d_{A_0}u_0\|_{L^2}$ gives $\|u_0\|_{L^4} \le
c(1+\nu_0^{-1/2})\|d_{A_0}u_0\|_{L^2}$ and thus
$$
\|d_{A_0}u_0\|_{L^2}^2 
\le (1+\nu_0^{-1/2})\|d_{A_0}^*d_{A_0}u_0\|_{L^{4/3}}\|d_{A_0}u_0\|_{L^2}. 
$$
Therefore, if $d_{A_0}u_0\not\equiv 0$, the 
preceding expression for $d_{A_0}^*d_{A_0}u_0$ yields
\begin{align*}
\|d_{A_0}u_0\|_{L^2} &\le c(1+\nu_0^{-1/2})\|d_{A_0}^*d_{A_0}u_0\|_{L^{4/3}} \\
&\le c(1+\nu_0^{-1/2})
\|d_{A_0}u_0\|_{L^2}(\|A_1-A_0\|_{L^4} + \|A_2-A_0\|_{L^4}),
\end{align*}
and so we have
$$
1 \le c(1+\nu_0^{-1/2})(\|A_1-A_0\|_{L^4} + \|A_2-A_0\|_{L^4}) 
\le c(1+\nu_0^{-1/2})\eps_0.
$$
which gives a contradiction for $\eps_0<c^{-1}(1+\nu_0^{-1/2})^{-1}$
\end{pf}

We now return to consider the local diffeomorphism $\bPsi$ of Claim
\ref{lem:LocalDiffeo}. Suppose $\bPsi[u_1,A_1] = \bPsi[u_2,A_2] \in \sA_E^k$,
where $[u_1,A_1], [u_2,A_2] \in \sG_E^{k+1}\times_{\Stab_{A_0}}
\bB^4_{A_0}(\eps_0)$, and so $u_1(A_1)=u_2(A_2)\in\sA_E^k$ and hence 
$[A_1]=[A_2]\in\sB_E^k$. Provided $\eps_0$ also satisfies the constraints
of Claim \ref{lem:injective}, we have $u_2^{-1}u_1=\ga\in\Stab_{A_0}$
and $\ga(A_1)=A_2$. Hence $[u_2,A_2] = [u_1\ga^{-1},\ga(A_1)] =
[u_1,A_1]$, so $\bPsi$ is injective and therefore a diffeomorphism onto
$\sA_E^k$. 

The map $\pi:\bB^4_{A_0}(\eps_0)/\Stab_{A_0}\to \sB_E^k$ can be factored as
the composition of the inclusion $A\mapsto (\id_E,A)$ of
$\bB^4_{A_0}(\eps_0)$ into 
$\sG_E^{k+1}\times\bB^4_{A_0}(\eps_0)$, the projection onto the
$\Stab_{A_0}$-quotient
$\sG_E^{k+1}\times_{\Stab_{A_0}}\bB^4_{A_0}(\eps_0)$, the diffeomorphism
$\bPsi$ of $\sG_E^{k+1}\times_{\Stab_{A_0}}\bB^4_{A_0}(\eps_0)$ with
$\sA_E^k$ and the projection from $\sA_E^k$ onto the $\sG_E^{k+1}$-quotient
$\sB_E^k=\sA_E^k/\sG_E^{k+1}$. Hence,
$\pi$ is a homeomorphism onto an
open neighborhood of $[A_0]$ in $\sB_E^k$ and a diffeomorphism on the
open subset where $\Stab_{A_0}/\Center(G)$ acts freely.

\begin{claim}\label{claim:Hausdorff}
The quotient space $\sB_E^k$ is Hausdorff.
\end{claim}

\begin{pf}
Let $\Ga$ be the subspace
$\{\{A,u(A)\}:A\in \sA_E^k\text{
and }u\in \sG_E^{k+1}\}$ of $\sA_E^k\times
\sA_E^k$. If $\{(A_\al),u_\al(A_\al)\}$ is a
sequence in $\Ga$ which converges in $L^2_k$ to a point $\{A_\8,B_\8\}$,
then Lemma \ref{lem:GaugeSequence} implies that there is a subsequence
$\{\al'\}\subset\{\al\}$ such that $\{u_\al\}$ converges in $L^2_{k+1}$ to
$u_\8\in\sG_E^{k+1}$ and $u_\8(A_\8)=B_\8$. Thus, $\Ga$ is closed and the
quotient $\sA_E^k/\sG_E^{k+1}$ is Hausdorff.
\end{pf}

Claim \ref{claim:Hausdorff} gives Assertion (1) of the proposition and
Assertions (2), (3), and (4) now follow from the preceding arguments and
Lemma \ref{lem:injective}. This completes the proof of the proposition.
\end{pf}

\section{Critical-exponent Sobolev norms}
\label{sec:SharpSobolev}
We now describe the basic properties of the 
critical-exponent norms and corresponding
Banach spaces introduced by Taubes in
\cite{TauPath,TauFrame,TauStable,TauConf}. In particular, we give the basic
embedding, multiplication, and composition lemmas we need to complete the
proof of our slice theorem. We shall make frequent use of the pointwise
Kato inequality, $|d|v||\le |\cov_Av|$ for $v\in \Om^0(E)$, so that the norms
of the embedding and multiplication maps depend at most on the Riemannian
manifold $(X,g)$. Moreover, for simplicity, we confine our attention
to the case of closed four-manifolds: there are obvious analogues of the
Sobolev lemmas described here for any
$n$-manifold, with $n>2$. Similarly, extensions are possible to the case of
complete manifolds bounded geometry (bounded curvature and injectivity
radius uniformly bounded from below) --- see \cite{Adams, Aubin} for further
details for Sobolev embedding results in those situations and for the
construction of Green kernels. 
We refer the reader to the monograph of R. Adams \cite{Adams} for a
comprehensive treatment of Sobolev spaces and to that of E. Stein
\cite{Stein} for a treatment based on potential functions. 

Throughout this section, $A$, $B$ denote $C^\8$ orthogonal connections on
Riemannian vector bundles $E$, $F$ over $X$ with $C^\8$ sections $u$, $v$,
respectively.  We first have the following analogues of the $L^2$ and $L^4$
norms,
\begin{align}
\|u\|_{L^\sharp(X)} &= \sup_{x\in X}\|\dist^{-2}(x,\cdot)|u|\|_{L^1(X)}, \\
\|u\|_{L^{2\sharp}(X)} &= \sup_{x\in X}\|\dist^{-1}(x,\cdot)|u|\|_{L^2(X)},
\notag 
\end{align}
where $\dist(x,y)$ denotes the geodesic distance between points $x$ and $y$
in $X$ defined by the metric $g$; these norms have the same behavior as the
$L^2$ and $L^4$ norms with respect to constant rescalings of
the metric $g$ --- the $L^\sharp$ norm on two-forms and the $L^{2\sharp}$
norm on one-forms are {\em scale invariant\/}. Indeed, one sees this by
noting that if $g\mapsto \tg = \la^{-2}g$, then
$\dist_{\tg}(x,y)=\la^{-1}\dist_g(x,y)$ and $dV_{\tg}=\la^{-4}dV_g$, while
for any $a\in\Om^1(E)$ and $v\in\Om^2(E)$, we have
$|a|_{\tg}=\la|a|_g$, and $|v|_{\tg}=\la^2|v|_g$.

Next, we define analogues of the $L^2_1$ and $L^2_2$ norms
\begin{align*}
\|u\|_{L^2_{1,A}(X)} &=\|\cov_Au\|_{L^2(X)} + \|u\|_{L^2(X)}, \\
\|u\|_{L^2_{2,A}(X)} &=\|\cov_A^2u\|_{L^2(X)} + \|\cov_Au\|_{L^2(X)} +
\|u\|_{L^2(X)}, 
\end{align*}
and set 
\begin{align}
\|u\|_{L^\sharp_{1,A}(X)} 
&=\|\cov_Au\|_{L^\sharp(X)} + \|u\|_{L^{2\sharp}(X)} + \|u\|_{L^\sharp(X)}, \\
\|u\|_{L^\sharp_{2,A}(X)} 
&=\|\cov_A^*\cov_Au\|_{L^\sharp(X)} + \|u\|_{L^\sharp(X)}, \notag 
\end{align}
where $\cov_A^* = -*\cov_A*:\Om^1(E)\to\Om^0(E)$ is the $L^2$-adjoint of
the map $\cov_A:\Om^0(E)\to\Om^1(E)$.

Finally, we define analogues of the $C^0\cap L^2_2$ norm
$$
\|u\|_{C^0\cap L^2_{2,A}(X)} = \|u\|_{C^0(X)} + \|u\|_{L^2_{2,A}(X)},
$$
and set
\begin{align}
\|u\|_{L^{\sharp,2}(X)} &= \|u\|_{L^\sharp\cap L^2(X)} =
\|u\|_{L^\sharp(X)} + \|u\|_{L^2(X)}, \\
\|u\|_{L^{2\sharp,4}(X)} &= \|u\|_{L^{2\sharp}\cap L^4(X)} =
\|u\|_{L^{2\sharp}(X)} + \|u\|_{L^4(X)}, \notag  \\
\|u\|_{L^{\sharp,2}_{1,A}(X)} &= \|u\|_{L^\sharp_{1,A}\cap L^2_{1,A}(X)} =
\|u\|_{L^\sharp_{1,A}(X)} + \|u\|_{L^2_{1,A}(X)}, \notag  \\
\|u\|_{L^{\sharp,2}_{2,A}(X)} &= \|u\|_{L^\sharp_{2,A}\cap L^2_{2,A}(X)} =
\|u\|_{L^\sharp_{2,A}(X)} + \|u\|_{L^2_{2,A}(X)} \notag.
\end{align}
It might have appeared, at first glance, a little more natural to continue the
obvious pattern and instead define $\|u\|_{L^\sharp_{2,A}(X)}$ using
$\|\cov_A^2u\|_{L^\sharp(X)}$: as we shall below, though, the given
definition is most useful in practice. For related reasons, if $u\in\Om^1(E) =
\Om^0(\La^1\otimes E)$, 
it is convenient to define the norm $\|u\|_{L^\sharp_{1,A}(X)}$ by
\begin{equation}
\|u\|_{L^\sharp_{1,A}(X)} 
=\|\cov_A^*u\|_{L^\sharp(X)} + \|u\|_{L^{2\sharp}(X)} +
\|u\|_{L^\sharp(X)}.
\label{eq:LSharp1AdjointCov}
\end{equation}
Let $L^\sharp(X)$ be the Banach space completion of $C^\8(X)$ with respect to
the norm $\|\cdot\|_{L^\sharp}$ and similarly define the remaining Banach
spaces above.

We have the following extensions of the standard Sobolev embedding theorem
\cite{FU,Palais}: their proofs are given in the next section. See also 
\cite{DonApprox}, \cite{ParkerTaubes}, \cite{TauPath}, \cite[\S 6]{TauFrame}, 
\cite[Eq. (3.4) \& \S 5]{TauStable}, and \cite[Lemma 4.7]{TauConf}.

\begin{lem}\label{lem:EmbeddingSobolevInto*}
The following are continuous embeddings:
\begin{enumerate}
\item $L^p_k(E)\subset L^\sharp_k(E)$, for $k=0,1,2$ and all $p>2$;
\item $L^q(E)\subset L^{2\sharp}(E)$, for all $q>4$;
\item $L^2_1(E)\subset L^{2\sharp}(E)$.
\end{enumerate}
\end{lem}

In the reverse direction we have:

\begin{lem}\label{lem:Embedding*IntoSobolev}
The following are continuous embeddings:
\begin{enumerate}
\item $L^\sharp(E)\subset L^1(E)$ and $L^{2\sharp}(E)\subset L^2(E)$;
\item $L^\sharp_2(E)\subset C^0\cap L^2_1(E)$;
\end{enumerate}
\end{lem}

We next consider the extension of the standard Sobolev multiplication lemma
\cite{FU,Palais}. While there is no continuous multiplication map
$L^2_2\times L^2_2 \to L^2_2$, it is worth observing that there is a
continuous bilinear map $C^0\cap L^2_2(E)\times C^0\cap L^2_2(F)
\to C^0\cap L^2_2(E\otimes F)$ given by $(u,v)\mapsto u\otimes v$.
Note that for $u\in \Om^0(E)$ and $v\in \Om^0(F)$ we have
\begin{align}
\cov^2_{A\otimes B}(u\otimes v) 
&= (\cov_A^2 u)\otimes v + 2\cov_A u\otimes\cov_B v + u\otimes\cov _B^2v, 
\label{eq:SecondCovDerivLeibnitz}\\
\cov^*_{A\otimes B}\cov_{A\otimes B}(u\otimes v) 
&= (\cov_A^*\cov_A u)\otimes v + *((*\cov_Au)\wedge\cov_B v) \notag\\
&\qquad - *(\cov_A u\wedge *\cov_B v) + u\otimes\cov_B^*\cov_B v. \notag
\end{align}
Similarly, for $u\in\Om^0(\La^1\otimes E)$ and $v\in\Om^0(F)$, we have
\begin{equation}
\cov^*_{A\otimes B}(u\otimes v) 
= (\cov_A^*u)\otimes v + *(*u\wedge \cov_B v)
\label{eq:AdjointCovLeibnitz}
\end{equation}
In particular, we see that if $u,v\in\Om^0(\gl(E))$, then
\begin{align}
\cov^*_A\cov_A(uv) 
&= (\cov_A^*\cov_A u)v + *((*\cov_Au)\wedge\cov_A v) \\
&\qquad - *(\cov_A u\wedge (*\cov_A v)) + u(\cov_A^*\cov_A v), \notag
\label{eq:LaplaceLeibnitz}
\end{align}
an identity we will need in the next section.

\begin{lem}\label{lem:*Multiplication}
Let $\Om^0(E)\times\Om^0(F)\to \Om^0(E\otimes F)$ be given by $(u,v)\mapsto
u\otimes v$. Then the following hold.
\begin{enumerate}
\item The map $C^0(E)\otimes L^\sharp(F)\to
L^\sharp(E\otimes F)$ is continuous; 
\item The map $L^{2\sharp}(E)\otimes 
L^{2\sharp}(F) \to L^\sharp(E\otimes F)$ is continuous; 
\item The spaces $L^\sharp_1(F)$, $L^2_1(F)$, and
$L^\sharp_2(F)$ are $L^\sharp_2(E)$-modules;
\item The spaces $L^2_1(F)$, $L^{\sharp,2}_1(F)$, and
$L^{\sharp,2}_2(F)$ are $L^{\sharp,2}_2(E)$-modules;
\end{enumerate}
The conclusions continue to hold for $\Om^1(E)$ in place of $\Om^0(E)$ and
the norms on $L^{\sharp}_1(\La^1\otimes E)$ and
$L^{\sharp,2}_1(\La^1\otimes E)$ defined via 
\eqref{eq:LSharp1AdjointCov}. 
\end{lem}

\begin{pf}
Let $u\in C^\8(E)$ and $v\in C^\8(F)$ and denote the covariant derivatives
on $E$, $F$, and $E\otimes F$ by $\cov$. Using $\cov(u\otimes v)=(\cov
u)v + u\otimes\cov v$ and the embedding
$L^\sharp_2(E)\subset C^0(E)$, we see that
\begin{align*}
\|u\otimes v\|_{L^\sharp} &\le \|u\|_{C^0}\|v\|_{L^\sharp} \quad\text{and}\quad
\|u\otimes v\|_{L^\sharp} \le \|u\|_{L^{2\sharp}}\|v\|_{L^{2\sharp}}, \\
\|\cov(u\otimes v)\|_{L^\sharp} &\le \|\cov u\|_{L^{2\sharp}}\|v\|_{L^{2\sharp}} 
+ \|u\|_{C^0}\|\cov v\|_{L^\sharp} 
\le c\|u\|_{L^\sharp_2}\|v\|_{L^\sharp_1}, \\
\|\cov(u\otimes v)\|_{L^2} &\le \|\cov u\|_{L^4}\|v\|_{L^4}
+ \|u\|_{C^0}\|\cov v\|_{L^2} \le \|u\|_{C^0\cap L^2_2}\|v\|_{L^2_1},
\end{align*}
and hence the multiplication maps $C^0\times L^\sharp\to L^\sharp$,
$L^{2\sharp}\times L^{2\sharp}
\to L^\sharp$, and $L^\sharp_2\times L^\sharp_1 \to L^\sharp_1$ are
continuous. Moreover, 
$$
\|\cov(u\otimes v)\|_{L^2} \le \|\cov u\|_{L^4}\|v\|_{L^4} 
+ \|u\|_{C^0}\|\cov v\|_{L^2} 
\le c\|u\|_{C^0\cap L^2_2}\|v\|_{L^2_1},
$$
and so, using the embedding $L^\sharp_2\subset C^0$, the multiplication
$L^{\sharp,2}_2\times L^{\sharp,2}_1 \to L^{\sharp,2}_1$ is
continuous. Thus, $L^\sharp_2$ is an $L^\sharp_1$-module and
$L^{\sharp,2}_2$ is an $L^{\sharp,2}_1$-module.

Finally, to see that $L^\sharp_2$ and $L^{\sharp,2}_2$ are algebras, we use
the identities \eqref{eq:SecondCovDerivLeibnitz}, noting that
\begin{align*}
\|\cov^*\cov(u\otimes v)\|_{L^\sharp} 
&\le \|\cov^*\cov u\|_{L^\sharp}\|v\|_{C^0} 
+ 2\|\cov u\|_{L^{2\sharp}}\|\cov v\|_{L^{2\sharp}} 
+ \|u\|_{C^0}\|\cov^*\cov v\|_{L^\sharp} \\
&\le c\|u\|_{L^\sharp_2}\|v\|_{L^\sharp_2},
\end{align*}
so the multiplication
$L^\sharp_2\times L^\sharp_2 \to L^\sharp_2$ is continuous, while
\begin{align*}
\|\cov^2(u\otimes v)\|_{L^2} &\le \|\cov^2u\|_{L^2}\|v\|_{C^0} 
+ 2\|\cov u\|_{L^4}\|\cov v\|_{L^4} + \|u\|_{C^0}\|\cov^2 v\|_{L^2} \\
&\le c\|u\|_{C^0\cap L^2_2}\|v\|_{C^0\cap L^2_2}.
\end{align*}
The embedding $L^\sharp_2\subset C^0$ now implies that the multiplication
$L^{\sharp,2}_2\times L^{\sharp,2}_2 \to L^{\sharp,2}_2$ is continuous. 
\end{pf}

\section{Critical-exponent Sobolev embeddings and estimates for Green's
operators}  
\label{sec:Green}
We continue the notation and assumptions of Section \ref{sec:SharpSobolev}.
Our goal in this section is to prove the Sobolev embedding Lemmas
\ref{lem:EmbeddingSobolevInto*} and \ref{lem:Embedding*IntoSobolev}, and to
derive estimates for the Green's operator $G_A$ of the Laplacian
$\cov_A^*\cov_A$ on $\Om^0(E)$. The key estimates described in this section
are due to Taubes and they arise, in a variety of contexts, in
\cite{ParkerTaubes,TauPath,TauFrame,TauStable,TauConf}. However, we find it
convenient to collect them here --- together with some useful extensions
and generalizations --- both for the purposes of the present
article and applications in \cite{Daskal,FeehanLeness,MorganOzsvath}.

\subsection{Estimates for the covariant Laplacian $\cov_A^*\cov_A$}
\label{subsec:EstCovLaplacian}
Let $G(x,y)$ be the kernel function for
the Green's operator $(d^*d+1)^{-1}$ of the Laplacian $d^*d+1$ on
$C^\8(X)$. The kernel $G(x,y)$ of $(d^*d+1)^{-1}$ 
behaves like $\dist^{-2}(x,y)$ as $\dist(x,y)\to 0$ (see
\cite[Lemma 4.7]{TauConf} and \cite[\S 5]{TauStable}):

\begin{lem}\label{lem:Green}
The kernel $G(x,y)$ is a positive $C^\8$ function away from the diagonal
in $X\times X$ and as $\dist(x,y)\to 0$, 
$$
G(x,y) = {\frac{1}{4\pi^2\dist^2(x,y)}} + o(\dist^{-2}(x,y)).
$$
\end{lem}

\begin{pf}
These and other properties of $G$ are obtained by explicitly constructing
$G$ from an initial choice of parametrix $H$ for $d^*d+1$ using the method of
\cite[\S 4.2.2--3]{Aubin}, where the kernel for the Green's operator for
$d^*d$ is constructed. Recall from \cite[p. 132]{Stein} that the
kernel $G_0(x,y)$ for $(d^*d+1)^{-1}$ on $\RR^4$ with its standard metric
satisfies 
$$
G_0(x,y) = {\frac{1}{4\pi^2|x-y|^2}} + o(|x-y|^{-2}),\qquad |x-y|\to 0.
$$
The kernel $G$ is now constructed using $G_0$ by following the method
of \cite[\S 4.2.2--3]{Aubin}.
\end{pf}

Lemma \ref{lem:Green} implies that there is a constant $c$ depending
at most on $g$ such that for all $x\ne y$ in $X$,
\begin{equation}
c^{-1}\dist^{-2}(x,y) \le G(x,y) 
\le c\dist^{-2}(x,y). \label{eq:GreenKerEst}
\end{equation}
Consequently, for all $u\in \Om^0(E)$, we have
\begin{equation}
c^{-1}\|u\|_{L^\sharp(X)} \le \|G|u|\|_{C^0(X)} \le c\|u\|_{L^\sharp(X)}.
\label{eq:LinftyGreenEst}
\end{equation}
Lemma \ref{lem:Embedding*IntoSobolev} will follow from the next estimate;
a similar inequality is stated as Equation (3.4) in 
\cite{TauStable} for the case $p=2$; see \cite[Lemma 5.4]{ParkerTaubes} for
the proof when $p=2$ and $X=\RR^4$.

\begin{lem}\label{lem:Lp*Est}
Let $X$ be a closed, oriented four-manifold with metric $g$ and let $1\le
p<4$. Then there is a positive constant $C$ such that for all $f\in
L^p_1(X)$,
$$
\sup_{x\in X}\|\dist^{-1}(x,\cdot)f\|_{L^p} \le C\|f\|_{L^p_1}.
$$
\end{lem}

\begin{pf}
We may assume that $f\in C^\8(X)$.
Let $\varrho$ denote the injectivity radius of the Riemannian manifold
$(X,g)$ and let $\{B(x_i,\varrho)\}_{i=1}^N$ be a covering of $X$ by $N$ geodesic
balls with centers at points $x_i$ and radius $\varrho/2$. Let
$\{\chi_i\}_{i=1}^N$ be a partition of unity subordinate to this cover, so
that $\supp\chi_i\subset B(x_i,\varrho)$ and $\sum_{i=1}^N\chi_i=1$.
Using $f=\sum_{i=1}^N\chi_i f$, we have
$$
\sup_{x\in X}\|\dist^{-1}(x,\cdot)f\|_{L^p}^p 
\le \sum_{i=1}^N\sup_{x\in B(x_i,\varrho)}
\int_{B(x_i,\varrho)}\dist^{-p}(x,\cdot)|\chi_i f|^p\,dV.
$$
For any fixed $x\in X$, let $\exp_x^{-1}:B(x,\varrho)\to\RR^4$ be a
geodesic normal coordinate chart centered at the point $x\in X$, let
$\{y^\mu\}$ be the induced local coordinates centered at $x\in X$,
and let $r=\dist(x,y)$. 
We apply the divergence theorem and integration by parts to bound the 
integrals
$$
\int_{B(x_i,\varrho)}\dist^{-p}(x,\cdot)|\chi_i f|^p\,dV.
$$
For this purpose define vector fields $\eta_i=\sum_\mu\eta_i^\mu\rd/\rd
y^\mu$ supported in $B(x_i,\varrho)$ by setting
$$
\eta_i^\mu = y^\mu r^{-p}|\chi_i f|^p,\qquad i=1,\dots,N.
$$
If $g=\det(g_{\mu\nu})$ then, since $r^2=\sum_\mu(y^\mu)^2$, we have
$\sum_\mu y^\mu\rd/\rd y^\mu = r\rd/\rd r$ and so
\begin{align*}
\divg_g\eta_i &= {\frac{1}{\sqrt{g}}}\sum_\mu{\frac{\rd}{\rd y^\mu}}
\left(\eta_i^\mu \sqrt{g}\right) \\
&= \divg_{g_e}\eta_i 
+ \sum_\mu\eta_i^\mu{\frac{\rd}{\rd y^\mu}}\log\sqrt{g} \\
&= \divg_{g_e}\eta_i + r^{1-p}|\chi_i f|^p{\frac{\rd}{\rd r}}\log\sqrt{g},
\end{align*}
where $\divg_{g_e}\eta_i = \sum_\mu\rd\eta_i^\mu/\rd y^\mu$ is the
divergence of $\eta_i$ with respect to the pullback of the Euclidean metric
$g_e$ on $B(x_i,\varrho)$. From \cite[Theorem 1.53]{Aubin} we recall that
$\rd(\log\sqrt{g})/\rd r = O(r)$ and thus
$$
\divg_g\eta_i = \divg_{g_e}\eta_i + O(r^{2-p}|\chi_i f|^p).
$$
Using $\sum_\mu y^\mu\rd/\rd y^\mu = r\rd/\rd r$, we have
\begin{align*}
\divg_{g_e}\eta_i 
&= \sum_\mu{\frac{\rd}{\rd y^\mu}}\left(y^\mu r^{-p}|\chi_i f|^p\right) \\
&= (4-p) r^{-p}|\chi_i f|^p + pr^{1-p}|\chi_i f|^{p-1}
{\frac{\rd |\chi_if|}{\rd r}},
\end{align*}
and so as $p<4$,
\begin{align*}
r^{-p}|\chi_i f|^p &= {\frac{1}{4-p}}\divg_{g_e}\eta_i 
- {\frac{p}{4-p}}r^{1-p}|\chi_i f|^{p-1}{\frac{\rd|\chi_i f|}{\rd r}} \\
&\le {\frac{1}{4-p}}\divg_g\eta_i + {\frac{c}{4-p}}r^{2-p}|\chi_i f|^p
+ {\frac{p}{4-p}}r^{1-p}|\chi_i f|^{p-1}{\frac{\rd|\chi_i f|}{\rd r}}.
\end{align*}
Note that $\rd/\rd r$ is a unit vector with respect to the metric
$g$. The divergence theorem, H\"older's inequality, and the fact that $\eta_i$
has compact support in $B(x_i,\varrho)$ imply that
\begin{align*}
&\int_{B(x_i,\varrho)}r^{-p}|\chi_i f|^p\,dV \\
&\le {\frac{p}{4-p}}\int_{B(x_i,\varrho)}r^{1-p}
|\chi_i f|^{p-1}|d(\chi_i f)|\,dV 
+ {\frac{c}{4-p}}\int_{B(x_i,\varrho)}r^{2-p}|\chi_i f|^p\,dV \\
&\le {\frac{p}{4-p}}\left(\int_{B(x_i,\varrho)}
r^{-p}|\chi_i f|^p\,dV\right)^{(p-1)/p}
\left(\int_{B(x_i,\varrho)}|d(\chi_i f)|^p\,dV\right)^{1/p} \\
&\qquad + {\frac{c}{4-p}}
\left(\int_{B(x_i,\varrho)}r^{-p}|\chi_i f|^p\,dV\right)^{(p-1)/p}
\left(\int_{B(x_i,\varrho)}r^p|\chi_i f|^p\,dV\right)^{1/p},
\end{align*}
Consequently, 
$$
\|\dist^{-1}(x,\cdot)\chi_i f\|_{L^p}
\le C\left(\|d(\chi_i f)\|_{L^p} + \|\chi_i f\|_{L^p}\right)
\le C\left(\|d f\|_{L^p} + \|f\|_{L^p}\right),
$$
for each $i$ and so the desired result follows.
\end{pf}

\begin{proof}[Proof of Lemma~\ref{lem:EmbeddingSobolevInto*}]
Define $1\le p'<2$ by setting $1=1/p+1/p'$. Then
H\"older's inequality implies that
\begin{align*}
\|\dist^{-2}(x,\cdot)|u|\|_{L^1}
&\le \|\dist^{-2}(x,\cdot)\|_{L^{p'}}\|u\|_{L^p} \le C\|u\|_{L^p}, \\
\|\dist^{-1}(x,\cdot)|u|\|_{L^2}
&\le \|\dist^{-1}(x,\cdot)\|_{L^{2p'}}\|u\|_{L^{2p}} 
\le C\|u\|_{L^{2p}}, 
\end{align*}
which gives Assertions (1) and (2).
By Lemma \ref{lem:Lp*Est} and Kato's inequality, $|d|u||\le
|\cov_Au|$, we see that  
\begin{align*}
\sup_{x\in X}\|\dist^{-1}(x,\cdot)u\|_{L^p} 
&= \sup_{x\in X}\|\dist^{-1}(x,\cdot)|u|\|_{L^p}  \\
&\le C(\|d|u|\|_{L^p} + \|u\|_{L^p}) \\
&\le C(\|\cov_Au\|_{L^p} + \|u\|_{L^p}).
\end{align*}
Taking $p=2$ gives Assertion (3).
\end{proof}

Lemma \ref{lem:Embedding*IntoSobolev} will follow from the estimates below;
the key estimates (1) and (2) in Lemma \ref{lem:LinftyEstu} below and the
estimates (1), (2), and (3) in Lemma \ref{lem:L22Estu} are essentially
those of Lemma 6.2 in
\cite{TauFrame}, except that the dependence of the constant on
$\|F_A\|_{L^2}$ is made explicit, but the argument is the same as that of
\cite{TauFrame}.  

\begin{lem}\label{lem:LinftyEstu}
Let $X$ be a closed, oriented four-manifold with metric $g$. Then there is
a constant $c$ with the following significance. Let $E$ be a Riemannian
vector bundle over $X$ and let $A$ be an orthogonal
$L^2_2$ connection on $E$ with curvature $F_A$. Then  
$L^\sharp_2(E)\subset C^0\cap L^2_1(E)$ and the following estimates hold: 
\begin{align}
\|\cov_Au\|_{L^{2\sharp}(X)}+ \|u\|_{C^0(X)} 
&\le c\|\cov_A^*\cov_Au\|_{L^\sharp(X)} + \|u\|_{L^\sharp(X)}, \tag{1}\\
\|\cov_Au\|_{L^{2\sharp}(X)}+ \|u\|_{C^0(X)} 
&\le c\|\cov_A^*\cov_Au\|_{L^\sharp(X)} + \|u\|_{L^2(X)}, \tag{2}\\
\|u\|_{L^1(X)} &\le c\|u\|_{L^{\sharp}(X)}, \tag{3} \\
\|u\|_{L^2(X)} &\le c\|u\|_{L^{2\sharp}(X)}, \tag{4} \\
\|\cov_Au\|_{L^2(X)} &\le c\|\cov_Au\|_{L^{2\sharp}(X)}. \tag{5}
\end{align}
\end{lem}

\begin{pf}
For any $u\in C^\8(E)$ there is the following pointwise identity
\cite[p. 93]{FU},
$$
|\cov_Au|^2 + \half d^*d|u|^2 
= \langle\cov_A^*\cov_A u,u\rangle,
$$
and thus:
$$
|\cov_Au|^2 + \half(1+ d^*d)|u|^2 
= \langle\cov_A^*\cov_A u,u\rangle + \half|u|^2,
$$
Using this identity and the fact that 
$\int_X G(x,\cdot)(d^*d+1)|u|^2\,dV = |u|^2(x)$, we obtain
$$
\int_X G(x,\cdot)|\cov_Au|^2\,dV + \half |u|^2(x)
\le \int_X G(x,\cdot)|\langle \cov_A^*\cov_Au,u\rangle|\,dV
+ \half\int_X G(x,\cdot)|u|^2\,dV.
$$
Therefore, from \eqref{eq:GreenKerEst}, we have
\begin{align*}
\||\cov_Au|^2\|_{L^\sharp} + \||u|^2\|_{C^0}
&\le c\|\langle \cov_A^*\cov_Au,u\rangle\|_{L^\sharp} +
c\||u|^2\|_{L^\sharp} \\ 
&\le c\|\cov_A^*\cov_Au\|_{L^\sharp}\|u\|_{C^0} + c\|u\|_{L^\sharp}\|u\|_{C^0} 
\end{align*}
Consequently, using rearrangement with the last term, we see that
$$
\|\cov_Au\|_{L^{2\sharp}} +
\|u\|_{C^0} \le c\|\cov_A^*\cov_Au\|_{L^\sharp} + c\|u\|_{L^\sharp},
$$
giving (1).
Combining this estimate with the embedding and interpolation
inequalities,
$\|u\|_{L^\sharp} \le c\|u\|_{L^4}\le c\|u\|_{L^2}^{1/2}\|u\|_{C^0}^{1/2}$,
and again using rearrangement with the last term yields the bound in
(2). Since $X$ is closed, for all $x\ne y$ we have $\dist(x,y)\le M<\8$, so
$$
\int_X\dist^{-2}(x,\dot)|u|\,dV \ge M^{-2}\int_X|u|\,dV,
$$
and this gives the estimates in (3), (4), and (5).
\end{pf}

\begin{proof}[Proof of Lemma~\ref{lem:Embedding*IntoSobolev}]
{}From Lemma \ref{lem:LinftyEstu} we have the estimate 
$$
\|u\|_{C^0} \le c\|u\|_{L^\sharp_{2,A_0}},
$$
for any $u\in C^\8(E)$. Let $\{u_m\}$ be a sequence in $C^\8(E)$ converging
to $u\in L^\sharp_2(E)$. The sequence $\{u_m\}$ is Cauchy in 
$L^\sharp_2(E)$ and
applying the preceding estimate to the differences $u_{m_2}-u_{m_1}$, we
see that it is Cauchy in the Banach space $C^0(E)$ 
and so the limit $u$ lies in $C^0(E)$. The same argument, with estimates
(1) and (5) of Lemma \ref{lem:LinftyEstu}, shows that $u\in L^2_1(E)$ and
this yields Assertion (2) of the lemma. Assertion (1) follows in the same
manner. 
\end{proof}

\begin{lem}\label{lem:L22Estu}
Continue the hypotheses of Lemma \ref{lem:LinftyEstu}. Then for any $u\in
(C^0\cap L^2_2)(E)$, we have
\begin{align}
\|\cov_A^2u\|_{L^2(X)} &\le \|\cov_A^*\cov_Au\|_{L^2(X)} +
c\|F_A\|_{L^2(X)}^{1/2}\|\cov_Au\|_{L^4(X)} \tag{1}\\
&\qquad + \|F_A\|_{L^2(X)}\|u\|_{C^0(X)}, \notag\\
\|\cov_Au\|_{L^4(X)} &\le \|u\|_{C^0(X)}^{1/2}
\left(\|\cov_A^*\cov_Au\|_{L^2(X)}+2\|\cov_A^2u\|_{L^2(X)}\right)^{1/2}, 
\tag{2}\\
\|\cov_A^2u\|_{L^2(X)}
&\le 2\|\cov_A^*\cov_Au\|_{L^2(X)} + c\|F_A\|_{L^2(X)}\|u\|_{C^0(X)}. \tag{3}
\end{align}
\end{lem}

\begin{pf}
The Bochner-Weitzenb\"ock formula for the covariant Laplacian 
\cite[Appendix, Theorem II.1]{Lawson} asserts that
\begin{equation}
d_A^*d_A + d_Ad_A^* = \cov_A^*\cov_A + \{F_A,\cdot\},
\end{equation}
where we use $\{\cdot,\cdot\}$ to denote a certain bilinear map whose
precise form is unimportant here. Integrating by parts and noting that $d_A =
\cov_A$ and $d_A^*d_A = \cov_A^*\cov_A$ on $\Om^0(X,V)$ and
$F_A = d_A\circ d_A$ gives
\begin{align*}
\|\cov_A^2u\|_{L^2}^2 &= (\cov_A^*\cov_A\cov_Au,\cov_Au)_{L^2} \\
&= ((d_A^*d_A+d_Ad_A^*)d_Au,d_Au)_{L^2} - (\{F_A,d_Au\},\cov_Au)_{L^2} \\
&= (d_A^*F_Au,d_Au)_{L^2} + (d_A(d_A^*d_A)u,d_Au)_{L^2} 
- (\{F_A,\cov_Au\},\cov_Au)_{L^2} \\
&= (F_Au,F_Au)_{L^2} + (\cov_A^*\cov_Au,\cov_A^*\cov_Au)_{L^2} 
- (\{F_A,\cov_Au\},\cov_Au)_{L^2}.
\end{align*}
Therefore, applying H\"older's ineqality, we find that
$$
\|\cov_A^2u\|_{L^2}^2 \le \|\cov_A^*\cov_Au\|_{L^2}^2 
+ c\|F_A\|_{L^2}\|\cov_Au\|_{L^4}^2
+ \|F_A\|_{L^2}^2\|u\|_{C^0}^2,
$$
and taking square roots gives the desired bound in (1).

We now use integration by parts and Kato's inequality $|d|u||\le |\cov_Au|$
to obtain an $L^4$ bound on $d_Au$:
\begin{align*}
\|d_Au\|_{L^4}^4
&= \left(d_Au,|d_Au|^2d_Au\right)_{L^2} \\
&= \left(u,|d_Au|^2d_A^*d_Au\right)_{L^2} + 
2\left(u,|d_Au|d_Au\wedge d|d_Au|\right)_{L^2}, \\
&= \|u\|_{C^0}\|d_Au\|_{L^4}^2\|d_A^*d_Au\|_{L^2} +
2\|u\|_{C^0}\|d_Au\|_{L^4}^2\|\cov_Ad_Au\|_{L^2},
\end{align*}
and so, if $d_Au\not\equiv 0$,
$$
\|d_Au\|_{L^4}\le 
\|u\|_{C^0}^{1/2}\left(\|d_A^*d_Au\|_{L^2}+2\|\cov_A^2u\|_{L^2}\right)^{1/2},
$$
which gives the desired estimate in (2).

By combining the $L^4$ estimate for $\cov_Au$ with the $L^2$ estimate for
$\cov_A^2u$, we obtain
\begin{align*}
\|\cov_A^2u\|_{L^2} &\le \|\cov_A^*\cov_Au\|_{L^2} + \|F_A\|_{L^2}\|u\|_{C^0} \\
&\qquad + c\|F_A\|_{L^2}^{1/2}\|u\|_{C^0}^{1/2}
\left(\|\cov_A^*\cov_Au\|_{L^2}+\|\cov_A^2u\|_{L^2}\right)^{1/2}.
\end{align*}
We now use rearrangement with the last term above to give
$$
\|\cov_A^2u\|_{L^2}\le 2\|\cov_A^*\cov_Au\|_{L^2} + c\|F_A\|_{L^2}\|u\|_{C^0},
$$
and this establishes the desired bound in (3).
\end{pf}

\begin{lem}\label{lem:L22InfinityEstu}
Continue the hypotheses of Lemma \ref{lem:LinftyEstu}. Then for any $u\in
L^{\sharp,2}_2(E)$, we have:
$$
\|u\|_{L^2_{2,A}(X)} + \|u\|_{C^0(X)}  
\le c(1+\|F_A\|_{L^2(X)})\left(\|\cov_A^*\cov_Au\|_{L^{\sharp,2}(X)} +
\|u\|_{L^2(X)}\right). 
$$
\end{lem}

\begin{pf}
{}From Assertion (3) of Lemm \ref{lem:L22Estu} we have the estimate
$$
\|\cov_A^2u\|_{L^2} \le 2\|\cov_A^*\cov_Au\|_{L^2} +
c\|F_A\|_{L^2}\|u\|_{C^0}, 
$$
while integration by parts gives
$$
\|\cov_Au\|_{L^2} = (\cov_A^*\cov_Au,u)_{L^2}^{1/2} \le 
{\frac{1}{\sqrt{2}}}\left(\|\cov_A^*\cov_Au\|_{L^2} +\|u\|_{L^2}\right).
$$
According to Lemma \ref{lem:LinftyEstu} we have  
$$
\|u\|_{C^0} \le c\|\cov_A^*\cov_Au\|_{L^\sharp} + c\|u\|_{L^2}, 
$$
and therefore the desired bound follows by combining these estimates.
\end{pf}

The above lemmas lead to the following estimates for the Green's operator
$G_A:L^{\sharp,2}(E)\to L^{\sharp,2}_2(E)$ of the Laplacian
$\cov_A^*\cov_A:L^{\sharp,2}_2(E)\to L^{\sharp,2}(E)$. For $u\in\Om^0(E)$
define  
\begin{equation}
\|u\|_{\sL^{\sharp,2}_{2,A}(X)} 
=\|\cov_A^*\cov_Au\|_{L^{\sharp,2}(X)} + \|u\|_{L^{\sharp,2}(X)},
\end{equation}
and observe that this is equivalent to the $L^{\sharp,2}_{2,A}$ norm
defined in Section \ref{sec:SharpSobolev}, although the comparison depends
on the $L^2$ norm of the curvature $F_A$.

\begin{lem}\label{lem:L2*2GreenEst}
Continue the hypotheses of Lemma \ref{lem:LinftyEstu}. Let $\nu_0[A]$
be the least positive eigenvalue of the Laplacian 
$\cov_A^*\cov_A$. Then for any
$u\in L^{\sharp,2}\cap (\Ker \cov_A^*\cov_A)^\perp$, we have:
\begin{align}
\|G_Au\|_{L^\sharp_{2,A}(X)} 
&\le c(1+\nu_0[A]^{-1})\|u\|_{L^{\sharp,2}(X)},  \tag{1} \\
\|G_Au\|_{\sL^{\sharp,2}_{2,A}(X)} 
&\le c(1+\nu_0[A]^{-1})\|u\|_{L^{\sharp,2}(X)},  \tag{2} \\
\|G_Au\|_{L^{\sharp,2}_{2,A}(X)} 
&\le c(1+\nu_0[A]^{-1})(1+\|F_A\|_{L^2(X)})\|u\|_{L^{\sharp,2}(X)}.  \tag{3}
\end{align}
\end{lem}

\begin{pf}
The first and second
assertions follow from Lemma \ref{lem:LinftyEstu}, the fact that
$\cov_A^*\cov_A G_Au = u$ for $u\in (\Ker\cov_A^*\cov_A)^\perp$, and the
eigenvalue estimate $\|u\|_{L^2} 
\le \nu_0[A]^{-1}\|\cov_A^*\cov_Au\|_{L^2}$, while 
the third assertion follows from the first and Lemma
\ref{lem:L22InfinityEstu}.
\end{pf}

\subsection{Elliptic estimates for $d_A^++d_A^*$}
\label{subsec:EllipticEstimatedA+}
To illustrate their application and to point to possible extensions, we
note that the estimates of Section \ref{subsec:EstCovLaplacian} for the
covariant Laplacian $\cov_A^*\cov_A = d_A^*d_A$ on $\Omega^0(E)$ naturally
extend to give estimates for the covariant Laplacians $\cov_A^*\cov_A$ on
$\Omega^\ell(E)=\Om^0(\La^\ell\otimes E)$. Estimates for $\cov_A^*\cov_A$
on $\Om^0(\La^1\otimes E)$ and $\Om^0(\La^+\otimes E)$ are of particular
interest since these can in turn be profitably compared (via the
Bochner-Weitzenb\"ock formulas \cite[Equations (6.25) \& (6.26)]{FU}, as in
\cite{TauStable}) with the remaining Laplacians defined by the elliptic
deformation complex for the anti-self-dual equation \cite{DK,FU}, namely
$d_Ad_A^* + d_A^{+,*}d_A^+$ on $\Om^0(\La^1\otimes \fg_E)$ and
$d_A^+d_A^{+,*}$ on $\Om^0(\La^+\otimes \fg_E)$. Indeed, if $B_1$ and
$B_+$ are the Levi-Civita connections on $\La^1$ and $\La^+$ induced by the
Levi-Civita connection on $TX$ for the metric $g$, then the curvature
`$F_A$' in the estimates of the preceding subsection is simply replaced by
\cite[p. 165]{LM}
\begin{equation}
\label{eq:TensorConnection}
\begin{aligned}
F_{B_1\otimes A} &= F_{B_1}\otimes\id_{\fg_E} + \id_{\La^1}\otimes F_A, 
\\
F_{B_+\otimes A} &= F_{B_+}\otimes\id_{\fg_E} + \id_{\La^+}\otimes F_A,
\end{aligned}
\end{equation}
where $F_{B_1}$ and $F_{B_+}$ are expressed in terms of the Riemann
curvature tensor $\Rm$ and where we abuse notation slightly and denote the
connections on $E$ and $\fg_E$ both by $A$. (See \cite[Appendix C]{FU} and
\cite[Appendix II]{Lawson}.) In the interests of brevity we shall confine
our attention to the case of $L^p_\ell$ estimates with $p=2$, though the
methods can be modified to obtain estimates for $p\neq 2$ (some work is
required --- see \cite[p. 426]{DK} for hints). 

In order to compute the required elliptic estimates for $d_A^+$ we will
need the Bochner-Weitzenb\"ock formulas,
\begin{align}
d_Ad_A^* + 2d_A^*d_A^+ 
&= \cov_A^*\cov_A + \{\Ric,\cdot\}
- 2\{F_A^-,\cdot\}, \label{eq:BW1} \\
2d_A^+d_A^* &= \cov_A^*\cov_A - 2\{W^+,\cdot\} + \frac{R}{3} +
\{F_A^+,\cdot\}, \label{eq:BW+}
\end{align}
for the Laplacians on $\Om^1(\fg_E)$ and $\Om^+(\fg_E)$
\cite[Equations (6.25) \& (6.26)]{FU}; here, $\Ric$, $W^+$, and $R$ are the
Ricci, self-dual Weyl, and scalar curvatures of the Riemannian metric $g$ on
$X$.  In applications to the degeneration of anti-self-dual or `almost
anti-self-dual' connections $A$ as in
\cite{TauSelfDual,TauIndef,TauFrame,TauStable}, we can usually arrange to
have a uniform $L^\8$ bound on $F_A^+$, but not a uniform $L^p$ bound on
$F_A^-$ when $p>2$. We derive estimates in the remainder of this subsection
with such applications and assumptions in view.  To illustrate the nature
of the difficulty we first derive a naive $L^2_{1,A}$ estimate for $a\in
L^2_1(\fg_E)$ in terms of the $L^2$ norm of $(d_A^* + d_A^+)a$:

\begin{lem}\label{lem:L21AEsta}
Let $X$ be a closed, oriented four-manifold with metric $g$. Then there is
a constant $c$ with the following significance. Let $E$ be a Riemannian
vector bundle over $X$ and let $A$ be an orthogonal $L^2_4$ connection on
$E$ with curvature $F_A$. Then for any $a\in L^2_1(\Lambda^1\otimes\fg_E)$,
\begin{equation}
\label{eq:NaiveL21dA+Estimate}
\|a\|_{L^2_{1,A}(X)} 
\le \sqrt{2}\|(d_A^* + d_A^+)a\|_{L^2(X)}
+ c\left(1 + \|F_A^-\|_{C^0(X)}\right)^{1/2}\|a\|_{L^2(X)}. 
\end{equation}
If $a$ is $L^2$-orthogonal to $\Ker d_A^+$, so that $a=d_A^*v$ for some
$v\in L^2_2(\Lambda^+\otimes\fg_E)$, then
\begin{equation}
\label{eq:NaiveL21dA+OrthogKernelEstimate}
\begin{aligned}
\|d_A^*v\|_{L^2_{1,A}(X)} 
&\le \sqrt{2}\|d_A^+d_A^*v\|_{L^2(X)} 
+ c\left(1 + \|F_A^-\|_{C^0(X)}\right)^{1/2}\|v\|_{L^2(X)}
\\
&\quad + \|F_A^+\|_{C^0(X)}\|v\|_{L^2(X)}.
\end{aligned}
\end{equation}
\end{lem}

\begin{pf}
From the Bochner-Weitzenb\"ock formula for 
$d_Ad_A^*+2d_A^*d_A^+$ in
\eqref{eq:BW1} and integration by parts, we have:
\begin{align*}
\|\cov_Aa\|_{L^2}^2 &= \left(\cov_A^*\cov_Aa,a\right) \\
&=\left(d_Ad_A^*a,a\right) + 2\left(d_A^*d_A^+a,a\right)
- \left(\{\Ric,a\},a\right) 
+ 2\left(\{F_A^-,a\},a\right) \\
&\le \|d_A^*a\|_{L^2}^2+2\|d_A^+a\|_{L^2}^2
+c\left(1+\|F_A^-\|_{C^0}\right)\|a\|_{L^2}^2   
\end{align*}
which gives \eqref{eq:NaiveL21dA+Estimate}. If $a=d_A^*v$, then
$d_A^*d_A^*v = (d_A^+d_A)^*v = (F_A^+)^*v$, so that
$$
\|d_A^*d_A^*v\|_{L^2}\le \|F_A^+\|_{C^0}\|v\|_{L^2}.
$$
Thus, \eqref{eq:NaiveL21dA+OrthogKernelEstimate} 
follows from \eqref{eq:NaiveL21dA+Estimate} and the above inequality.
\end{pf}

Since $d_A^*+d_A^+$ is an elliptic operator, estimates of the above form
follow from the general theory of linear elliptic operators.  However, the
preceding elementary derivation using the Bochner-Weitzenb\"ock formula
gives us a constant whose dependence on the curvature terms $F_A^-$ and
$F_A^+$ is made explicit. In particular, we see that the estimate is only
useful when we have a uniform $C^0$ bound on $F_A^-$ independent of $A$,
which is not possible when $A$ bubbles. At the cost of introducing a
slightly stronger norm than the $L^2$ norm on the right hand side of the
estimate above, we can derive an $L^2_{1,A}$ bound for $a=d_A^*v$ with an
estimate constant depending on $\|F_A^-\|_{L^2(X)}$ rather than
$\|F_A^-\|_{C^0(X)}$. Specifically, Equation \eqref{eq:TensorConnection}
and Lemma \ref{lem:L22InfinityEstu} give the following $L^2_{2,A}$
estimates for sections of $\La^+\otimes\fg_E$:

\begin{lem}\label{lem:LinftyL22CovLapEstv}
Continue the hypotheses of Lemma \ref{lem:L21AEsta}. Then the following
estimate holds for any $v\in L^{\sharp,2}_2(\Lambda^+\otimes\fg_E)$: 
\begin{align*}
\|v\|_{L^2_{2,A}(X)} + \|v\|_{C^0(X)} 
&\le c(1+\|F_A\|_{L^2(X)})
(\|\cov_A^*\cov_Av\|_{L^{\sharp,2}(X)} + \|v\|_{L^2(X)}). 
\end{align*}
\end{lem}

We now replace the covariant Laplacian $\cov_A^*\cov_A$ in
the estimates of Lemma \ref{lem:LinftyL22CovLapEstv} by the
Laplacian $d_A^+d_A^*$ via the Bochner formula \eqref{eq:BW+} to
give: 

\begin{lem}\label{lem:LinftyL22Estv}
Continue the hypotheses of Lemma \ref{lem:LinftyL22CovLapEstv}. Then there
is a positive constant $\eps=\eps(c)$ such that the following holds. If
$\|F_A\|_{L^{\sharp,2}(X)} <\eps$, then
\begin{align*}
\|v\|_{L^2_{2,A}(X)} + \|v\|_{C^0(X)}
&\le c(1+\|F_A\|_{L^2(X)})
(\|d_A^+d_A^*v\|_{L^{\sharp,2}(X)} + \|v\|_{L^2(X)}).
\end{align*}
\end{lem}

\begin{proof}
{}From \eqref{eq:TensorConnection} and Lemma \ref{lem:*Multiplication} we have
$$
\|\cov_A^*\cov_Av\|_{L^{\sharp,2}}
\le 2\|d_A^+d_A^*v\|_{L^{\sharp,2}} 
+ c\|v\|_{L^{\sharp,2}}
+ c\|F_A^+\|_{L^{\sharp,2}}\|v\|_{C^0}.
$$
Combining the preceding estimate with that of Lemma
\ref{lem:LinftyL22CovLapEstv}, together with the embedding and
interpolation inequalities $\|v\|_{L^\sharp} \le c\|v\|_{L^4}\le
c\|v\|_{L^2}^{1/2}\|v\|_{C^0}^{1/2}$, and using rearrangement with the last
term yields the desired bound.  In particular, by choosing $\eps(c)$ small
enough that $c\|F_A^+\|_{L^{\sharp,2}}\|v\|_{C^0} \leq 1/2$, we may use
rearrangement to bring the right-hand term $\|v\|_{C^0}$ to the left-hand
side.
\end{proof}

Since $\|d_A^*v\|_{L^2_{1,A}}\le \|v\|_{L^2_{2,A}}$, 
Lemma \ref{lem:LinftyL22Estv} yields an $L^2_{1,A}$ estimate for
$d_A^*v$:

\begin{cor}\label{cor:L21AEstdA*v}
Continue the hypotheses of Lemma \ref{lem:LinftyL22Estv}. Then: 
\begin{align*}
\|d_A^*v\|_{L^2_{1,A}(X)}
&\le
c(1+\|F_A\|_{L^2(X)})
(\|d_A^+d_A^*v\|_{L^{\sharp,2}(X)} + \|v\|_{L^2(X)}).
\end{align*}
\end{cor}

Note that if $a\in\Om^1(\fg_E)$ is $L^2$-orthogonal to $\Ker d_A^+$, so
that $a=d_A^*v$ for some $v\in \Om^+(\fg_E)$, 
and $\Ker d_A^+d_A^{+,*} = 0$, then the estimate of
Corollary \ref{cor:L21AEstdA*v} can be written in the more familiar form
\begin{equation}
\label{eq:UniformFirstOrderEllipticEst}
\begin{aligned}
\|a\|_{L^2_{1,A}(X)}
&\le
c(1+\|F_A\|_{L^2(X)})
(\|d_A^+a\|_{L^{\sharp,2}(X)} + \nu_2[A]^{-1/2}\|a\|_{L^2(X)}),
\end{aligned}
\end{equation}
where we make use of the eigenvalue estimate $\|v\|_{L^2} \le
\nu_2[A]^{-1/2}\|d_A^*v\|_{L^2}$; the term $d_A^+a$ above can be
replaced by $(d_A^++d_A^*)a$ without changing the estimate constants. Here, $\nu_2[A]$ is the least positive eigenvalue of the Laplacian $d_A^+d_A^{+,*}$.

\section{Existence of gauge transformations via the method of continuity}
\label{sec:Continuity}
In this section we complete the proof of Theorem \ref{thm:GaugeFixing}
using the method of continuity. The strategy broadly follows that of
Uhlenbeck's proof of Theorem 2.1 in \cite{UhlLp}. The main new technical
difficulty, not present in \cite{UhlLp}, is the need to compare distances
in the Coulomb-gauge slice $\bS_{A_0}\subset\sA_E^k$ through the
connection $A_0$ and gauge-invariant distances in $\sB_E^k$ from the point
$[A_0]$. It is at this stage of the method of continuity (in proving
openness --- see Lemma \ref{lem:CoulombToDistance}) --- that we need to
employ the special norms and Green's operator estimates described in
Sections \ref{sec:SharpSobolev} and \ref{sec:Green} in order to achieve the
requisite $C^0$ control of gauge transformations; the proof of closedness
works, as one would expect, with standard Sobolev $L^4$ and $L^2_1$ norms.
In \cite{UhlLp}, the $L^2$ norm of the curvature $F_A$ essentially
serves as a gauge-invariant $L^2_1$ measure of distance from $[A]$ to
$[\Ga]$, where $\Ga$ is the product connection on the product $G$ bundle
over the unit ball. 

For $k\ge 2$, we define the following distance functions on $\sB_E^k$,
\begin{align*}
\dist_{L^4}([A],[A_0])
&= \inf_{u\in\sG_E^{k+1}}\|u(A)-A_0\|_4, \\
\dist_{\sL^{\sharp,2}_{1,A_0}}([A],[A_0])
&= \inf_{u\in\sG_E^{k+1}}\left(\|u(A)-A_0\|_{L^{2\sharp,4}}
+ \|d_{A_0}^*(u(A)-A_0)\|_{L^{\sharp,2}}\right), \\
\dist_{L^{\sharp,2}_{1,A_0}}([A],[A_0])
&= \inf_{u\in\sG_E^{k+1}}\left(\|u(A)-A_0\|_{L^2_{1,A_0}}
+ \|d_{A_0}^*(u(A)-A_0)\|_{L^{\sharp,2}}\right), \\
\dist_{L^2_{\ell,A_0}}([A],[A_0])
&= \inf_{u\in\sG_E^{k+1}}\|u(A)-A_0\|_{L^2_{\ell,A_0}}, \qquad 1\le\ell\le k,
\end{align*}
and the following balls with center at $[A_0]$ and radius $\eps$:
\begin{align*}
\barB_{[A_0]}^{1,*,2}(\eps) 
&= \{[A]\in\sB_E^k:\dist_{\sL^{\sharp,2}_{1,A_0}}([A],[A_0]) \le \eps\},\\
\barB_{[A_0]}^{1,\sharp,2}(\eps) 
&= \{[A]\in\sB_E^k:\dist_{L^{\sharp,2}_{1,A_0}}([A],[A_0]) \le \eps\}, \\
\barB_{[A_0]}^{\ell,2}(\eps) 
&= \{[A]\in\sB_E^k:\dist_{L^2_{\ell,A_0}}([A],[A_0]) \le \eps\}.
\end{align*}
Their open counterparts are denoted $B_{[A_0]}^{1,*,2}(\eps)$,
$B_{[A_0]}^{1,\sharp,2}(\eps)$, and $B_{[A_0]}^{\ell,2}(\eps)$,
respectively. 
Our goal in this section is to establish the following:

\begin{thm}\label{thm:GlobalUhlenbeck}
Let $X$ be a closed, smooth four-manifold with metric $g$
and let $G$ be a compact Lie
group. Then there are positive constants $c,z$ with the following
significance.  Let $E$ be a $G$ bundle over $X$ and suppose that $k\ge 2$
is an integer.  Given a point $[A_0]$ in $\sB_E^k$, let $\nu_0[A_0]$ be the
least positive eigenvalue of the Laplacian $\cov_{A_0}^*\cov_{A_0}$
on $\Om^0(\fg_E)$ and set $K_0 =
(1+\nu_0[A_0]^{-1})(1+\|F_{A_0}\|_{L^2})$. Let $\eps_1$ be a 
constant satisfying $0<\eps_1\le zK_0^{-2}(1+\nu_0[A_0]^{-1/2})^{-1}\}$. 
Then the following hold:
\begin{enumerate}
\item For any $[A]\in\sB^k_E$ with
$\dist_{L^{\sharp,2}_{1,A_0}}([A],[A_0]) < \eps_1$,
there is a gauge transformation $u\in \sG_E^{k+1}$ such that
\begin{enumerate}
\item $d_{A_0}^*(u(A)-A_0) = 0$, 
\item $\|u(A)-A_0\|_{L^{2\sharp,4}} \le
cK_0\dist_{\sL^{\sharp,2}_{1,A_0}}([A],[A_0])$. 
\end{enumerate}
\item For any $[A]\in\sB^k_E$ with
$\dist_{L^{\sharp,2}_{1,A_0}}([A],[A_0]) < \eps_1$,
there is a gauge transformation $u\in \sG_E^{k+1}$ such that
\begin{enumerate}
\item $d_{A_0}^*(u(A)-A_0) = 0$, 
\item $\|u(A)-A_0\|_{L^{2\sharp,4}} \le
cK_0\dist_{\sL^{\sharp,2}_{1,A_0}}([A],[A_0])$, 
\item $\|u(A)-A_0\|_{L^2_{1,A_0}} \le
cK_0\dist_{L^{\sharp,2}_{1,A_0}}([A],[A_0])$. 
\end{enumerate}
\end{enumerate}
\end{thm}

Our first proof of Theorem \ref{thm:GlobalUhlenbeck}, via the method of
continuity, occupies Sections \ref{subsec:DistFunction},
\ref{subsec:Closedness} and \ref{subsec:Openness}. A rather different
proof, via a direct application of the inverse function theorem using
$L^{\sharp,2}_2$ gauge transformations, is given in Section \ref{sec:SliceII}.

\subsection{Distance functions on the quotient space}
\label{subsec:DistFunction}
Our first task is to verify the existence of minimizing gauge
transformations $u\in\sG_E^{k+1}$ for the family of distance functions on
$\sB_E^k$ defined above: 
this is established in Lemma \ref{lem:Metric} and the proof uses
the following version of Uhlenbeck's weak compactness theorem.

\begin{prop}\label{prop:StrongUhlenbeck}
Let $X$ be a closed, smooth, Riemannian four-manifold, let $G$ be a compact
Lie group, let $M$ be a positive constant, 
let $A_0$ be an $L^2_2$
connection on a $G$ bundle $E$ over $X$. If $\{A_\al\}$ is a sequence of
$L^2_2$ connections on $E$ such that $\|F_{A_\al}\|_{L^2_{1,A_0}} \le M$,
then there is a subsequence $\{\al'\}\subset \{\al\}$ and a sequence of
$L^2_3$ gauge transformations $\{u_{\al'}\}$ such that
$u_{\al'}(A_{\al'})$ 
converges weakly in $L^2_{2,A_0}$ and strongly in $L^p_{1,A_0}$, for
$1\le p<4$, to an $L^2_2$ connection $A_\8$ on $E$.
\end{prop}

\begin{pf}
{}From the Sobolev embedding $L^2_1\subset L^p$, $2<p<4$, we obtain a
uniform $L^p$ bound $\|F_{A_\al}\|_{L^p}\le cM$ and so,
according to \cite[Theorem 3.6]{UhlLp}, there is a subsequence 
$\{\al'\}\subset \{\al\}$ and a sequence of $L^p_2$ gauge transformations
$\{u_{\al'}\}$ such that $u_{\al'}(A_{\al'})$ converges weakly in $L^p_{1,A_0}$
to an $L^p_1$ connection $A_\8$ on $E$. The stronger conclusion above is
obtained simply by reworking the proof of Theorem 3.6 in \cite{UhlLp},
using the following local estimate for the connections $A_\al$ over small balls
$B\subset X$. Theorem 2.1 of \cite{UhlLp} provides a sequence of
local trivializations $v_\al:P|_B\to B\times G$ such that $a_\al =
v_\al(A_\al)-\Ga$ satisfies $d^*a_\al = 0$ and
$$
\|a_\al\|_{L^p_1(B)} \le c\|F_{A_\al}\|_{L^p(B)}, \qquad 2\le p<4,
$$
where $\Ga$ is the product connection.  Now $F_{A_\al} = da_\al +
a_\al\wedge a_\al$, so
\begin{align*}
\|a_\al\|_{L^2_2(B)} &\le \|da_\al\|_{L^2_1(B)} + \|a_\al\|_{L^2_1(B)} \\
&\le \|a_\al\wedge a_\al\|_{L^2_1(B)} +
\|F_{A_\al}\|_{L^2_1(B)} + \|a_\al\|_{L^2_1(B)}.
\end{align*}
Now, using the multiplication $L^6\times L^3\to L^2$, the embeddings 
$L^3_1\subset L^{12/5}_1\subset L^6$ and $d(a_\al\wedge a_\al) =
da_\al\wedge a_\al - a_\al\wedge da_\al$, we have 
$$
\|d(a_\al\wedge a_\al)\|_{L^2} \le c\|da_{\al}\|_{L^3}\|a_\al\|_{L^6}
 \le c\|a_\al\|_{L^3_1}^2,
$$
while $\|a_\al\wedge a_\al\|_{L^2}\le \|a_\al\|_{L^4(B)}^2 \le
c\|a_\al\|_{L^2_1(B)}^2$. Hence, we obtain
\begin{align*}
\|a_\al\|_{L^2_2(B)} &\le c\|F_{A_\al}\|_{L^2_1(B)} 
+ c\|a_\al\|_{L^3_1(B)}^2 + \|a_\al\|_{L^2_1(B)} \\
&\le c\|F_{A_\al}\|_{L^2_1(B)}(1 + \|F_{A_\al}\|_{L^2_1(B)}) \\
&\le c\|F_{A_\al}\|_{L^2_{1,A_0}(B)}
(1 + \|F_{A_\al}\|_{L^2_{1,A_0}(B)})(1+\|A_0-\Ga\|_{L^2_1}).
\end{align*}
In particular, the sequence of Coulomb-gauge, local connection matrices
$\{a_\al\}$ is bounded in $L^2_2(B)$, so we can extract a weakly
$L^2_2(B)$-convergent and strongly $L^p_1(B)$-convergent subsequence, via
the compactness of embedding $L^2_2(B)\subset L^p_1(B)$ when $1\le p<4$. The
patching argument used to complete the proof of Uhlenbeck's theorem now
proceeds exactly as in \cite{UhlLp} to give the desired conclusion. 
\end{pf}

The proposition is used to extract the desired convergence in the next
lemma.

\begin{lem}\label{lem:Metric}
For any points $[A_0],[A]$ in $\sB_E^k$ there are gauge transformations
such that the following equalities hold:
\begin{align*}
\dist_{L^4}([A],[A_0])
&= \|u(A)-A_0\|_{L^4(X)}, \quad u\in\sG_E^3\tag{1} \\
\dist_{\sL^{\sharp,2}_{1,A_0}}([A],[A_0])
&= \|v(A)-A_0\|_{\sL^{\sharp,2}_{1,A_0}}, \quad v\in\sG_E^3,\tag{2} \\
\dist_{L^{\sharp,2}_{1,A_0}}([A],[A_0])
&= \|w(A)-A_0\|_{L^{\sharp,2}_{1,A_0}}, \quad w\in\sG_E^3,\tag{3} \\
\dist_{L^2_{\ell,A_0}}([A],[A_0])
&= \|w_\ell(A)-A_0\|_{L^2_{\ell,A_0}}, \qquad \ell=1\text{ and }
3\le \ell\le k, \tag{4}
\end{align*}
where $w_1\in \sG_E^3$ and $w_\ell\in \sG_E^{\ell+1}$ in (4).
Moreover, the above distance functions (including the $\ell=2$ distance
function in (4)) are continuous with respect to the quotient $L^2_k$
topology on $\sB_E^k$.
\end{lem}

\begin{pf}
Consider (1).
Let $\{u_\al\}$ be a minimizing sequence in $\sG_E^{k+1}$, so
$\|u_\al(A)-A_0\|_{L^4}$ converges to
$\dist_{L^4}([A],[A_0])$ as $\al\to\8$. Setting $B_\al = u_\al(A) =
A - (d_Au_\al)u_\al^{-1}\in\sA_E^k$, we see that
$B_\al u_\al = Au_\al - d_Au_\al = Au_\al - d_{A_0}u_\al -
[A-A_0,u_\al]$, and thus
\begin{equation}
d_{A_0}u_\al = u_\al(A-A_0) - (B_\al-A_0)u_\al. \label{eq:EllipticGaugeEqnSeq}
\end{equation}
Therefore, as $\|u_\al\|_{C^0}\le c(G)$, we have
$$
\|\cov_{A_0}u_\al\|_{L^4} \le c(\|A-A_0\|_{L^4} + \|B_\al-A_0\|_{L^4}), 
$$
so the sequence $\{u_\al\}\subset
L^2_{k+1}(\gl(E))$ is bounded in $L^4_{1,A_0}(\gl(E)$. So, passing to a
subsequence, we may suppose that $\{u_\al\}$ converges weakly in
$L^4_{1,A_0}(\gl(E))$ and strongly in $L^q(\gl(E))$, via the compact embedding
$L^4_1\subset L^q$, for any $1\le q<\8$, to a limit $u\in L^4_1(\gl(E))$.

We also have $F_{B_\al} = F_{u_\al(A)} = u_\al F_A
u_\al^{-1}$, so $\|F_{B_\al}\|_{L^2} = \|F_A\|_{L^2}$ and as 
$$
\cov_{A_0}F_{B_\al} = (\cov_{A_0}u_\al)\otimes F_A u_\al^{-1} 
+ u_\al (\cov_{A_0}F_A) u_\al^{-1} 
- u_\al F_A \otimes u_\al^{-1}(\cov_{A_0}u_\al)u_\al^{-1},
$$
we see that
\begin{align}
\|\cov_{A_0}F_{B_\al}\|_{L^2} &\le c(\|\cov_{A_0}u_\al\|_{L^4}\|F_A\|_{L^4}
+ \|\cov_{A_0}F_A\|_{L^2}) \label{eq:L21CurvBound}\\
&\le c(1+\|u_\al\|_{L^4_{1,A_0}})\|F_A\|_{L^2_{1,A_0}}. \notag
\end{align}
Hence, the sequence of $L^2_k$ connections $\{B_\al\}$ has curvature
uniformly bounded in $L^2_{1,A_0}$: Proposition \ref{prop:StrongUhlenbeck}
implies, after passing to a subsequence, that the sequence $\{B_\al\}$
converges weakly in $L^2_{2,A_0}$ and strongly in $L^p_{1,A_0}$, for $1\le
p<4$, to an $L^2_2$ $G$ connection $B$ on $E$. {}From
\eqref{eq:EllipticGaugeEqnSeq} we obtain
\begin{equation}
d_{A_0}u = u(A-A_0) - (B-A_0)u, \label{eq:EllipticGaugeEqn}
\end{equation}
a first-order linear elliptic equation in $u$ with $L^2_2$
coefficients. Therefore, $u\in L^2_3(\gl(E))$ and $B = u(A) = A -
(d_Au)u^{-1}$ lies in $\sA_E^2$. It is not {\em a priori\/} clear that the
limit $u$ actually lies in $\sG_E^3$ (since the convergence was only
weakly $L^4_{1,A_0}(\gl(E))$ and strongly $L^q(\gl(E))$): however,
the argument of the last paragraph in the proof of Lemma 4.2.4 in
\cite[p. 130]{DK} applies (using the compactness of the structure group
$G$) and shows that the limit gauge transformation $u$ lies in
$\sG_E^3$. Since $B_\al=u_\al(A)$ converges strongly in $L^p_{1,A_0}$ to
$u(A)$ we now have
$$
\dist_{L^4}([A],[A_0])
= \lim_{\al\to\8}\|u_\al(A)-A_0\|_{L^4}
= \|u(A)-A_0\|_{L^4},
$$
as required in (1). The same argument proves Assertions (2) and (3) and
Assertion (4) when $\ell =1$. The case $\ell\ge 3$ in (4) is straightforward
as we can now apply Lemma
\ref{lem:GaugeSequence} to obtain the desired convergence.

It remains to check $L^2_k$ continuity. We just consider (1), as the
remaining cases are identical.  If $[A_\al]\in\sB_E^k$ is a sequence
converging to $[A_\8]\in\sB_E^k$, then there is a
sequence of gauge transformations $s_\al\in\sG_E^{k+1}$ such that
$s_\al(A_\al)$ converges in $L^2_{k,A_0}$ to $A_\8\in\sA_E^k$ and, in
particular, in $L^4$. But then 
\begin{align*}
&|\dist_{L^4}([A_\al],[A_0]) - \dist_{L^4}([A_\8],[A_0])| \\
&\qquad = |\dist_{L^4}([s_\al(A_\al)],[A_0]) - \dist_{L^4}([A_\8],[A_0])| \\
&\qquad \le \dist_{L^4}([s_\al(A_\al)],[A_\8]) \le \|s_\al(A_\al)-A_\8\|_{L^4},
\end{align*}
and so
$$
\lim_{\al\to\8}\dist_{L^4}([A_\al],[A_0]) 
= \dist_{L^4}([A_\8],[A_0]),
$$
as desired. 
\end{pf}

\subsection{Closedness}
\label{subsec:Closedness}
Let $\fB\subset\barB^{1,\sharp,2}_{[A_0]}(\eps)$ be the subset of points $[A]$
such that there exists a gauge transformation
$u\in\sG_E^{k+1}$ satisfying the conclusions of Assertion (2) of Theorem
\ref{thm:GlobalUhlenbeck}; let $\fB^*\subset\barB^{1,*,2}_{[A_0]}(\eps)$
be the subset of points $[A]$ 
such that there exists a gauge transformation $u\in\sG_E^{k+1}$ satisfying
the conclusions of Assertion (1). As in the proof of Theorem 2.1 in
\cite{UhlLp}, we apply the method of continuity to show that
$\fB^*=\barB^{1,*,2}_{[A_0]}(\eps)$ and
$\fB=\barB^{1,\sharp,2}_{[A_0]}(\eps)$ for small enough $\eps$. Not
surprisingly, we have:

\begin{lem}\label{lem:ConnectedBall}
The balls $\barB_{[A_0]}^{1,*,2}(\eps)$ and
$\barB_{[A_0]}^{1,\sharp,2}(\eps)$ are connected.
\end{lem}

\begin{pf}
If $[A]\in \barB^{1,\sharp,2}_{[A_0]}(\eps)$, there is a 
gauge transformation $u\in\sG_E^{k+1}$
such that $\|u(A)-A_0\|_{L^{\sharp,2}_{1,A_0}}\le \eps$. Then $A_t=A_0 +
t(u(A)-A_0)$, $t\in [0,1]$, 
is a path in $\sA_E^k$ joining $A_0$ to $u(A)$ and
$\|A_t-A_0\|_{L^{\sharp,2}_{1,A_0}} = 
t\|u(A)-A_0\|_{L^{\sharp,2}_{1,A_0}} \le t\eps$, so the path $[A_t]$ lies in
$\barB^{1,\sharp,2}_{[A_0]}(\eps)$ and joins $[A_0]$ to $[A]$. Similarly
for $\barB_{[A_0]}^{1,*,2}(\eps)$.
\end{pf}

Our task then reduces to showing that $\fB^*$ is an open and closed subspace
of $\barB_{[A_0]}^{1,*,2}(\eps)$ and that $\fB$ is an open and closed
subspace of $\barB_{[A_0]}^{1,\sharp,2}(\eps)$. 
First we consider closedness: 

\begin{lem}\label{lem:ClosedSubspace}
The subspaces $\fB^*\subset \barB_{[A_0]}^{1,*,2}(\eps)$ and
$\fB\subset \barB_{[A_0]}^{1,\sharp,2}(\eps)$ are closed.
\end{lem}

\begin{pf}
It suffices to consider the second assertion as the same argument yields
the first.  Suppose $[A_\al]$ is a sequence of points in $\fB$ which
converges in $\sB_E^k$ to a point $[B_\8]$. We may suppose, without loss of
generality, that $A_\al\in\sA_E^k$ is the corresponding sequence of
connections, representing the gauge-equivalence classes $[A_\al]$, which
satisfy the defining conditions for $\fB$:
\begin{align}
d_{A_0}^*(A_\al-A_0) &= 0, \label{eq:SeqUhlCond}\\
\|A_\al-A_0\|_{L^{2\sharp,4}} &\le
cK_0\dist_{\sL^{\sharp,2}_{1,A_0}}([A_\al],[A_0]), \notag\\ 
\|A_\al-A_0\|_{L^2_{1,A_0}} &\le
cK_0\dist_{L^{\sharp,2}_{1,A_0}}([A_\al],[A_0]). \notag
\end{align}
Since $[A_\al]$
converges in $\sB_E^k$ to $[B_\8]$, there is a sequence of gauge
transformations $u_\al\in\sG_E^{k+1}$ such that $B_\al := u_\al(A_\al)$
converges in $L^2_{k,A_0}$ to $B_\8\in \sA_E^k$. Since $B_\al =
u_\al(A_\al)$ and $d_{A_0}^*(A_\al-A_0)=0$, we have
\begin{align}
d_{A_0}u_\al &= u_\al(A_\al-A_0) - (B_\al-A_0)u_\al, 
\label{eq:SeqFirstOrder}\\
d_{A_0}^*d_{A_0}u_\al &= -*(d_{A_0}u_\al\wedge *(A_\al-A_0)) 
- (d_{A_0}^*(B_\al-A_0))u_\al \label{eq:SeqSecondOrder}\\
&\qquad - *(*(B_\al-A_0)\wedge d_{A_0}u_\al), \notag
\end{align}
and so, as $\|u_\al\|_{C^0}\le 1$,
\begin{align*}
\|d_{A_0}u_\al\|_{L^{2\sharp,4}} 
&\le \|A_\al-A_0\|_{L^{2\sharp,4}} + \|B_\al-A_0\|_{L^{2\sharp,4}} \\
\|d_{A_0}^*d_{A_0}u_\al\|_{L^{\sharp,2}}
&\le \|d_{A_0}u_\al\|_{L^{2\sharp,4}}\|A_\al-A_0\|_{L^{2\sharp,4}}
+ \|d_{A_0}^*B_\al-A_0\|_{L^{\sharp,2}} \\
&\qquad \|B_\al-A_0\|_{L^{2\sharp,4}}\|d_{A_0}u_\al\|_{L^{2\sharp,4}}.
\end{align*}
Therefore, the sequence $u_\al$ is bounded in 
$L^2_{2,A_0}(\gl(E))$ and so,
passing to a subsequence, we may suppose that $u_\al$ converges weakly in
$L^2_{2,A_0}(\gl(E))$ (and strongly in $L^p_{1,A_0}$, for any
$p<4$ via the compact embedding $L^2_2\subset L^p_1$) to a limit 
$u_\8\in L^\8\cap L^2_{2,A_0}(\gl(E))$. 

On the other hand, using $A_\al=u_\al^{-1}(B_\al)$, we have
$\|F_{A_\al}\|_{L^2} = \|F_{B_\al}\|_{L^2}$ and
the derivation of \eqref{eq:L21CurvBound} gives
$$
\|\cov_{A_0}F_{A_\al}\|_{L^2} 
\le c(1+\|u_\al\|_{L^4_{1,A_0}})\|F_{B_\al}\|_{L^2_{1,A_0}},
$$
so the sequence $A_\al$ has curvature uniformly bounded in
$L^2_{1,A_0}$. Thus, after passing to a subsequence we may assume by
Proposition \ref{prop:StrongUhlenbeck} that the sequence $A_\al$
converges weakly in $L^2_{2,A_0}$ and strongly in $L^p_{1,A_0}$, $2\le
p<4$, to a limit $A_\8\in\sA_E^2$. 

Taking weak limits in
\eqref{eq:SeqFirstOrder} and \eqref{eq:SeqSecondOrder} yields
\begin{equation}
d_{A_0}u_\8 = u_\8(A_\8-A_0) - (B_\8-A_0)u_\8.
\label{eq:LimitFirstOrder}
\end{equation}
The equation \eqref{eq:LimitFirstOrder} is first order, linear, elliptic
in $u_\8\in L^\8\cap L^2_2$ with $L^2_2$ coefficients. Hence, $u_\8$ is in
$L^2_3(\gl(E))$ and in particular, in $\sG_E^3$, while
$B_\8=u_\8(A_\8)$. {}From \eqref{eq:LimitFirstOrder} we see that
$$
A_\8-A_0 = u_\8^{-1}(B_\8-A_0)u_\8 + u_\8^{-1}d_{A_0}u_\8 
$$
and so, as $d_{A_0}^*(A_\8-A_0=0$, we have
\begin{equation}
d_{A_0}^*(u_\8^{-1}d_{A_0}u_\8+u_\8^{-1}(B_\8-A_0)u_\8) = 0.
\label{eq:UhlCondBootstrap}
\end{equation}  
This is a second-order elliptic equation for $u_\8\in \sG_E^3$ with $L^2_k$
coefficients: in particular, since $u_\8\in L^p_2$ for $2\le p\le 4$, a
standard 
elliptic bootstrapping argument then implies that $u_\8\in L^2_{k+1}$ (see,
for example, the proof of Proposition 3.3 in \cite{FL1}) and therefore
$A_\8 = u_\8^{-1}(B_\8) \in \sA_E^k$.  

Now, taking weak limits in \eqref{eq:SeqUhlCond}, we have
\begin{align*}
d_{A_0}^*(A_\8-A_0) &= \lim_{\al\to\8} d_{A_0}^*(A_\al-A_0) = 0, \\
\|A_\8-A_0\|_{L^{2\sharp,4}} 
&= \lim_{\al\to\8}\|A_\al-A_0\|_{L^{2\sharp,4}} 
\le \lim_{\al\to\8} cK_0\dist_{\sL^{\sharp,2}_{1,A_0}}([A_\al],[A_0]), \\
\|A_\al-A_0\|_{L^{\sharp,2}_{1,A_0}} 
&= \lim_{\al\to\8}\|A_\al-A_0\|_{L^{\sharp,2}_{1,A_0}} 
\le \lim_{\al\to\8} cK_0\dist_{L^{\sharp,2}_{1,A_0}}([A_\al],[A_0]).
\end{align*}
Moreover, as $B_\8=u_\8(A_\8)$ and $u_\8\in\sG_E^{k+1}$,
\begin{align*}
\lim_{\al\to\8} \dist_{\sL^{\sharp,2}_{1,A_0}}([A_\al],[A_0]) 
&= \dist_{\sL^{\sharp,2}_{1,A_0}}([B_\8],[A_0])
= \dist_{\sL^{\sharp,2}_{1,A_0}}([A_\8],[A_0]),\\
\lim_{\al\to\8} \dist_{L^{\sharp,2}_{1,A_0}}([A_\al],[A_0]) 
&= \dist_{L^{\sharp,2}_{1,A_0}}([B_\8],[A_0]),
= cK_0\dist_{L^{\sharp,2}_{1,A_0}}([A_\8],[A_0]),
\end{align*}
where the $L^2_k$ continuity of the distance functions is given by Lemma
\ref{lem:Metric}. 
Therefore, $[B_\8]=[A_\8]\in\fB$. Thus, $\fB$ is closed in $\sB_E^k$ and in
particular, closed in $\barB_{[A_0]}^{1,\sharp,2}(\eps)$, as desired.
\end{pf}

\subsection{Openness}
\label{subsec:Openness}
We must first compare distances from the connection $A_0$ in
the Coulomb slice through $A_0$ in $\sA_E^k$ and gauge-invariant distances
in $\sB_E^k$ from the point $[A_0]$:

\begin{lem}\label{lem:CoulombToDistance}
Let $(X,g)$ be a closed, smooth, Riemannian four-manifold. Then there are
positive constants $c,z$ with the following significance.
Let $K_0=(1+\nu_0[A_0]^{-1})(1+\|F_{A_0}\|_{L^2})$. If  
$A\in\sA_E^k$ satisfies
\begin{itemize}
\item $d_{A_0}^*(A-A_0) = 0$, 
\item $\|A-A_0\|_{L^{2\sharp,4}} \le zK_0^{-1}$,
\end{itemize}
then the following hold:
\begin{enumerate}
\item If $\dist_{\sL^{\sharp,2}_{1,A_0}}([A],[A_0])\le zK_0^{-1}$, then
$$
\|A-A_0\|_{L^{2\sharp,4}} 
\le cK_0\dist_{\sL^{\sharp,2}_{1,A_0}}([A],[A_0]);
$$
\item If $\dist_{L^{\sharp,2}_{1,A_0}}([A],[A_0])\le zK_0^{-1}$, then
\begin{align*}
\|A-A_0\|_{L^{2\sharp,4}} 
&\le cK_0\dist_{\sL^{\sharp,2}_{1,A_0}}([A],[A_0]), \\
\|A-A_0\|_{L^2_{1,A_0}} &\le cK_0\dist_{L^{\sharp,2}_{1,A_0}}([A],[A_0]).
\end{align*}
\end{enumerate}
\end{lem}

\begin{pf}
Recall that for either distance function, minimizing gauge transformations
in $\sG_E^3$ exist by Lemma \ref{lem:Metric}; for convenience, we denote both
by $u\in\sG_E^3$ although they need not {\em a priori\/} coincide. 
Since $B := u(A)=A-(d_Au)u^{-1}\in\sA_E^2$, we have
$$
u(A)-A_0 = u(A-A_0)u^{-1} - (d_{A_0}u)u^{-1}.
$$
Our task, in essence, is to estimate the second term on the right above.
Rewriting this equality gives a first-order, linear elliptic equation in
$u$ with $L^2_2$ coefficients:
\begin{equation}
d_{A_0}u = u(A-A_0) - (B-A_0)u.
\label{eq:FirstEllipticu}
\end{equation}
Let $u_0\in L^2_3(\gl(E))$ be the $L^2$ orthogonal projection of
$u\in\sG_E^3\subset L^2_3(\gl(E))$ onto $\Ker
(d_{A_0}|_{L^2_3})^\perp$, so $u=u_0+\ga$, where $\ga\in\Ker
d_{A_0}\subset\Om^0(\gl(E))$. 
Thus, as $d_{A_0}^*(A-A_0)=0$ and $d_{A_0}u=d_{A_0}u_0$, we see that
\begin{align*}
d_{A_0}^*d_{A_0}u_0 &= -*(d_{A_0}u\wedge *(A-A_0)) + ud_{A_0}^*(A-A_0) \\
&\qquad  - (d_{A_0}^*(B-A_0))u  - *(*(B-A_0)\wedge d_{A_0}u) \\
&= -*(d_{A_0}u_0\wedge *(A-A_0)) 
 - (d_{A_0}^*(B-A_0))u  - *(*(B-A_0)\wedge d_{A_0}u_0).
\end{align*}
Therefore, using the bound $\|u\|_{C^0}\le 1$ for any $u\in\sG_E^3$ (as
the representation for $G$ is orthogonal), we have
\begin{align*}
\|\De_{A_0}^0u_0\|_{L^{\sharp,2}} 
&\le \|d_{A_0}u_0\|_{L^{2\sharp,4}}\|A-A_0\|_{L^{2\sharp,4}} 
+ \|d_{A_0}^*(B-A_0)\|_{L^{\sharp,2}}\|u\|_{C^0}  \\
&\qquad + \|B-A_0\|_{L^{2\sharp,4}}\|d_{A_0}u_0\|_{L^{2\sharp,4}} \\
&\le C\left(\|A-A_0\|_{L^{2\sharp,4}} +
\|B-A_0\|_{L^{2\sharp,4}}\right)\|d_{A_0}^*d_{A_0}u_0\|_{L^{\sharp,2}} \\
&\qquad + \|d_{A_0}^*(B-A_0)\|_{L^{\sharp,2}},
\end{align*}
where $C=cK_0$. Now $\|B-A_0\|_{L^{2\sharp,4}} \le
c\|B-A_0\|_{L^2_{1,A_0}}$ via the embedding $L^2_1\subset L^{2\sharp,4}$
of Lemma \ref{lem:EmbeddingSobolevInto*}. For either 
$\dist_{\sL^{\sharp,2}_{1,A_0}}([A],[A_0])\le \quarter C^{-1}$ or
$\dist_{L^{\sharp,2}_{1,A_0}}([A],[A_0])\le \quarter C^{-1}$
and $\|A-A_0\|_{L^{2\sharp,4}}\le \quarter C^{-1}$, rearrangement yields
\begin{equation}
\|\De_{A_0}^0u_0\|_{L^{\sharp,2}} 
\le 2\|d_{A_0}^*(B-A_0)\|_{L^{\sharp,2}}.  \label{eq:Lapu0Est}
\end{equation}
On the other hand, from Lemma \ref{lem:L2*2GreenEst} we have
\begin{align}
\|u_0\|_{L^{\sharp,2}_{2,A_0}} 
&\le C\|\De_{A_0}^0u_0\|_{L^{\sharp,2}}, \label{eq:u0L22andL4Est}\\
\|u_0\|_{L^4_{1,A_0}} 
&\le C\|\De_{A_0}^0u_0\|_{L^{\sharp,2}}, \notag
\end{align}
where $C=cK_0$ and the second bound
follows from the embedding $L^2_2\subset L^4_1$.
So, combining \eqref{eq:Lapu0Est} and \eqref{eq:u0L22andL4Est} yields:
\begin{align}
\|u_0\|_{L^{\sharp,2}_{2,A_0}} 
&\le C\|d_{A_0}^*(B-A_0)\|_{L^{\sharp,2}}, \label{eq:u0L22andL4EstDist}\\
\|u_0\|_{L^4_{1,A_0}} 
&\le C\|d_{A_0}^*(B-A_0)\|_{L^{\sharp,2}}. \notag
\end{align}
Consequently, using $d_{A_0}u=d_{A_0}u_0$ and \eqref{eq:FirstEllipticu}
rewritten in the form,
\begin{equation}
u^{-1}(B-A_0)u = (A-A_0) - u^{-1}d_{A_0}u_0, \label{eq:BACompare}
\end{equation}
we obtain
\begin{align}
\|A-A_0\|_{L^{2\sharp,4}} 
&\le \|u^{-1}(B-A_0)u\|_{L^{2\sharp,4}} 
+ \|u^{-1}d_{A_0}u_0\|_{L^{2\sharp,4}},
\label{eq:CoulombL2sharp4Est} \\
\|A-A_0\|_{L^2_{1,A_0}} 
&\le \|u^{-1}(B-A_0)u\|_{L^2_{1,A_0}} + \|u^{-1}d_{A_0}u_0\|_{L^2_{1,A_0}}.
\label{eq:CoulombL21Est}
\end{align}
{}From \eqref{eq:CoulombL2sharp4Est} and
\eqref{eq:u0L22andL4EstDist}, we see that 
\begin{align*}
\|A-A_0\|_{L^{2\sharp,4}} 
&\le \|B-A_0\|_{L^{2\sharp,4}} + \|d_{A_0}u_0\|_{L^{2\sharp,4}} \\
&\le \dist_{\sL^{\sharp,2}_{1,A_0}}([A],[A_0]) 
+  C\|d_{A_0}^*(B-A_0)\|_{L^{\sharp,2}}\\
&\le (1+C)\dist_{\sL^{\sharp,2}_{1,A_0}}([A],[A_0]), 
\end{align*}
giving the desired $L^{2\sharp,4}$ estimate for $A-A_0$. 

Considering the first term in \eqref{eq:BACompare}, we have 
\begin{align*}
\cov_{A_0}(u^{-1}(B-A_0)u)
&= -u^{-1}(\cov_{A_0}u)u^{-1}\otimes(B-A_0)u 
+ u^{-1}(\cov_{A_0}(B-A_0))u \\
&\qquad + u^{-1}(B-A_0)\otimes \cov_{A_0}u,
\end{align*}
and so applying \eqref{eq:u0L22andL4EstDist},
noting that $\cov_{A_0}u=\cov_{A_0}u_0$ and $\|u\|_{C^0}\le 1$, we have
\begin{align*}
\|\cov_{A_0}(u^{-1}(B-A_0)u)\|_{L^2}
&\le \|\cov_{A_0}u_0\|_{L^4}\|B-A_0\|_{L^4}
+ \|\cov_{A_0}(B-A_0)\|_{L^2} \\
&\le C\dist_{L^{\sharp,2}_{1,A_0}}^2([A],[A_0])
+ \dist_{L^{\sharp,2}_{1,A_0}}([A],[A_0]).
\end{align*}
Thus, if $\dist_{L^{\sharp,2}_{1,A_0}}([A],[A_0])\le\quarter C^{-1}$,
say, we obtain
\begin{equation}
\|\cov_{A_0}(u^{-1}(B-A_0)u)\|_{L^2}
\le 2\dist_{L^{\sharp,2}_{1,A_0}}([A],[A_0]).
\label{eq:BACompareL21Est1}
\end{equation}
Similarly, considering the second term in \eqref{eq:BACompare}, we have
$$
\cov_{A_0}(u^{-1}d_{A_0}u_0)
= -u^{-1}(\cov_{A_0}u)u^{-1}\otimes d_{A_0}u
+ u^{-1}\cov_{A_0}d_{A_0}u
$$
and therefore, by \eqref{eq:u0L22andL4EstDist}, we see that
\begin{align*}
\|\cov_{A_0}(u^{-1}d_{A_0}u_0)\|_{L^2}
&\le \|\cov_{A_0}u_0\|_{L^4}^2 + \|\cov_{A_0}^2u_0\|_{L^2} \\
&\le C\dist_{L^{\sharp,2}_{1,A_0}}([A],[A_0])
\left(1+C\dist_{L^{\sharp,2}_{1,A_0}}([A],[A_0])\right).
\end{align*}
Provided $\dist_{L^{\sharp,2}_{1,A_0}}([A],[A_0])\le\quarter C^{-1}$, we
obtain
\begin{equation}
\|\cov_{A_0}(u^{-1}d_{A_0}u_0)\|_{L^2}
\le 2C\dist_{L^{\sharp,2}_{1,A_0}}([A],[A_0]).
\label{eq:BACompareL21Est2}
\end{equation}
Taking the $L^2$ norm of \eqref{eq:BACompare} and applying
\eqref{eq:u0L22andL4EstDist} to estimate the second term gives
\begin{align*}
\|A-A_0\|_{L^2} 
&\le \|B-A_0\|_{L^2} + \|d_{A_0}u_0\|_{L^2} \\
&\le \dist_{L^{\sharp,2}_{1,A_0}}([A],[A_0]) 
+ C\|d_{A_0}^*(B-A_0)\|_{L^{\sharp,2}}, 
\end{align*}
and so
\begin{equation}
\|A-A_0\|_{L^2}\le (1+C)\dist_{L^{\sharp,2}_{1,A_0}}([A],[A_0]). 
\label{eq:BACompareL21Est3}
\end{equation}
Combining the estimates \eqref{eq:CoulombL21Est},
\eqref{eq:BACompareL21Est1}, \eqref{eq:BACompareL21Est2}, and 
\eqref{eq:BACompareL21Est3} yields 
$$
\|A-A_0\|_{L^2_{1,A_0}}
\le 2(1+C)\dist_{L^{\sharp,2}_{1,A_0}}([A],[A_0]),
$$
giving the desired $L^2_{1,A_0}$ estimate for $A-A_0$. 
\end{pf}

Naturally, a comparison --- going in the other direction --- 
of distances from $A_0$ in the Coulomb slice in
$\sA_E^k$ through $A_0$ and gauge-invariant distances
in $\sB_E^k$ from the point $[A_0]$ is elementary: If $A\in\bS_{A_0}$ and 
$\|A-A_0\|_{L^2_{k,A_0}}<\de$, say, then Lemma
\ref{lem:EmbeddingSobolevInto*} implies that
\begin{align}
\dist_{\sL^{\sharp,2}_{1,A_0}}([A],[A_0])
&\le c\|A-A_0\|_{L^2_{k,A_0}} < c\de, \qquad k\ge 1, 
\label{eq:DistanceToCoulomb}\\
\dist_{L^{\sharp,2}_{1,A_0}}([A],[A_0])
&\le c\|A-A_0\|_{L^2_{k,A_0}} < c\de, \qquad k\ge 2, \notag
\end{align}
for some positive constant $c(X,g,k)$. The observation is used in concluding
that $\fB^*$, $\fB$ are open subspaces of $\barB^{1,*,2}_{[A_0]}(\eps_1)$,
$\barB^{1,\sharp,2}_{[A_0]}(\eps_1)$, respectively:

\begin{lem}\label{lem:OpenSubspace}
Let $(X,g)$ be a closed, smooth, Riemannian four-manifold and let $G$ be a
compact Lie group.  Then there is a positive constant $z$ with the
following significance. Let $K_0 =
(1+\nu_0[A_0]^{-1})(1+\|F_{A_0}\|_{L^2})$.  If $\eps_1<
zK_0^{-2}(1+\nu_0[A_0]^{-1/2})^{-1}\}$, then the following hold:
\begin{itemize}
\item $\fB^*$ is an open subspace of $\barB^{1,*,2}_{[A_0]}(\eps_1)$;
\item $\fB$ is an open subspace of $\barB^{1,\sharp,2}_{[A_0]}(\eps_1)$.
\end{itemize}
\end{lem}

\begin{pf}
It suffices to consider the second assertion, as the argument for the
first is identical.
Suppose $[A]\in \fB$ and that $A\in\sA_E^k$ is a representative satisfying
the defining conditions for $\fB$. Then $A$ satisfies $d_{A_0}^*(A-A_0)=0$
and the estimates
\begin{align*}
\|A-A_0\|_{L^{2\sharp,4}} &\le c_0K_0\dist_{\sL^{\sharp,2}_{1,A_0}}([A],[A_0]),
\le c_0K_0\eps_1\\
\|A-A_0\|_{L^2_{1,A_0}} &\le c_0K_0\dist_{L^{\sharp,2}_{1,A_0}}([A],[A_0])
\le c_0K_0\eps_1,
\end{align*}
while $\|A-A_0\|_{L^{2\sharp,4}}\le c_1\|A-A_0\|_{L^2_{1,A_0}}$ via the Sobolev
embedding $L^2_1\subset L^{2\sharp,4}$ and Kato's inequality.
Consequently, if $c_1c_0K_0\eps_1\le\half\eps_0$, then
$A\in\bB_{A_0}^4(\eps_0)\subset\sA_E^k$ and we see that
$$
\barB^{1,\sharp,2}_{[A_0]}(\eps_1)\subset \pi(\bB_{[A_0]}^4(\eps_0)).
$$
Lemma \ref{lem:LocalDiffeo} implies that the map
$\pi:\bB_{[A_0]}^4(\eps_0)/\Stab_{A_0}\to\sB_E^k$ 
given by $A'\mapsto [A']$ is a local homeomorphism onto its
image $\pi(\bB_{[A_0]}^4(\eps_0))$ for any
$\eps_0<z(1+\nu_0[A_0]^{-1/2})^{-1}$. In particular, if $A'\in
\bB_{[A_0]}^4(\eps_0)$ and $\|A'-A\|_{L^2_{k,A_0}}<\de$, then 
$A'\in B^{1,\sharp,2}_{[A_0]}(\eps_1)\subset
\barB^{1,\sharp,2}_{[A_0]}(\eps_1)$ for small enough $\de$. 

The embedding $L^2_1\subset L^{2\sharp,4}$ and 
Lemma \ref{lem:CoulombToDistance} imply that if $\|A'-A_0\|_{L^2_{1,A_0}} \le
zK_0^{-1}$ 
and $\dist_{L^{\sharp,2}_{1,A_0}}([A'],[A_0]) \le zK_0^{-1}$, then
\begin{align*}
\|A'-A_0\|_{L^{2\sharp,4}} 
&\le cK_0\dist_{\sL^{\sharp,2}_{1,A_0}}([A'],[A_0]) \le cK_0\eps_1, \\
\|A'-A_0\|_{L^2_{1,A_0}} &\le cK_0\dist_{L^{\sharp,2}_{1,A_0}}([A'],[A_0])
\le cK_0\eps_1.
\end{align*}
These inequalities are satisfied by $A$; moreover
$\dist_{L^{\sharp,2}_{1,A_0}}([A],[A_0])\le\eps_1$ and  
$\|A-A_0\|_{L^2_{1,A_0}}\le c_0K\eps_1$. Require that $\eps_1\le \half
zK_0^{-1}$ and $c_0K_0\eps_1\le \half zK_0^{-1}$, so 
$\eps_1\le \half z\min\{1,c_0\}K_0^{-2}$.
Hence, if $A'$ is $L^2_{k,A_0}$-close enough to $A$ (where $k\ge 2$),
we can ensure $[A']$
obeys the last two defining conditions for $\fB$ and so $[A']\in\fB$. Thus,
$\fB\subset \barB^{1,\sharp,2}_{[A_0]}(\eps_1)$ is open, as desired.
\end{pf}

We can now conclude the proof of Theorem \ref{thm:GlobalUhlenbeck}:

\begin{proof}[Proof of Theorem~\ref{thm:GlobalUhlenbeck}]
Lemmas \ref{lem:ClosedSubspace} 
and \ref{lem:OpenSubspace} imply that $\fB$ is an open and closed subset of
the connected space $\barB^{1,\sharp,2}_{[A_0]}(\eps_1)$, so 
$\fB = \barB^{1,\sharp,2}_{[A_0]}(\eps_1)$. 
Similarly for $\fB^*$ and 
$\barB^{1,*,2}_{[A_0]}(\eps_1)$ and hence the result follows.
\end{proof}

Similarly, we conclude the proof of Theorem \ref{thm:GaugeFixing}:

\begin{proof}[Proof of Theorem~\ref{thm:GaugeFixing}]
Given Theorem \ref{thm:GlobalUhlenbeck},
the only assertion left unaccounted for is the uniqueness of the gauge
transformation $u\in\sG_E^{k+1}$, modulo $\Stab_{A_0}$. But this follows
from Lemma \ref{lem:injective}, just as in the paragraph immediately
following the proof of that lemma.
\end{proof}

\section{Critical-exponent Sobolev norms and the group of gauge
transformations} 
\label{sec:SharpGaugeGroup}
We now define an $L^{\sharp,2}_2$ space of gauge transformations,
by analogy with the definition of $\sG_E^{k+1}$ when $k\ge 2$, and set
$$
\sG_E^{2,\sharp,2} := \{u\in L^{\sharp,2}_2(\gl(E)): u\in G\text{ a.e.}\} 
\subset L^2_k(\gl(E)).
$$
It is not entirely clear {\em a priori\/}
that $\sG_E^{2,\sharp,2}$ is a Banach Lie group.
In the case of its counterpart, $\sG_E^{k+1}$, the manifold structure
follows from the fact that the exponential map 
$$
\Exp:T_{\id_E}\sG_E=\Om^0(\fg_E)\to \sG_E, \quad \zeta\mapsto \Exp\zeta,
$$
extends to a smooth map $\Exp:L^2_{k+1}(\fg_E)\to L^2_{k+1}(\fg_E)$ and
defines a system of smooth coordinate charts for $\sG_E^{k+1}$. Here, $\Exp$
is defined pointwise at $u\in\sG_E$ for $\zeta\in T_{\id_E}\sG_E$ by setting 
$$
(\Exp_u\zeta)(x) := \exp_{u(x)}(\zeta(x)), \qquad x\in X,
$$
where $\exp:\fg\to G$ is the usual, $C^\8$ exponential map for the Lie
group $G$ on the right-hand side \cite[Appendix A]{FU}. 

To verify that $\sG_E^{2,\sharp,2}$ is in fact a Banach Lie group
we will need estimates for the covariant derivatives of the exponential
map. The estimates below follow by reworking 
the usual proof of the Sobolev lemma for
left composition of Sobolev sections by smooth vector bundle maps
\cite[Lemma 9.9]{Palais}. The difference
here is that we keep track of the dependence of the constants on the
geometric data: this precision is required for the implicit function
argument in the next section in order to complete the proof of our slice
theorem.

For $\chi,\zeta,\xi\in\Om^0(\fg_E)$, the differentials
\begin{align*}
&(D\Exp)_\chi: \Om^0(\fg_E)\to T_{\Exp\chi}\sG_E, 
\qquad \zeta\mapsto (D\Exp)_\chi\zeta, \\
&(D^2\Exp)_{\chi,\zeta}:\Om^0(\fg_E) \to T_{\Exp\chi}\sG_E, 
\qquad \xi\mapsto (D^2\Exp)_{\chi,\zeta}\xi, 
\end{align*}
are defined pointwise by setting 
\begin{align*}
(D\Exp)_\chi\zeta|_x &= (D\exp)_{\chi(x)}\zeta(x), \\ 
(D^2\Exp)_{\chi,\zeta}\xi|_x &= (D\exp)_{\chi(x),\zeta(x)}\xi(x), 
\end{align*}
for any $x\in X$. When writing the differential
$(D^2\Exp)_{\chi,\zeta}$ above, we have identified 
$T_{(D\exp)_\chi\zeta}(T_{\Exp\chi}\sG_E)$ with $T_{\Exp\chi}\sG_E$.

The maps $(D\Exp)_\chi:\Om^0(\fg_E)\to \Om^0(\fg_E)$ and
$(D^2\Exp)_{\chi,\zeta}:\Om^0(\fg_E)\to \Om^0(\fg_E)$ extend linearly to
maps
\begin{align*}
&(D\Exp)_\chi:C^\8(\otimes^\ell(T^*X)\otimes \fg_E)\to
C^\8(\otimes^\ell(T^*X)\otimes\fg_E), \\
&(D^2\Exp)_\chi:C^\8(\otimes^\ell(T^*X)\otimes \fg_E)\to
C^\8(\otimes^\ell(T^*X)\otimes\fg_E),
\end{align*}
for $\ell\ge 1$, by setting
\begin{align*}
(D\Exp)_\chi(\theta\otimes\zeta) 
&= \theta\otimes(D\Exp)_\chi\zeta, \\
(D^2\Exp)_{\chi,\zeta}(\theta\otimes\xi) 
&= \theta\otimes(D^2\Exp)_{\chi,\zeta}\xi, 
\end{align*}
for $\theta\in\otimes^\ell(T^*X)$ and $\xi\in \Om^0(\fg_E)$.
As usual, we embed $\sG_E\subset\Om^0(\fg_E)$ in order to compute the
covariant derivatives of sections $u\in\sG_E$.

\begin{lem}\label{lem:PtExpEstimates}
Let $G$ be a compact Lie group. Then there is a
positive constant $c(G)$ with the following significance. 
Let $X$ be a closed, smooth, Riemannian four-manifold.  If 
$A$ is a $C^\8$ connection on a $G$ bundle $E$, and $\chi\in
\Om^0(\fg_E)$, then $e^\chi\in\sG_E$ satisfies the following
pointwise bounds:
\begin{align}
|\cov_A e^\chi| &\le |\cov_A\chi| + c|\chi||\cov_A\chi|, \tag{1} \\
|\cov_A^2 e^\chi| &\le c(|\chi| + |\cov_A\chi|)|\cov_A\chi| 
+ c(1+|\chi|)|\cov_A^2\chi|, \tag{2} \\
|\cov_A^*\cov_A e^\chi| &\le c(|\chi| + |\cov_A\chi|)|\cov_A\chi| 
+ c(1+|\chi|)|\cov_A^*\cov_A\chi|. \tag{3}
\end{align}
\end{lem}

\begin{pf}
We have
$$
\cov_A e^\chi = \cov_A(\Exp\chi) =
(D\Exp)_\chi\circ d_A\chi \in \Om^1(\fg_E). 
$$
Since $(D\exp)_0=\id_\EE$ and $\exp:\fg\to G$ is analytic,
we have the pointwise bound
$|(D\exp)_{\chi(x)}-\id_\EE|\le c(G)|\chi(x)|$ and thus a pointwise bound
$$
|(D\Exp)_\chi-\id_E|\le c|\chi|,
$$
noting that $(D\Exp)_0 = \id_E$. Therefore, we have
$$
|\cov_A e^\chi| \le |\cov_A\chi| + c|\chi||\cov_A\chi|,
$$
which gives the first assertion.

Define $\Phi(\chi,\zeta) = (D\Exp)_\chi(\zeta)\in \Om^1(\fg_E)$, for
$\chi\in\Om^0(\fg_E)$ and $\zeta\in\Om^1(\fg_E)$, noting that $\Phi$ is
nonlinear in $\chi$, but linear in $\zeta$. Thus,
$$
\cov_A^2u = (D_1\Phi)_{(\chi,\cov_A\chi)}(\cov_A\chi) 
+ (D_2\Phi)_{\chi}(\cov_A^2\chi),
$$
where $D_i\Phi$, $i=1,2$, denote the partial derivatives of $\Phi$ with
respect to first and second variables. Since $(D\Phi)_{(0,0)} =
(D^2\Exp)_{(0,0)} = \id_E$, as $(D^2\exp)_{0,0} = \id_\EE$, 
and $\exp:\fg\to G$ is analytic we have the pointwise bound
$$
|\cov_A^2u| \le c(|\chi| + |\cov_A\chi|)|\cov_A\chi| 
+ c(1+|\chi|)|\cov_A^2\chi|,
$$
giving the second assertion.
Similarly, as $*\Phi(\chi,\zeta) = \Phi(\chi,*\zeta)$ and $\cov_A^*\cov_Au
= -*\cov_A*\cov_Au$, we have
$$
|\cov_A^*\cov_Au| \le c(|\chi| + |\cov_A\chi|)|\cov_A\chi| 
+ c(1+|\chi|)|\cov_A^*\cov_A\chi|, 
$$
giving the third assertion.
\end{pf}

The preceding pointwise bounds for $\cov_Au$, $\cov_A^2u$, and
$\cov_A^*\cov_Au$ yield the following estimates for the exponential map:

\begin{lem}\label{lem:SobolevExpEstimates}
Let $G$ be a compact Lie group. Then there is a positive constant $c(G)$ with
the following significance.  Let $X$ be a closed, smooth, Riemannian
four-manifold. If $k\ge 2$ is an integer (so $L^2_{k+1}\subset C^0$), $A$
is an $L^2_k$ connection on a $G$ bundle $E$, and $\chi\in
L^2_{k+1}(\fg_E)$, then $e^\chi\in\sG_E^{k+1}$ satisfies
\begin{align}
\|\cov_Ae^\chi\|_{L^2(X)} 
&\le \|\cov_A\chi\|_{L^2(X)} + c\|\chi\|_{C^0(X)}\|\cov_A\chi\|_{L^2(X)},
\tag{1} \\
\|\cov_Ae^\chi\|_{L^\sharp(X)} 
&\le \|\cov_A\chi\|_{L^\sharp(X)} 
+ c\|\chi\|_{C^0(X)}\|\cov_A\chi\|_{L^\sharp(X)},
\tag{2} \\
\|\cov_A^2e^\chi\|_{L^2(X)} &\le c\|\chi\|_{C^0(X)}\|\cov_A\chi\|_{L^2(X)} 
+ \|\cov_A\chi\|_{L^4(X)}^2 \tag{3}\\
&\qquad + c(|1+\|\chi\|_{C^0(X)})\|\cov_A^2\chi\|_{L^2(X)}, \notag\\  
\|\cov_A^*\cov_Ae^\chi\|_{L^\sharp(X)} 
&\le c\|\chi\|_{C^0(X)}\|\cov_A\chi\|_{L^\sharp(X)} 
+ \|\cov_A\chi\|_{L^{2\sharp}(X)}^2 \tag{4}\\
&\qquad + c(|1+\|\chi\|_{C^0(X)})\|\cov_A^*\cov_A\chi\|_{L^\sharp(X)}, 
\notag
\end{align}
Moreover, the preceding bounds continue to hold for $\chi\in
L^{\sharp,2}_2(\fg_E)\subset C^0(\fg_E)$, with $A$ an $L^{\sharp,2}_1$
connection on $E$, and $\Exp:\Om^0(\fg_E)\to\sG_E$ extends to a continuous
map $\Exp:L^{\sharp,2}_2(\fg_E)\to\sG_E^{2,\sharp,2}$.
\end{lem}

Let $\sA_E^{1,\sharp,2}={A_0}+L^{\sharp,2}_{1,A_0}(\La^1\otimes\fg_E)$, for any
$C^\8$ reference connection $A_0$ on $E$. Recall that we have an embedding
$L^{\sharp,2}_2(\fg_E)\subset C^0(\fg_E)$ and that the space
$L^{\sharp,2}_2(\fg_E)$ is an algebra, while
$L^{\sharp,2}_1(\La^1\otimes\fg_E)$ 
and $L^2_1(\La^1\otimes\fg_E)$ are
$L^{\sharp,2}_2(\fg_E)$-modules. Therefore, the 
proofs of Propositions (A.2) and (A.3) in \cite{FU} extend easily to give
the following analogue for $\sG_E^{2,\sharp,2}$ in place of $\sG_E^{k+1}$:

\begin{lem}\label{lem:SharpGaugeGroup}
Let $X$ be a closed Riemannian four-manifold and let $E$
be a Hermitian vector bundle over $X$. Then the following hold.
\begin{enumerate}
\item The space $\sG_E^{2,\sharp,2}$ is a Banach Lie group with
Lie algebra $T_{\id_E}\sG_E^{2,\sharp,2}=L^{\sharp,2}_2(\fg_E)$;
\item The action of $\sG_E^{2,\sharp,2}$ on $\sA_E^1$ and 
$\sA_E^{1,\sharp,2}$ is smooth;
\item For $A\in \sA_E^{1,\sharp,2}$,
the differential, at the identity $\id_E\in \sG_E^{2,\sharp,2}$, 
of the map $\sG_E^{2,\sharp,2}\to \sA_E^{1,\sharp,2}$
given by $u\mapsto u(A)=A-(d_Au)u^{-1}$ is 
$\zeta\mapsto -d_A\zeta$ as a map
$L^{2,\sharp}_2(\fg_E)\to 
L^{\sharp,2}_1(\La^1\otimes\fg_E)$, and similarly for $A\in\sA_E^1$.
\end{enumerate}
\end{lem}

\section{Existence of gauge transformations via the inverse function theorem} 
\label{sec:SliceII}
Our goal in this section is to give an alternative, `direct' proof of
Theorem \ref{thm:GlobalUhlenbeck} via the inverse function theorem. A direct
argument --- due to our overarching constraint of ensuring that the
constants given there 
ultimately depend at most on the $L^2$ norm of the curvature and the least
positive eigenvalue $\nu_0[A_0]$ --- appears to be difficult within the
standard framework of $L^p_2$ ($p>2$) gauge transformations acting on $L^p_1$
connections; of course, if this constraint is dropped then a direct proof
is standard. On the other hand, with the results of the last section, it is
fairly straightforward within the framework of $L^{\sharp,2}_2$ gauge
transformations.

We already know that $\pi(\bB^4_{A_0}(\eps_0))$ is open in $\sB_E^k$, so it
necessarily contains an $L^2_{k,A_0}$-ball centered at $[A_0]$.  Via the
inverse function theorem we estimate the radii of $\sL^{\sharp,2}_{1,A_0}$
and $L^{\sharp,2}_{1,A_0}$ balls, $B_{[A_0]}^{1,\sharp,2}(\eps)$ and
$B_{[A_0]}^{1,\sharp,2}(\eps)$, which are contained in
$\pi(\bB^4_{A_0}(\eps_0))$. Let us first dispose of the question of
regularity for solutions to the second-order gauge-fixing equation:

\begin{lem}\label{lem:GaugeFixRegularity}
Let $X$ be a closed, Riemannian four-manifold. Then there is a constant
$\eps$ with the following significance. Let $G$ be compact Lie
group and let $k\ge 2$ be an integer. Suppose that $A_0$ is an $L^2_k$
connection on a $G$ bundle $E$, that $a\in L^2_k(\La^1\otimes\fg_E)$ and
$\chi\in 
L^\sharp_2(\fg_E)$, and that $u=e^\chi$ is a solution to
$$
d_{A_0}^*\left((d_{A_0}u)u^{-1}-uau^{-1}\right)=0.
$$
If $\|d_{A_0}u\|_{L^4}<\eps$ then
$\chi\in L^2_{k+1}(\fg_E)$ and $u=e^\chi\in\sG_E^{k+1}$.
\end{lem} 

\begin{pf}
Differentiation and right multiplication by $u$ yields
\begin{align}
&d_{A_0}^*d_{A_0}u + *((*d_{A_0}u)\wedge u^{-1}d_{A_0}u)
+ *(d_{A_0}u\wedge *a) + ud_{A_0}^*a \label{eq:LapGaugeFix} \\
&\qquad + *(ua\wedge *u^{-1}d_{A_0}u) = 0. \notag
\end{align}
From Lemma \ref{lem:Embedding*IntoSobolev} we know that $u\in C^0\cap
L^2_1$ and so the last four terms in \eqref{eq:LapGaugeFix} are in $L^2$.
Hence, $d_{A_0}^*d_{A_0}u$ is in $L^2$ and so $u\in L^2_2$ by elliptic
regularity for $d_{A_0}^*d_{A_0}$. The Sobolev embedding $L^2_1\subset L^4$
and multiplication $L^4\times L^q\to L^p$ for $2\le p<4$ and $1/p=1/4+1/q$
(so $4\le q<\8$) now show that the last three terms in
\eqref{eq:LapGaugeFix} are in $L^p$, so the equation takes the
simpler form
\begin{equation}
d_{A_0}^*d_{A_0}u + *((*d_{A_0}u)\wedge u^{-1}d_{A_0}u) = v,
\label{eq:SimpleLapGaugeFix}
\end{equation}
where $v\in L^p(\fg_E)$ is the tautologically defined right-hand side and
$u\in L^\8\cap L^2_2$. Setting $b=d_{A_0}u$ and noting that $d_{A_0}b =
F_{A_0}u$, with $F_{A_0}\in L^2_{k-1}(\La^2\otimes\fg_E)\subset
L^2_1(\La^2\otimes\fg_E)$ and 
$F_{A_0}u\in L^2_1(\La^2\otimes\fg_E)$. Thus, we may conveniently rewrite
\eqref{eq:SimpleLapGaugeFix} as a first-order elliptic equation in $b\in
L^2_1(\La^1\otimes\fg_E)$, 
\begin{equation}
(d_{A_0}^*+d_{A_0})b + *((*b)\wedge u^{-1}b) 
= v' \in L^p(\fg_E)\oplus L^p(\La^2\otimes\fg_E),
\label{eq:SimpleFirstOrdGaugeFix}
\end{equation}
where $2<p<4$ and $v'=F_{A_0}u + v$. Finally, 
\eqref{eq:SimpleFirstOrdGaugeFix} can be rewritten as a local equation by
writing $A_0=\Ga + a_0$, where 
$\Ga$ is the product connection in a local trivialization for $E$ over a
small ball $U\subset X$. Thus,
the operator $d_{A_0}^*+d_{A_0}$ is replaced by $d^*+d$ in
\eqref{eq:SimpleLapGaugeFix} and the additional terms are absorbed into the
$L^p$ inhomogeneous term $v'$ to give:
\begin{equation}
(d^*+d)b + *((*b)\wedge u^{-1}b) = v'' \in 
L^p(U,\fg_E)\oplus L^p(U,\La^2\otimes\fg_E).
\label{eq:VSimpleFirstOrdGaugeFix}
\end{equation}
This is a first-order, elliptic
equation with a quadratic non-linearity and Proposition 3.10 in \cite{FL1}
implies that the solution $b=d_{A_0}u\in L^2_1(U,\La^1\otimes\fg_E)$ is
necessarily in $L^p_1(U',\La^1\otimes\fg_E)$ for $U'\Subset U$, provided
$\|b\|_{L^4(U)}<\eps(g,p,U)$, and so $u\in L^p_2(U',\fg_E)$. In particular,
we find that $b\in L^p_1(X,\La^1\otimes\fg_E)$ and
$u\in L^p_2(X,\fg_E)$ for any $2<p<4$, provided
$\|d_{A_0}u\|_{L^4}<\eps(g,p,X)$. The bootstrapping argument of Proposition
3.3 in \cite{FL1} now implies that $d_{A_0}u\in
L^2_k(X,\La^1\otimes\fg_E)$. Thus 
$u\in \sG_E^{k+1}$ and $\chi\in L^2_{k+1}(X,\fg_E)$, as desired.
\end{pf}

We can now proceed to the main argument:

\begin{thm}\label{thm:CoulombConn}
Let $X$ be a closed, Riemannian four-manifold and let $G$ be compact Lie
group. Then there are positive constants $c,z$ with the following
significance.  Let $E$ be a $G$ bundle over $X$ and suppose that 
that $A_0\in\sA_E^2$, let $K_0[A_0]=(1+\nu_0[A_0]^{-1})(1+\|F_{A_0}\|_{L^2}$
and let $\eps_1$ be a constant satisfying
$$
0 < \eps_1 \le zK_0^{-2}.
$$
Then for any $A\in\sA^2_E$ such that
$\|A - A_0\|_{L^{\sharp,2}_{1,A_0}} < \eps_1$
there is a gauge transformation $u\in \sG_E^3$ with the following
properties: 
\begin{itemize}
\item $d_{A_0}^*(u(A)-A_0) = 0$;
\item $\|u(A) - A_0\|_{L^2_{1,A_0}} 
\le cK_0\|A - A_0\|_{L^{\sharp,2}_{1,A_0}}$;
\item $\|u - \id_E\|_{L^{\sharp,2}_{2,A_0}}< 
cK_0\|A - A_0\|_{L^{\sharp,2}_{1,A_0}}$.
\end{itemize}
\end{thm}

\begin{pf}
The argument is broadly similar to that of Lemma \ref{lem:LocalDiffeo},
except that we can show $\bPsi$ is a diffeomorphism directly --- rather
than just a local diffeomorphism --- using the slightly stronger norms now
at our disposal. Moreover, on this occasion we seek precise bounds on the
solutions so we keep track of the dependence of constants on the curvature
$F_{A_0}$ and the least positive eigenvalue $\nu_0=\nu_0[A_0]$ of the
Laplacian $\De_{A_0}=d_{A_0}^*d_{A_0}$.

Write $A=A_0+a$ and observe that
$$
u(A)-A_0 = A-A_0 - (d_Au)u^{-1} = uau^{-1} - (d_{A_0}u)u^{-1}.
$$
Recall that we have an $L^2$-orthogonal decomposition 
$$
\Om^0(\fg_E) = (\Ker d_{A_0})^\perp\oplus\Ker d_{A_0} 
= \Imag d_{A_0}^*\oplus \Ker d_{A_0},
$$
and that $d_{A_0}^*:L^2_1(\La^1\otimes\fg_E)\to L^2(\fg_E)$ has closed
range; this gives 
\begin{align*}
L^{\sharp,2}_{2;A_0}(\fg_E) 
&= (\Ker d_{A_0}|_{L^{\sharp,2}_{2;A_0}})^\perp
\oplus \Ker d_{A_0}|_{L^{\sharp,2}_{2;A_0}} \\
&= (\Ker d_{A_0}^*|_{L^{\sharp,2}_{2;A_0}})^\perp
\oplus (\Imag d_{A_0}^*|_{L^{\sharp,2}_{1,A_0}}).
\end{align*}
We have a similar $L^2$-orthogonal decomposition 
$$
\Om^1(\fg_E) = \Imag d_{A_0}\oplus\Ker d_{A_0}^* 
= (\Ker d_{A_0}^*)^\perp\oplus\Ker d_{A_0}^*,
$$
and $d_{A_0}:L^2_1(\fg_E)\to L^2(\La^1\otimes\fg_E)$ has closed range; this
leads to the $L^2$-orthogonal decomposition
\begin{align*}
L^{\sharp,2}_{1,A_0}(\La^1\otimes\fg_E) 
&= (\Imag d_{A_0}|_{L^{\sharp,2}_{2;A_0}})
\oplus (\Ker d_{A_0}^*|_{L^{\sharp,2}_{1,A_0}}) \\
&= (\Ker d_{A_0}^*|_{L^{\sharp,2}_{2;A_0}})^\perp
\oplus (\Ker d_{A_0}^*|_{L^{\sharp,2}_{1,A_0}}).
\end{align*}
We now define a map 
\begin{align}
\bPsi: &(\Ker(d_{A_0}|_{L^{\sharp,2}_2}))^\perp\oplus
\Ker(d_{A_0}^*|_{L^{\sharp,2}_1})
\rightarrow L^{\sharp,2}_{1,A_0}(\La^1\otimes\fg_E), \label{eq:Psi}\\
&(\chi,a)  \mapsto uau^{-1} - (d_{A_0}u)u^{-1}, \notag
\end{align}
where $u=e^\chi$ and the differential at $(\chi,a)$ given by
\begin{align}
(D\bPsi)_{(\chi,a)}: &(\Ker(d_{A_0}|_{L^{\sharp,2}_2}))^\perp
\oplus \Ker (d_{A_0}^*|_{L^{\sharp,2}_1})
\to L^{\sharp,2}_1(\La^1\otimes\fg_E), \label{eq:DPsi}\\
&(\zeta,b) \mapsto u(-d_A\oplus \iota)u^{-1}(\zeta,b) = 
u(-d_A\zeta + b)u^{-1}, \notag
\end{align}
since $(D\bPsi)_{(0,a)}(\zeta,b)=-d_A\zeta + b$ and $\bPsi$ is
$\sG_E$-equivariant. Moreover, we have
\begin{equation}
(D^2\bPsi)_{(\chi,a)}((\zeta,b),(\eta,\al)) 
= u[\eta,-d_A\zeta+b]u^{-1} + u[\al,\zeta]u^{-1}. \label{eq:D2Psi}
\end{equation}
for $(\zeta,b),(\eta,\al)\in (\Ker(d_{A_0}|_{L^{\sharp,2}_2}))^\perp
\oplus \Ker (d_{A_0}^*|_{L^{\sharp,2}_1})$. 

We now verify that the conditions of the inverse function theorem
(Theorem \ref{thm:InverseFT}) hold for
suitable constants $K$ and $\de$. The operator 
$$
d_{A_0}:(\Ker(d_{A_0}|_{L^{\sharp,2}_2}))^\perp 
\to (\Ker(d_{A_0}^*|_{L^{\sharp,2}_1}))^\perp 
$$
has a two-sided inverse 
$$
G^0_{A_0}d^*_{A_0}: (\Ker(d_{A_0}^*|_{L^{\sharp,2}_1}))^\perp
\to (\Ker(d_{A_0}|_{L^{\sharp,2}_2}))^\perp.
$$
Indeed, for $b\in (\Ker(d_{A_0}^*|_{L^{\sharp,2}_1}))^\perp$, we have
$$
\|G^0_{A_0}d^*_{A_0}b\|_{L^{\sharp,2}_{2,A_0}} 
\le c_0K_0\|d^*_{A_0}b\|_{L^{\sharp,2}}
\le c_0K_0\|b\|_{L^{\sharp,2}_{1,A_0}},
$$
and so $G^0_{A_0}d^*_{A_0}$ has 
$\Hom(L^{\sharp,2}_{1,A_0},L^{\sharp,2}_{2,A_0})$ operator norm bound
$$
\|G^0_{A_0}d^*_{A_0}\| \le c_0K_0.
$$
In particular, we see that 
$(D\bPsi)_{(0,0)}^{-1} = G^0_{A_0}d^*_{A_0}\oplus \id$ satisfies
\begin{equation}
\|(D\bPsi)_{(0,0)}^{-1}\| \le c_0K_0
\label{eq:FirstDPsiCond}
\end{equation}
the first of the conditions we need to verify for $(D\bPsi)_{(0,0)}$ in
order to apply the inverse function theorem.

It remains to compare the differentials $(D\bPsi)_{(\chi,a)}$ and
$(D\bPsi)_{(0,0)}$ using the mean value theorem,
\begin{equation}
(D\bPsi)_{(\chi,a)}(\zeta,b) - (D\bPsi)_{(0,0)}(\zeta,b)
= \int_0^1(D^2\bPsi)_{(t\chi,ta)}((\zeta,b),(\chi,a))\,dt.
\label{eq:MVT}
\end{equation}
Thus, we need an estimate for $D^2\bPsi$:

\begin{claim}\label{claim:D2PsiEst}
There is a universal polynomial function $f(x,y)$, depending only on
$(X,g)$ and $G$, with $f(0,0)=0$, such that the following holds.
For any $t\in[0,1]$ we have:
\begin{align*}
&\|(D^2\bPsi)_{(t\chi,ta)}((\zeta,b),(\chi,a))\|_{L^{\sharp,2}_{1,A_0}} \\
&\qquad\le f(\|\chi\|_{L^{\sharp,2}_{2,A_0}},\|a\|_{L^{\sharp,2}_{1,A_0}})
\left(\|\zeta\|_{L^{\sharp,2}_{2,A_0}} 
+ \|b\|_{L^{\sharp,2}_{1,A_0}}\right).
\end{align*}
\end{claim}

\begin{pf}
From \eqref{eq:D2Psi} we have the $L^{\sharp,2}$ estimate 
\begin{align*}
&\|(D^2\bPsi)_{(t\chi,ta)}((\zeta,b),(\chi,a))\|_{L^{\sharp,2}} \\
&\qquad\le c\|\chi\|_{C^0}\left(\|d_{A_0}\zeta\|_{L^{\sharp,2}} 
+ \|a\|_{L^{\sharp,2}}\|\zeta\|_{C^0}
+ \|b\|_{L^{\sharp,2}}\right)
+ c\|a\|_{L^{\sharp,2}}\|\zeta\|_{C^0}, 
\end{align*}
and thus:
\begin{align}
&\|(D^2\bPsi)_{(t\chi,ta)}((\zeta,b),(\chi,a))\|_{L^{\sharp,2}} 
\label{eq:D2PsiL2Est}\\
&\le c\left(\|\chi\|_{L^{\sharp,2}_{2,A_0}}
+ \|a\|_{L^{\sharp,2}}\|\chi\|_{L^{\sharp,2}_{2,A_0}}
+ \|a\|_{L^{\sharp,2}}\right)
\left(\|\zeta\|_{L^{\sharp,2}_{2,A_0}} + 
+ \|b\|_{L^{\sharp,2}}\right). \notag
\end{align}
The $L^2$ estimate of $\cov_{A_0}(D^2\bPsi)_{(t\chi,ta)}((\zeta,b),(\chi,a))$
is given by
\begin{align*}
&\|\cov_{A_0}(D^2\bPsi)_{(t\chi,ta)}((\zeta,b),(\chi,a))\|_{L^2} \\
&\le c\left(\|\cov_{A_0}u\|_{L^4}\|\chi\|_{C^0} 
+ \|\cov_{A_0}\chi\|_{L^4}\right)
\left(\|d_{A_0}\zeta\|_{L^4}+\|a\|_{L^4}\|\zeta\|_{C^0} +
\|b\|_{L^4}\right)  \\
&\qquad + c\|\chi\|_{C^0}\left(\|\cov_{A_0}^2\zeta\|_{L^2}
+ \|\cov_{A_0}a\|_{L^2}\|\zeta\|_{C^0}
+ \|a\|_{L^4}\|\cov_{A_0}\zeta\|_{L^4} + \|\cov_{A_0}b\|_{L^2}\right) \\
&\qquad + c\|\cov_{A_0}u\|_{L^4}\|a\|_{L^4}\|\zeta\|_{L^4}
+ c\|\cov_{A_0}a\|_{L^2}\|\zeta\|_{C^0}
+ c\|a\|_{L^4}\|\cov_{A_0}\zeta\|_{L^4}, 
\end{align*}
and hence, using Lemma \ref{lem:SobolevExpEstimates} 
to estimate $u=e^\chi$ in terms of $\chi$,
\begin{align}
&\|\cov_{A_0}(D^2\bPsi)_{(t\chi,ta)}((\zeta,b),(\chi,a))\|_{L^2} 
\label{eq:D2PsiL2CovEst}\\
&\qquad\le f_1\left(\|\chi\|_{L^{\sharp,2}_{2,A_0}},
\|a\|_{L^2_{1,A_0}}\right)
\left(\|\zeta\|_{L^{\sharp,2}_{2,A_0}}  
+ \|b\|_{L^2_{1,A_0}}\right), \notag
\end{align}
where $f_1(x,y)$ is a polynomial function with $f_1(0,0)=0$. 

Noting that $d_{A_0}^*a=0$, we have
\begin{align}
d_{A_0}^*[a,\zeta] &= d_{A_0}^*(a\zeta-\zeta a) 
\label{eq:dA0*azeta}\\
&= (d_{A_0}^*a)\zeta - a\wedge d_{A_0}\zeta 
- *(d_{A_0}\zeta \wedge *a) - \zeta(d_{A_0}^*a) \notag\\
&= - a\wedge d_{A_0}\zeta - *(d_{A_0}\zeta \wedge *a). \notag
\end{align}
and similarly for $[\chi,b]$ since $d_{A_0}^*b=0$. For any $\be\in
L^2_1(\La^1\otimes\fg_E)$ we have
\begin{align}
d_{A_0}^*(u\be u^{-1}) &= -*d_{A_0}(u(*\be)u^{-1})) \label{eq:dA0*ubetau}\\
&= -*(d_{A_0}u\wedge *\be u^{-1}) + u(d^*_{A_0}\be)u^{-1} 
- *u((*\be)\wedge u(d_{A_0}u)u^{-1}). \notag
\end{align} 
Therefore, equations \eqref{eq:D2Psi}, \eqref{eq:dA0*azeta}, and
\eqref{eq:dA0*ubetau} and the estimates for $u=e^\chi$ in 
Lemma \ref{lem:SobolevExpEstimates} yield
\begin{align}
&\|d_{A_0}^*(D^2\bPsi)_{(\chi,a)}((\zeta,b),(\chi,a))\|_{L^{\sharp,2}} 
\label{eq:D2PsiLSharp2dA0*Est}\\
&\qquad\le \left\|d_{A_0}^*\left(u[\chi,-d_A\zeta+b]u^{-1} 
+ u[a,\zeta]u^{-1}\right)\right\|_{L^{\sharp,2}} \notag \\
&\qquad\le f_2(\|\chi\|_{L^{\sharp,2}_{2,A_0}},\|a\|_{L^{2\sharp,4}})
\left(\|\zeta\|_{L^{\sharp,2}_{2,A_0}} 
+ \|b\|_{L^{2\sharp,4}}\right), \notag
\end{align}
where $f_2(x,y)$ is a polynomial function with $f_2(0,0)=0$. The claim
now follows by combining 
\eqref{eq:D2PsiL2Est}, \eqref{eq:D2PsiL2CovEst}, and 
\eqref{eq:D2PsiLSharp2dA0*Est}. 
\end{pf}

Therefore, from Claim \ref{claim:D2PsiEst} and \eqref{eq:MVT} we have
\begin{align}
&\|(D\bPsi)_{(\chi,a)}(\zeta,b) 
- (D\bPsi)_{(0,0)}(\zeta,b)\|_{L^{\sharp,2}_{1,A_0}} 
\label{eq:CompareDPsiEst}\\
&\qquad\le f(\|\chi\|_{L^{\sharp,2}_{2,A_0}},\|a\|_{L^{\sharp,2}_{1,A_0}})
\left(\|\zeta\|_{L^{\sharp,2}_{2,A_0}}  
+ \|b\|_{L^{\sharp,2}_{1,A_0}}\right). \notag
\end{align}
Consequently, with respect to the $\Hom(L^{\sharp,2}_{2,A_0},
L^{\sharp,2}_{1,A_0})$ operator norm, \eqref{eq:CompareDPsiEst} yields
the bound
\begin{equation}
\|(D\bPsi)_{(\chi,a)}-(D\bPsi)_{(0,0)}\| \le \half c_0^{-1}K_0^{-1},
\label{eq:SecondDPsiCond}
\end{equation}
where $c_0K_0=K$ is the constant of \eqref{eq:FirstDPsiCond}, 
provided $(\chi,a)$ satisfies the constraint
\begin{equation}
\|\chi\|_{L^{\sharp,2}_{2,A_0}} + \|a\|_{L^{\sharp,2}_{1,A_0}}
\le c_1K_0^{-1}=\de.
\label{eq:ChiADistCond}
\end{equation} 
Define balls centered at the origins in
$(\Ker(d_{A_0}|_{L^{\sharp,2}_2}))^\perp$ and
$\Ker(d_{A_0}^*|_{L^{\sharp,2}_1})$ by setting
\begin{align*}
B^{\perp;2,\sharp,2}_0(\de) &= \{\chi\in 
(\Ker(d_{A_0}|_{L^{\sharp,2}_2}))^\perp:
\|\chi\|_{L^{\sharp,2}_{2,A_0}} < \de\}, \\
\bB^{1,\sharp,2}_0(\de) &= \{a\in \Ker(d_{A_0}^*|_{L^{\sharp,2}_1}):
\|a\|_{L^{\sharp,2}_{1,A_0}} < \de\}. 
\end{align*}
Hence, Theorem \ref{thm:InverseFT} implies that the map
$$
\bPsi: B^{\perp;2,\sharp,2}_0(\de) \times
\bB_0^{1,\sharp,2}(\de) 
\to \sA_E^{1,\sharp,2}
$$
is injective, its image is an open subset of $\sA_E^{1,\sharp,2}$ and
contains the ball $B_{A_0}^{1,\sharp,2}(\de/(2K))$, the inverse map
$\bPsi^{-1}$ is a diffeomorphism from $B_{A_0}^{1,\sharp,2}(\de/(2K))$ onto
its image, and if $(\chi_1,A_1)$, $(\chi_2,A_2)$ are points in 
$B^{\perp,2,\sharp,2}_0(\de)\times \bB^{1,\sharp,2}_0(\de)$, then
$$
\|\chi_1-\chi_2\|_{L^{\sharp,2}_{2,A_0}} + \|A_1-A_2\|_{L^{\sharp,2}_{1,A_0}}
\le 2K\|u_1(A_1)-u_2(A_2)\|_{L^{\sharp,2}_{1,A_0}},
$$
where $u_i=e^{\chi_i}$, $i=1,2$.
In particular, setting $(\chi_2,A_2-A_0)=(0,0)$,
we see that if $A$ is a point in $\sA^{1,\sharp,2}_E$
such that $\|A-A_0\|_{L^{\sharp,2}_{1,A_0}}<\de/(2K)$, then there is a unique
solution $(\chi,u^{-1}(A))=\bPsi^{-1}(A)$ in
$B^{\perp,2,\sharp,2}_0(\de)\times \bB^{1,\sharp,2}_0(\de)$.
Here, $u=e^\chi$ is a gauge transformation with
$\chi\in B^{\perp;2,\sharp,2}_0(\de)$ such that
\begin{align}
d^*_{A_0}(u^{-1}(A)-A_0) &= 0, \label{eq:ThmCoulombToDistanceEst}\\
\|\chi\|_{L^{\sharp,2}_{2,A_0}} 
+ \|u^{-1}(A)-A_0\|_{L^{\sharp,2}_{1,A_0}} 
&\le 2K\|A-A_0\|{L^{\sharp,2}_{1,A_0}}. \notag
\end{align}
Lemma \ref{lem:SobolevExpEstimates} implies that $u=e^\chi$ satisfies
\begin{equation}
\|u-\id_E\|_{L^{\sharp,2}_{2,A_0}} \le f_3(\|\chi\|_{L^{\sharp,2}_{2,A_0}})
\le c\|\chi\|_{L^{\sharp,2}_{2,A_0}} \le c_2\de,
\label{eq:uTochiEst}
\end{equation} 
where $f_3(x)$ is a polynomial with coefficients depending only on $(X,g)$
and $G$ such that $f_3(0)=0$. Noting that $K=c_0K_0$, $\de=c_1K_0^{-1}$, and
$\de/(2K) = \half c_0c_1K_0^{-2}$, the desired estimates follows from  
\eqref{eq:ThmCoulombToDistanceEst} and \eqref{eq:uTochiEst}.
Finally, Lemma \ref{lem:GaugeFixRegularity} implies that $u\in \sG_E^3$
and this completes the proof of the theorem.
\end{pf}

While the $L^2_1$ estimate of Theorem \ref{thm:CoulombConn} suffices for
most practical purposes, it is occasionally useful to have the slightly
weaker $L^{2\sharp,4}$ bound at hand. Recall from Section
\ref{sec:SharpSobolev} that we defined
$$
\|a\|_{\sL^{\sharp,2}_{1,A_0}}
= \|a\|_{L^{2\sharp,4}} + \|d_{A_0}^*a\|_{L^{\sharp,2}},
\qquad a\in\Om^1(\fg_E).
$$
A slight modification of the proof of Theorem \ref{thm:CoulombConn} yields:

\begin{thm}\label{thm:CoulombConnWeak}
Continue the hypotheses of Theorem \ref{thm:CoulombConn}.
Then for any $A\in\sA^2_E$ such that
$\|A - A_0\|_{\sL^{\sharp,2}_{1,A_0}} < \eps_1$
there is a gauge transformation $u\in \sG_E^3$ with the following
properties: 
\begin{itemize}
\item $d_{A_0}^*(u(A)-A_0) = 0$;
\item $\|u(A) - A_0\|_{L^{2\sharp,4}} 
\le cK_0\|A - A_0\|_{\sL^{\sharp,2}_{1,A_0}}$;
\item $\|u - \id_E\|_{L^{\sharp,2}_{2,A_0}}< 
 cK_0\|A - A_0\|_{\sL^{\sharp,2}_{1,A_0}}$.
\end{itemize}
\end{thm}

\begin{pf}
The first difference in the argument is that the map $\bPsi$ in
\eqref{eq:Psi} is replaced by
\begin{align}
\bPsi: &(\Ker(d_{A_0}|_{L^{\sharp,2}_2}))^\perp\oplus
\Ker(d_{A_0}^*|_{\sL^{\sharp,2}_1})
\rightarrow \sL^{\sharp,2}_{1,A_0}(\La^1\otimes\fg_E), \label{eq:PsiWeak}\\
&(\chi,a)  \mapsto uau^{-1} - (d_{A_0}u)u^{-1}. \notag
\end{align}
The second difference is that Claim \ref{claim:D2PsiEst} is replaced by:

\begin{claim}\label{claim:D2PsiWeakEst}
There is a universal polynomial function $f(x,y)$, depending only on
$(X,g)$ and $G$, with $f(0,0)=0$, such that the following holds.
For any $t\in[0,1]$ we have:
\begin{align*}
&\|(D^2\bPsi)_{(t\chi,ta)}((\zeta,b),(\chi,a))\|_{\sL^{\sharp,2}_{1,A_0}} \\
&\qquad\le f(\|\chi\|_{L^{\sharp,2}_{2,A_0}},\|a\|_{L^{2\sharp,4}})
\left(\|\zeta\|_{L^{\sharp,2}_{2,A_0}} 
+ \|b\|_{L^{2\sharp,4}}\right).
\end{align*}
\end{claim}

\begin{pf}
From \eqref{eq:D2Psi} we now have the $L^{2\sharp,4}$ estimate 
\begin{align*}
&\|(D^2\bPsi)_{(t\chi,ta)}((\zeta,b),(\chi,a))\|_{L^{2\sharp,4}} \\
&\qquad\le c\|\chi\|_{C^0}\left(\|d_{A_0}\zeta\|_{L^{2\sharp,4}} 
+ \|a\|_{L^{2\sharp,4}}\|\zeta\|_{C^0}
+ \|b\|_{L^{2\sharp,4}}\right)
+ c\|a\|_{L^{2\sharp,4}}\|\zeta\|_{C^0}, 
\end{align*}
and thus:
\begin{align}
&\|(D^2\bPsi)_{(t\chi,ta)}((\zeta,b),(\chi,a))\|_{L^{2\sharp,4}} 
\label{eq:D2PsiL2Sharp4Est}\\
&\le c\left(\|\chi\|_{L^{\sharp,2}_{2,A_0}}
+ \|a\|_{L^{2\sharp,4}}\|\chi\|_{L^{\sharp,2}_{2,A_0}}
+ \|a\|_{L^{2\sharp,4}}\right)
\left(\|\zeta\|_{L^{\sharp,2}_{2,A_0}} + 
+ \|b\|_{L^{2\sharp,4}}\right). \notag
\end{align}
Combining \eqref{eq:D2PsiLSharp2dA0*Est} and \eqref{eq:D2PsiL2Sharp4Est}
yields the claim. 
\end{pf}

The rest of the argument proceeds exactly as before and completes the
proof of the theorem.
\end{pf}

We now have our second proof of Theorem \ref{thm:GlobalUhlenbeck} via
Theorems \ref{thm:CoulombConn} and \ref{thm:CoulombConnWeak}: 

\begin{proof}[Proof of Theorem~\ref{thm:GlobalUhlenbeck}]
{}From the hypotheses we have $A_0\in\sA_E^k$ and $[A]\in\sB_E^k$ with
$k\ge 2$. According to Lemma \ref{lem:Metric}, there is gauge transformation
$w\in\sG_E^3$ such that 
$$
\dist_{L^{\sharp,2}_{1,A_0}}([A],[A_0]) 
= \|w(A)-A_0\|_{L^{\sharp,2}_{1,A_0}},
$$
where $A\in\sA_E^k$, so Theorems \ref{thm:CoulombConn} and the argument of
\ref{thm:CoulombConnWeak} imply that there is a gauge transformation
$v\in\sG_E^3$ so that $u(A)$ satisfies the conclusions of Assertion (2)
with $u=vw\in\sG_E^3$. Since $d_{A_0}^*(u(A)-A_0)=0$ and $u\in\sG_E^3$ and
$A,A_0\in\sA_E^k$, a standard bootstrapping argument implies that
$u\in\sG_E^{k+1}$.  
 
Similarly, by Lemma \ref{lem:Metric}, there is gauge transformation 
$w\in\sG_E^3$ such that 
$$
\dist_{\sL^{\sharp,2}_{1,A_0}}([A],[A_0]) 
= \|w(A)-A_0\|_{L^{\sharp,2}_{1,A_0}},
$$
so Assertion (1) follows from Theorem \ref{thm:CoulombConnWeak} in the same
manner. 
\end{proof}





\begin{thebibliography}{99}

\bibitem{Adams}
R. A. Adams, {\em Sobolev Spaces\/}, Academic Press, Orlando, FL, 1975. 

\bibitem{AHS} 
M. F. Atiyah, N. J. Hitchin, and I. M. Singer, {\em
Self-duality in four-dimensional Riemannian geometry\/}, Proc. Royal
Soc. London {\bf A 362} (1978) 425--461.

\bibitem{Aubin}
T. Aubin, {\em Nonlinear Analysis on Complex Manifolds. Monge-Amp\`ere
Equations\/}, Springer, New York, 1982.

\bibitem{Bradlow}
S. B. Bradlow, {\em Special metrics and stability for holomorphic bundles
with global sections\/}, J. Differential Geom. {\bf 33} (1991) 169--213.

\bibitem{BradlowDask}
S. B. Bradlow and G. D. Daskalopoulos, {\em Moduli of stable pairs for
holomorphic bundles over Riemann surfaces\/}, Internat. J. Math. {\bf
2} (1991) 477--513.

\bibitem{Daskal}
G. Daskalopoulos, private communication.

\bibitem{DonNS}
S. K. Donaldson, {\em A new proof of a theorem of Narasimhan and
Seshadri\/},  J. Differential Geom. {\bf 18} (1983) 269--277.

\bibitem{DonASD} 
\bysame, {\em Anti-self-dual Yang-Mills connections over complex
algebraic surfaces and stable vector bundles\/}, Proc. London
Math. Soc. {\bf 50} (1985) 1--26.

\bibitem{DonConn} 
\bysame, {\em Connections, cohomology and the
intersection forms of four-manifolds\/}, J. Differential Geom.
{\bf 24} (1986) 275--341.

\bibitem{DonHCobord}
\bysame, {\em Irrationality and the $h$-cobordism conjecture\/},
J. Differential Geom. {\bf 26} (1987), 141--168.

\bibitem{DonPoly} 
\bysame, {\em Polynomial invariants for
smooth 4-manifolds}, Topology {\bf 29}, 257-315 (1990).

\bibitem{DonApprox}
\bysame, {\em The approximation of instantons\/},
Geom. Funct. Anal. {\bf 3} (1993) 179--200.

\bibitem{DonSW}
\bysame, {\em The Seiberg-Witten equations and four-manifold
topology\/}, Bull. Amer. Math. Soc. {\bf 33} (1996) 45--70.

\bibitem{DK}
S. K. Donaldson and P. B. Kronheimer,  {\em The Geometry of Four-Manifolds\/},
Oxford University Press, Oxford, 1990.

\bibitem{DS} 
S. K. Donaldson and D. P. Sullivan, {\em Quasi-conformal
four-manifolds\/}, Acta Math. {\bf 163} (1990) 181--252.

\bibitem{FL1} 
P. M. N. Feehan and T. G. Leness, 
{\em {PU(2)} monopoles. {I}: {R}egularity, {U}hlenbeck compactness,
  and transversality}, J. Differential Geom. {\bf 49} (1998), 265--410,
  \texttt{dg-ga/9710032}.

\bibitem{FLGeorgia} 
\bysame, {\em {PU(2)} monopoles and relations between four-manifold
  invariants}, Topology Appl. {\bf 88} (1998), 111--145,
  \texttt{dg-ga/9709022}.

\bibitem{FL2} 
\bysame, {\em $\PU(2)$ monopoles, II: Highest-level singularities and 
relations between four-manifold invariants\/}, 
submitted to a print journal, December 8, 1997; \texttt{dg-ga/9712005}.

\bibitem{FL3} 
\bysame, {\em $\PU(2)$ monopoles, III:
Existence of gluing and obstruction maps\/}, in preparation. 

\bibitem{FL4} 
\bysame, {\em $\PU(2)$ monopoles, IV:
Surjectivity of gluing maps\/}, in preparation.

\bibitem{FeehanLeness}
\bysame, {\em in preparation}.

\bibitem{FU}
D. Freed and K. K. Uhlenbeck, {\em Instantons and Four-Manifolds\/}, 2nd
ed., Springer, New York, 1991.

\bibitem{FrM} 
R. Friedman and J. W. Morgan, {\em Smooth Four-manifolds and 
Complex Surfaces}, Springer, Berlin, 1994.

\bibitem{GT} 
D. Gilbarg and N. Trudinger, {\em Elliptic Partial Differential
Equations of Second Order\/}, 2nd ed., Springer, New York, 1983.

\bibitem{Goettsche}
L. G\"ottsche, {\em Modular forms and Donaldson invariants for 4-manifolds
with $b_+=1$\/}, J. Amer. Math. Soc. {\bf 9} (1996), 827--843;
\texttt{alg-geom/9506018}.  

\bibitem{KoM}
D. Kotschick and J. W. Morgan, {\em $\SO(3)$ invariants for four-manifolds
with $b^+=1$. II\/}, J. Differential Geom. {\bf 39} (1994) 433--456.

\bibitem{KMStructure}
P. B. Kronheimer and T. S. Mrowka, {\em Embedded surfaces and the
structure of Donaldson's polynomial invariants\/}, J. Differential Geom. 
{\bf 43} (1995) 573--734.

\bibitem{Lawson}
H. B. Lawson, {\em The Theory of Gauge Fields in Four Dimensions\/},
Conf. Board Math. Sci. {\bf 58}, Amer. Math. Soc., Providence, RI, 1985.

\bibitem{Li}
J. Li, {\em Algebraic geometric interpretation of Donaldson's polynomial
invariants\/}, J. Differential Geom. {\bf 37} (1993) 417--466. 

\bibitem{LM}
H. B. Lawson and M-L. Michelsohn, {\em Spin Geometry\/}, Princeton Univ.
Press, Princeton, NJ, 1988.

\bibitem{MitterViallet}
P. K. Mitter and C. M. Viallet, {\em On the bundle of connections and the
gauge orbit manifold in Yang-Mills theory\/}, Comm. Math. Phys. {\bf 79}
(1981) 457--472. 

\bibitem{MorganGTNotes}
J. W. Morgan, {\em Gauge Theory and the
Topology of Smooth Four-Manifolds\/}, Lecture Notes, Harvard University,
1988. 

\bibitem{MorganComparison}
\bysame, {\em Comparison of the Donaldson polynomial invariants with their
algebro-geometric analogues\/}, Topology {\bf 32} (1993) 449--488.

\bibitem{MorganOzsvath}
J. W. Morgan and P. S. Ozsv\'ath, private communication.

\bibitem{OTQuaternion}
C. Okonek and A. Teleman, {\em Quaternionic monopoles\/},
Comm. Math. Phys. {\bf 180} (1996) 363--388; \texttt{alg-geom/9505029}.

\bibitem{PalaisLS}
R. S. Palais, {\em Ljusternik-Snirelman theory on Banach manifolds\/},
Topology {\bf 5} (1966) 115--132.

\bibitem{Palais}
\bysame, {\em Foundations of Global Non-Linear Analysis\/},
Benjamin, New York, 1968.

\bibitem{PalaisCP}
\bysame, {\em Critical-point theory and the min-max principle\/},
Proc. Symp. Pure Math. {\bf 15}, Amer. Math. Soc., Providence, RI, 1970.

\bibitem{Parker}
T. H. Parker, {\em Gauge theories on four dimensional Riemannian manifolds\/},
Comm. Math. Phys. {\bf 85} (1982) 563--602.

\bibitem{ParkerTaubes}
T. H. Parker and C. H. Taubes, {\em On Witten's proof of the positive
energy theorem\/}, Comm. Math. Phys. {\bf 84} (1982) 223--238.

\bibitem{PTCambridge}
V. Y. Pidstrigach, Lecture at the Newton Institute, Cambridge, December 1994.

\bibitem{PTLocal}
V. Y. Pidstrigach and A. N. Tyurin, {\em Localization of Donaldson
invariants along the Seiberg-Witten classes\/}, \texttt{dg-ga/9507004}.

\bibitem{Singer}
I. M. Singer, {\em Some remarks on the Gribov ambiguity\/},
Comm. Math. Phys. {\bf 60} (1978) 7--12.

\bibitem{Stein} 
E. Stein, {\em Singular integral operators and differentiability properties
of functions\/}, Princeton University Press, Princeton, NJ, 1970.

\bibitem{TauSelfDual}
C. H. Taubes, {\em Self-dual Yang-Mills
connections on non-self-dual 4-manifolds\/}, J. Differential Geom.
{\bf 17} (1982) 139--170.

\bibitem{TauStab}
\bysame, {\em Stability in Yang-Mills theories\/},
Comm. Math. Phys. {\bf 91} (1983) 235--263. 

\bibitem{TauIndef} 
\bysame, {\em Self-dual connections on
4-manifolds with indefinite intersection matrix\/}, 
J. Differential Geom. {\bf 19} (1984) 517--560.

\bibitem{TauPath} 
\bysame, {\em Path-connected Yang-Mills moduli
spaces\/}, J. Differential Geom. {\bf 19} (1984) 337--392.

\bibitem{TauFrame} 
\bysame, {\em A framework for Morse theory for
the Yang-Mills functional}, Invent. Math. {\bf 94} (1988) 327--402. 

\bibitem{TauStable} 
\bysame, {\em The stable topology of
self-dual moduli spaces\/}, J. Differential Geom. {\bf 29} (1989) 162--230.

\bibitem{TauConf} 
\bysame, {\em The existence of anti-self-dual
conformal structures}, J. Differential Geom. {\bf 36} (1992) 163-253.

\bibitem{TauGluing} 
\bysame, {\em Metrics, Connections, and Gluing Theorems\/},
Conf. Board Math. Sci. {\bf 89}, Amer. Math. Soc., Providence, RI, 1996.

\bibitem{TelemanNonabelian}
A. Teleman, {\em Non-abelian Seiberg-Witten theory and stable oriented
pairs\/}, Internat. J. Math. {\bf 8} (1997) 507--535;
\texttt{alg-geom/9609020}. 

\bibitem{UhlLp}
K. K. Uhlenbeck, {\em Connections with $L^p$ bounds on curvature\/},
Comm. Math. Phys. {\bf 83} (1982) 31--42.

\bibitem{Witten}
E. Witten, {\em Monopoles and four-manifolds}\/, Math. Res. Lett. {\bf
1} (1994) 769--796.

\end{thebibliography}
\end{document}